\documentclass{svjour3}                   
\smartqed  

\usepackage{mathptmx}      

\usepackage[T1]{fontenc}
\usepackage[utf8]{inputenc}
\usepackage{babel}

\usepackage{amsmath}
\usepackage{amssymb}
\usepackage{graphicx}
\usepackage{color}
\usepackage{bm}
\usepackage{fancybox}		
\usepackage{array}
\usepackage{booktabs}		
\usepackage[figuresright]{rotfloat}	

\usepackage{stmaryrd}
\usepackage{stackrel}

\usepackage{flushend}

\usepackage[numbers,sort]{natbib}

\providecommand{\tabularnewline}{\\}

\definecolor{note_fontcolor}{rgb}{1, 0.667969, 0}

\usepackage{environ}
\usepackage[colorinlistoftodos,prependcaption,textsize=small]{todonotes}
\RenewEnviron{lyxgreyedout}{\todo[inline]{\protect\BODY}}

\usepackage{etoolbox}

\spnewtheorem*{defn*}{Definition}{\bfseries}{\rmfamily}

\renewcommand{\vec}[1]{\bm{#1}}

\usepackage{hyperref}

\makeatother

\newcommand{\dd}{\mathrm{d}}

\newcommand{\R}{\mathbb{R}}

\newcommand{\Ls}{\mathrm{L}}
\newcommand{\Cs}{\mathrm{C}}

\journalname{arXiv}
\begin{document}

\title{Numerical Optimization of the Dirichlet Boundary Condition in the Phase Field Model with an Application to Pure Substance Solidification}

\titlerunning{Optimization of the Dirichlet B.C. in the Phase Field Model of Pure Substance Solidification}

\author{Ale\v{s} Wodecki		\and
	Pavel Strachota         	\and
        Tom\'{a}\v{s} Oberhuber		\and
	Kate\v{r}ina \v{S}kardov\'{a}	\and
	Monika Bal\'{a}zsov\'{a}
}

\authorrunning{Wodecki, Strachota, Oberhuber, \v{S}kardov\'{a}, Bal\'{a}zsov\'{a}} 

\institute{A. Wodecki \at
              Department of Mathematics, Faculty of Nuclear Sciences and Physical Engineering,
              Czech Technical University in Prague \\
              Tel.: +420-733-522-679\\
              \email{wodecale@fjfi.cvut.cz}           
           \and
	   P. Strachota \at
              Department of Mathematics, Faculty of Nuclear Sciences and Physical Engineering,
              Czech Technical University in Prague \\
              Tel.: +420-22435-8563\\
              \email{pavel.strachota@fjfi.cvut.cz}           
           \and
           T. Oberhuber, K. \v{S}kardov\'{a}, M. Bal\'{a}zsov\'{a} \at
              Department of Mathematics, Faculty of Nuclear Sciences and Physical Engineering,
              Czech Technical University in Prague \\
              \email{tomas.oberhuber@fjfi.cvut.cz, katerina.skardova@fjfi.cvut.cz }           
}

\date{Received: date / Accepted: date}

\maketitle

\begin{abstract}
As opposed to the distributed control of parabolic PDE's, very few
contributions currently exist pertaining to the Dirichlet boundary
condition control for parabolic PDE's. This motivates our interest
in the Dirichlet boundary condition control for the phase field model
describing the solidification of a pure substance from a supercooled
melt. In particular, our aim is to control the time evolution of the
temperature field on the boundary of the computational domain in order
to achieve the prescribed shape of the crystal at the given time.
To obtain efficient means of computing the gradient of the cost functional,
we derive the adjoint problem formally. The gradient is then used
to perform gradient descent. The viability of the proposed optimization
method is supported by several numerical experiments performed in
one and two spatial dimensions. Among other things, these experiments
show that a linear reaction term in the phase field equation proves
to be insufficient in certain scenarios and so an alternative reaction
term is considered to improve the models behavior.


\keywords{
adjoint method \and Allen-Cahn equation \and crystal growth \and optimization \and functional analytic framework
}


\subclass{35K51 \and 35K57 \and 49M41 \and 35Q93 \and 65K10}
\end{abstract}

\newlength \figwidth
\setlength \figwidth {0.6\columnwidth}



\section{Introduction}

Phase field models have found utility in various areas including phase
transition modeling \cite{Benes-Anisotropic_phase_field_model,PF-memory,Physics-first-order-phase,Karma_Rappel-Quant_phase_field_modeling},
multi component flow simulations \cite{Alpak-Phase-Field-TwoPhase-Porous,Amiri-Level_set-vs-Phase_Field-in-TwoPhase_Flow}
and fracture mechanics \cite{WU20201,PF-fracture-crack-propagation}.
They can take the form of an initial boundary-value problem for a
system of partial differential equations (PDE's). In this case, one
of the PDE's is the phase field equation, which is derived by the
minimization of the Allen-Cahn or Cahn-Hilliard functional \cite{Allen-Cahn-orig,free-energy-cahn-hil}.
To describe the solidification of a pure supercooled liquid \cite{Benes-Math_comp_aspects_solid,PF-Focusing-Latent-Heat}
the phase equation is augmented by the heat equation. This formulation
is linked to materials science and practical applications, where controlling
solidification is of interest \cite{Guo_Xiong-3D_multi_dendritic_growth_coarsening,ISPMA14-Pavel_Ales-ActaPhPoloA,Jeong_Goldenfeld_Dantzig-PF-FEM-3D-Flow,PF-cystal-growth-2009}.

The optimal control phase field models (PFM) has been studied extensively
both from a theoretical and numerical simulations perspective. In
particular, the distributed control of this type of PFM with homogeneous
Neumann or Dirichlet boundary conditions has been addressed by a large
number of publica\-tions \cite{Colli2019OptimalCO,2163-2480_2018_1_95,PF-Hoff-92,Zonghong22}.
Another type of control that may be considered for PFM is the Neumann
or Robin Control (NoR) \cite{CHRYSAFINOS2006891,SenerNeumann}. To
the authors' best knowledge however, the Dirichlet boundary control
of parabolic equations has not been studied as extensively. Hinze
et. al. \cite{Hinze-control-very-weak-sol-16} and Kunisch \cite{Kun-parabolic-opt-cont-07}
have utilized the very weak formulation to derive optimality conditions
and perform calculations using the finite element method (FEM) \cite{Hinze-control-very-weak-sol-16,Kun-parabolic-opt-cont-07}.
More recently, Gudi et. al. have utilized the weak formulation along
with an alternative control space formulation to derive theoretical
results and perform FEM simulations as well \cite{Gudi_parabolic_control_21}.
All these contributions however, only address a single parabolic PDE
and do not apply to the optimization of the Dirichlet boundary condition
in the PFM.

In this article, the numerical optimization of the Dirichlet boundary
condition in a PFM that describes the solidification of a supercooled
pure substance is discussed. Since we restrict ourselves to numerical
simulations, the strong form can be used to formally derive the adjoint
equations. These are then used to propose an efficient method of gradient
computation. Gradient descent is then performed using an initial guess
to eventually arrive at a locally optimal control. A number of numerical
simulations in one and two dimensions have been performed to show
the effectiveness of this method (see Section \ref{sec:Numerical-Results-1}).
Some of these simulations show that the optimization is quantitatively
and qualitatively affected by the inadequacies of the linear reaction
term in the phase field equation (see Section \ref{subsec:Moving-a-Gap}).
An alternative reaction term, proposed in \cite{PF-Focusing-Latent-Heat},
is then used to remove these limitations.

\section{General Problem Formulation and Gradient Computation\label{sec:General-Problem-Formulation}}

For the general theory, we formally follow \cite{Opt-Book-Hind}.
Let $Y,U$ be Banach spaces and $Z$ be a Hilbert space. We call $Y$
the solution space and $U$ the control space.  Consider a map

\[
e:Y\times U\rightarrow Z.
\]
We call $e$ the state map, implicitly defining the dependence of
the solution $y\in Y$ on the control $u\in U$ by the state equation

\begin{equation}
e\left(y,u\right)=0.\label{eq:state-equation}
\end{equation}
For example, the state equation can assume the form of a system of
ordinary differential equations (ODE's) with boundary conditions or
a system of PDE's with boundary and initial conditions. In this setting,
$u$ is involved in the formulation of the problem for the ODE's or
PDE's and $y$ represents its solution. In particular, $u$ can appear
in the source term or an initial or boundary condition.

Assume that there exists a map $S:U\rightarrow Y$ such that
\[
e\left(S\left(u\right),u\right)=0.
\]
We call $S$ the solution operator. The existence of the map $S$
is equivalent to the existence and uniqueness of the solution $y=S\left(u\right)$
of (\ref{eq:state-equation}) for any control $u\in U$.

We introduce the cost functional $J:Y\times U\rightarrow\mathbb{R}$.
The fundamental minimization problem then reads

\begin{align}
\min_{u\in U} & J\left(y,u\right)\label{eq:FP - functional to min}\\
\text{s.t. } & e\left(y,u\right)=0\text{ where }y\in Y,u\in U.\label{eq:FP - state equation}
\end{align}

To obtain the derivative of the cost functional $J$ with respect
to $u$, it is possible to either perform direct computation (sensitivity
analysis) or use adjoint methods \cite{Opt-Book-Hind}. We opt for
the latter approach since it's computationally more efficient when
performing optimization with respect to a large amount of parameters.

First, define the reduced cost functional associated to the fundamental
problem (\ref{eq:FP - functional to min})-(\ref{eq:FP - state equation})
as

\[
\hat{J}\left(u\right)=J\left(y\left(u\right),u\right),
\]
where the term $y\left(u\right)$ is used in place of $S\left(u\right)$
to denote the solution $y$ of the state equation (\ref{eq:state-equation})
for the given control $u$. Then the fundamental problem (\ref{eq:FP - functional to min})-(\ref{eq:FP - state equation})
simply becomes
\begin{align}
\underset{u\in U}{\text{min}}\hat{J}\left(u\right) & .\label{eq:FP - functional to min-RED}
\end{align}
We compute the Fréchet derivative of $\hat{J}$ at $u$ formally.
Define the Lagrange function $L:Y\times U\times Z\rightarrow\mathbb{R}$
as

\begin{equation}
L\left(y,u,\lambda\right)\equiv J\left(y,u\right)+\left\langle \lambda,e\left(y,u\right)\right\rangle _{Z},\label{eq:Lag theory - general Lagrangian}
\end{equation}
where $\lambda$ is called the adjoint variable. We notice that since
$e\left(y\left(u\right),u\right)=0$ for all $u\in U$, we have

\begin{equation}
L\left(y\left(u\right),u,\lambda\right)=\hat{J}\left(u\right).\label{eq:lagrangian and reduced functional relationship}
\end{equation}
By differentiating (\ref{eq:lagrangian and reduced functional relationship}),
we arrive at
\begin{align}
\hat{J}^{\prime}\left(u\right)s & =L_{y}\left(y\left(u\right),u,\lambda\right)y^{\prime}\left(u\right)s+L_{u}\left(y\left(u\right),u,\lambda\right)s,\label{eq:derivative calculation}
\end{align}
where $s\in U$. The equation
\begin{equation}
L_{y}\left(y\left(u\right),u,\lambda\right)=0\label{eq:lagrangian adjoint condition}
\end{equation}
is called the adjoint problem. Finding a particular $\lambda_{0}\in Z$
that solves (\ref{eq:lagrangian adjoint condition}) leads to a simplified
calculation of the derivative (\ref{eq:derivative calculation}) and
we no longer need to compute $y^{\prime}\left(u\right)$ in order
to evaluate $\hat{J}^{\prime}\left(u\right)s$ for a given direction
$s\in U$.

\section{\label{subsec:Adjoint-Method-Derivative-PF}Adjoint Problem for Optimizing
the Solution of the Phase Field Problem}

Using strong formalism, the procedure shown in Section \ref{sec:General-Problem-Formulation}
is applied to an optimization problem where the state equation (\ref{eq:FP - state equation})
is the isotropic phase field problem in a simplified dimensionless
form \cite{PF-Focusing-Latent-Heat,Benes-Math_comp_aspects_solid,Benes-Asymptotics,CaginalpAnalysis86}.
Let $\Omega\subset\mathbb{R}^{n}$ be a bounded domain and let $T>0$.
The course of solidification (or melting) of a pure material in $\Omega$
is described by the evolution of the phase field $\tilde{y}$ with
values between $0$ and $1$, identifying the solid subdomain $\Omega_{s}\left(t\right)$,
the liquid subdomain $\Omega_{l}\left(t\right)$ and the phase interface
$\Gamma\left(t\right)$ by
\begin{align}
\Omega_{s}\left(t\right) & =\left\{ x\in\Omega;\tilde{y}\left(t,x\right)>\frac{1}{2}\right\} ,\label{eq:solid-subdomain}\\
\Omega_{l}\left(t\right) & =\left\{ x\in\Omega;\tilde{y}\left(t,x\right)<\frac{1}{2}\right\} ,\label{eq:liquid-subdomain}\\
\Gamma\left(t\right) & =\left\{ x\in\Omega;\tilde{y}\left(t,x\right)=\frac{1}{2}\right\} .\label{eq:phase-interface}
\end{align}
The phase field problem governs the evolution of both the phase field
$\tilde{y}$ and the temperature field $y$ by modeling phase transitions
and heat transfer inside $\Omega$. Our aim is to obtain such solution
that approaches the prescribed crystal shape at the given time $T$
(i.e., the shape of the solid subdomain $\Omega_{s}\left(T\right)$).
This is done by controlling the time-dependent Dirichlet boundary
condition for the temperature $y$.

Consider the problem
\begin{equation}
\min_{u\in U}J\left(y,\tilde{y},u\right)\equiv\frac{1}{2}\underset{\Omega}{\int}\left|\tilde{y}\left(T,x\right)-\tilde{y}_{f}\left(x\right)\right|^{2}\dd x+\frac{\alpha}{2}\stackrel[0]{T}{\int}\underset{\partial\Omega}{\int}\left|u\left(t,x\right)\right|^{2}\dd S\dd t\label{eq:PF-functional}
\end{equation}
s.t.
\begin{align}
y_{t}= & \Delta y+H\tilde{y}_{t}, & \text{in }\left(0,T\right)\times\Omega,\label{eq:PF-heat-equation}\\
y\left|_{\partial\Omega}\right.= & u & \text{on }\left[0,T\right]\times\partial\Omega,\label{eq:PF-heat-BC}\\
y\left|_{t=0}\right.= & y_{\text{ini}} & \text{in }\Omega,\label{eq:PF-heat-IC}\\
\gamma\xi^{2}\tilde{y}_{t}= & \xi^{2}\Delta\tilde{y}+f_{0}\left(y,\tilde{y};\xi\right) & \text{in }\left(0,T\right)\times\Omega,\label{eq:PF-Allen-Cahn-equation}\\
\tilde{y}\left|_{\partial\Omega}\right.= & \tilde{y}_{\text{bc}} & \text{on }\left[0,T\right]\times\partial\Omega,\label{eq:PF-PF-BC}\\
\tilde{y}\left|_{t=0}\right.= & \tilde{y}_{\text{ini}} & \text{in }\Omega,\label{eq:PF-PF-IC}
\end{align}
where $\tilde{y}_{f}\in\Ls^{2}\left(\Omega\right)$ in (\ref{eq:PF-functional})
is the target profile of the phase field $\tilde{y}$, $\alpha$ denotes
the strength of the regularization term, and
\begin{equation}
f_{0}\left(y,\tilde{y};\xi\right)=\tilde{y}\left(1-\tilde{y}\right)\left(\tilde{y}-\frac{1}{2}\right)-\beta\xi\left(y-y_{\text{mt}}\right).\label{eq:reaction-term}
\end{equation}
Equation (\ref{eq:PF-heat-equation}) is the heat equation with release
of the latent heat of fusion $H$. The corresponding Dirichlet boundary
condition for the temperature $y$ is given by the control $u\in\Cs\left(\left[0,T\right]\times\partial\Omega\right)$
and the initial temperature field is determined by (\ref{eq:PF-heat-IC}).

Next, (\ref{eq:PF-Allen-Cahn-equation}) is the Allen-Cahn (phase
field) equation with a simple linear form of the reaction term \cite{Kobayashi-PF-Dendritic,PF-Focusing-Latent-Heat},
containing the melting temperature $y_{\text{mt}}$, the parameter
$\xi$ related to the interface thickness, and the dimensionless model
parameters $\gamma,\beta$. For simplicity, we consider a constant
Dirichlet boundary condition (\ref{eq:PF-PF-BC}) for $\tilde{y}$.
Finally, the initial condition (\ref{eq:PF-PF-IC}) describes the
initial shape of the crystal $\Omega_{s}\left(0\right)$.

In order for the phase interface to form correctly and the solution
to have physical interpretation, the reaction term on the right hand
side of (\ref{eq:PF-Allen-Cahn-equation}) must have three roots in
terms of $\tilde{y}$ and thus satisfy the condition

\begin{equation}
\beta\xi\left(y-y_{\text{mt}}\right)\in\left(-\frac{\sqrt{3}}{36},\frac{\sqrt{3}}{36}\right),\label{eq:PF physicality condition}
\end{equation}
as detailed in \cite{PF-Focusing-Latent-Heat,Benes-Math_comp_aspects_solid,CaginalpAnalysis86}.
As will be shown further in Section \ref{sec:Numerical-Results-1},
the presented formulation of the optimization problem may yield the
optimal control $u$ such that the solution $\left(y,\tilde{y}\right)$
violates (\ref{eq:PF physicality condition}).

Let us proceed to the derivation of the adjoint problem. Consider
the setting:
\begin{itemize}
\item $Y\equiv\Cs^{2}\left(\left[0,T\right]\times\overline{\Omega}\right)^{2},$
\item $U\equiv\Cs\left(\partial\Omega\times\left[0,T\right]\right),$
\item $Z\equiv\Ls^{2}\left(\left[0,T\right]\times\overline{\Omega}\right)^{2}\times\Ls^{2}\left(\partial\Omega\times\left[0,T\right]\right)^{2}\times\Ls^{2}\left(\Omega\right)^{2}$.
\end{itemize}
Define the state equation operator component-wise as

\begin{align}
e^{1}\left(y,\tilde{y}\right)= & y_{t}-\Delta y-H\tilde{y}_{t},\nonumber \\
e^{2}\left(y,u\right)= & y\left|_{t=0}\right.-y_{\text{ini}},\nonumber \\
e^{3}\left(y\right)= & y\left|_{\partial\Omega}-u\right.,\nonumber \\
e^{4}\left(y,\tilde{y}\right)= & \gamma\xi^{2}\tilde{y}_{t}-\xi^{2}\Delta\tilde{y}-\tilde{y}\left(1-\tilde{y}\right)\left(\tilde{y}-\frac{1}{2}\right)+\beta\xi\left(y-y_{\text{mt}}\right),\nonumber \\
e^{5}\left(\tilde{y}\right)= & \tilde{y}\left|_{t=0}\right.-\tilde{y}_{\text{ini}},\nonumber \\
e^{6}\left(\tilde{y}\right)= & \tilde{y}\left|_{\partial\Omega}\right.-\tilde{y}_{\text{bc}}.\label{eq: PF strong defintion of operators}
\end{align}
Let $\lambda\equiv\left(p_{1},p_{2},p_{3},q_{1},q_{2},q_{3}\right)\in Z$,
then the Lagrangian for the problem (\ref{eq:PF-functional})-(\ref{eq:PF-PF-IC})
reads

\begin{align}
L\left(y,\tilde{y},u,\lambda\right)= & J\left(y,\tilde{y},u\right)+\lambda\left(e\left(y,\tilde{y},u\right)\right)\nonumber \\
= & \frac{1}{2}\underset{\Omega}{\int}\left|\tilde{y}\left|_{t=T}\right.-\tilde{y}_{f}\right|^{2}\dd x+\frac{\alpha}{2}\stackrel[0]{T}{\int}\underset{\partial\Omega}{\int}\left|u\right|^{2}\dd S\dd t\nonumber \\
+ & \stackrel[0]{T}{\int}\underset{\Omega}{\int}\left(y_{t}-\Delta y-H\tilde{y}_{t}\right)p_{1}\dd x\dd t+\underset{\Omega}{\int}\left(y\left|_{t=0}\right.-y_{\text{ini}}\right)p_{2}\dd x+\stackrel[0]{T}{\int}\underset{\partial\Omega}{\int}\left(y\left|_{\partial\Omega}\right.-u\right)p_{3}\dd S\dd t\nonumber \\
+ & \stackrel[0]{T}{\int}\underset{\Omega}{\int}\left(\gamma\xi^{2}\tilde{y}_{t}-\xi^{2}\Delta\tilde{y}-\tilde{y}\left(1-\tilde{y}\right)\left(\tilde{y}-\frac{1}{2}\right)+\beta\xi\left(y-y_{\text{mt}}\right)\right)q_{1}\dd x\dd t\nonumber \\
+ & \underset{\Omega}{\int}\left(\tilde{y}\left|_{t=0}\right.-\tilde{y}_{\text{ini}}\right)q_{2}\dd x+\stackrel[0]{T}{\int}\underset{\partial\Omega}{\int}\left(\tilde{y}\left|_{\partial\Omega}\right.-\tilde{y}_{\text{bc}}\right)q_{3}\dd S\dd t.\label{eq:PF Lagrangian}
\end{align}
Consider $\left(v,\hat{v}\right)\in Y$, then
\begin{align}
L_{\left(y,\tilde{y}\right)}\left(y,\tilde{y},u,\lambda\right)\left[\left(v,\hat{v}\right)\right]= & \underset{\Omega}{\int}\left(\tilde{y}\left|_{t=T}\right.-\tilde{y}_{f}\right)\hat{v}\left|_{t=T}\right.\dd x+\overset{I.}{\overbrace{\stackrel[0]{T}{\int}\underset{\Omega}{\int}\left(v_{t}-\Delta v-H\hat{v}_{t}\right)p_{1}\dd x\dd t}}\nonumber \\
+ & \underset{\Omega}{\int}v\left|_{t=0}\right.p_{2}\dd x+\stackrel[0]{T}{\int}\underset{\partial\Omega}{\int}v\left|_{\partial\Omega}\right.p_{3}\dd S\dd t\nonumber \\
+ & \overset{II.}{\overbrace{\stackrel[0]{T}{\int}\underset{\Omega}{\int}\left(\gamma\xi^{2}\hat{v}_{t}-\xi^{2}\Delta\hat{v}+3\tilde{y}^{2}\hat{v}-3\tilde{y}\hat{v}+\frac{1}{2}\hat{v}+\beta\xi v\right)q_{1}\dd x\dd t}}\nonumber \\
+ & \underset{\Omega}{\int}\hat{v}\left|_{t=0}\right.q_{2}\dd x+\stackrel[0]{T}{\int}\underset{\partial\Omega}{\int}\hat{v}\left|_{\partial\Omega}\right.q_{3}\dd S\dd t.\label{eq:PF - Frechet Derivative}
\end{align}
Condition (\ref{eq:lagrangian adjoint condition}) is satisfied if
and only if

\[
L_{\left(y,\tilde{y}\right)}\left(y,\tilde{y},u,\lambda\right)\left[\left(v,\hat{v}\right)\right]=0\text{ for all }\left(v,\hat{v}\right)\in Y.
\]
To this end, we use Greens formula to offload the derivatives in expressions
$I.$ and $II.$ Expression $I.$ then becomes

\begin{align}
I.= & -\stackrel[0]{T}{\int}\underset{\Omega}{\int}\left(p_{1}\right)_{t}v\dd x\dd t+\underset{\Omega}{\int}\left(p_{1}v\right)\left|_{t=T}\right.-\left(p_{1}v\right)\left|_{t=0}\right.\dd x+\stackrel[0]{T}{\int}\underset{\Omega}{\int}\nabla p_{1}\cdot\nabla v\dd x\dd t\nonumber \\
- & \stackrel[0]{T}{\int}\underset{\partial\Omega}{\int}p_{1}\nabla v\cdot\vec{n}\dd S\dd t+\stackrel[0]{T}{\int}\underset{\Omega}{\int}H\left(p_{1}\right)_{t}\hat{v}\dd x\dd t+\underset{\Omega}{\int}-H\left(p_{1}\hat{v}\right)\left|_{t=T}\right.+H\left(p_{1}\hat{v}\right)\left|_{t=0}\right.\dd x\nonumber \\
= & -\stackrel[0]{T}{\int}\underset{\Omega}{\int}\left(p_{1}\right)_{t}v\dd x\dd t+\underset{\Omega}{\int}\left(p_{1}v\right)\left|_{t=T}\right.-\left(p_{1}v\right)\left|_{t=0}\right.\dd x-\stackrel[0]{T}{\int}\underset{\Omega}{\int}\Delta p_{1}v\dd x\dd t+\stackrel[0]{T}{\int}\underset{\partial\Omega}{\int}v\nabla p_{1}\cdot\vec{n}\dd S\dd t\nonumber \\
- & \stackrel[0]{T}{\int}\underset{\partial\Omega}{\int}p_{1}\nabla v\cdot\vec{n}\dd S\dd t+\stackrel[0]{T}{\int}\underset{\Omega}{\int}H\left(p_{1}\right)_{t}\hat{v}\dd x\dd t+\underset{\Omega}{\int}-H\left(p_{1}\hat{v}\right)\left|_{t=T}\right.+H\left(p_{1}\hat{v}\right)\left|_{t=0}\right.\dd x.\label{eq:PF - I.}
\end{align}
Analogously, expression $II.$ can be rewritten as

\begin{align}
II.= & -\stackrel[0]{T}{\int}\underset{\Omega}{\int}\gamma\xi^{2}\left(q_{1}\right)_{t}\hat{v}\dd x\dd t+\underset{\Omega}{\int}\gamma\xi^{2}\left(q_{1}\hat{v}\right)\left|_{t=T}\right.-\gamma\xi^{2}\left(q_{1}\hat{v}\right)\left|_{t=0}\right.\dd x\nonumber \\
+ & \xi^{2}\left[\stackrel[0]{T}{\int}\underset{\Omega}{\int}-\Delta q_{1}\hat{v}\dd x\dd t-\stackrel[0]{T}{\int}\underset{\partial\Omega}{\int}q_{1}\nabla\hat{v}\cdot\vec{n}\dd S\dd t+\stackrel[0]{T}{\int}\underset{\partial\Omega}{\int}\hat{v}\nabla q_{1}\cdot\vec{n}\dd S\dd t\right]\nonumber \\
+ & \stackrel[0]{T}{\int}\underset{\Omega}{\int}\left(3\tilde{y}^{2}q_{1}-3\tilde{y}q_{1}+\frac{1}{2}q_{1}\right)\hat{v}\dd x\dd t+\stackrel[0]{T}{\int}\underset{\Omega}{\int}\beta\xi q_{1}v\dd x\dd t.\label{eq:PF - II.}
\end{align}
Using (\ref{eq:PF - I.}), (\ref{eq:PF - II.}) along with (\ref{eq:PF - Frechet Derivative})
results in

\begin{align}
L_{\left(y,\tilde{y}\right)}\left(y,\tilde{y},u,\lambda\right)\left[\left(v,\hat{v}\right)\right]= & \underset{\Omega}{\int}\left(\tilde{y}\left|_{t=T}\right.-\tilde{y}_{f}-Hp_{1}\left|_{t=T}\right.+\gamma\xi^{2}q_{1}\left|_{t=T}\right.\right)\hat{v}\left|_{t=T}\right.\dd x\nonumber \\
+ & \stackrel[0]{T}{\int}\underset{\Omega}{\int}\left(-\left(p_{1}\right)_{t}-\Delta p_{1}+\beta\xi q_{1}\right)v\dd x\dd t+\underset{\Omega}{\int}\left(p_{2}-p_{1}\right)\left|_{t=0}\right.v\left|_{t=0}\right.\dd x\nonumber \\
+ & \stackrel[0]{T}{\int}\underset{\partial\Omega}{\int}v\left(\nabla p_{1}\cdot\vec{n}+p_{3}\right)\dd S\dd t-\stackrel[0]{T}{\int}\underset{\partial\Omega}{\int}p_{1}\nabla v\cdot\vec{n}\dd S\dd t+\underset{\Omega}{\int}\left(p_{1}v\right)\left|_{t=T}\right.\dd x\nonumber \\
+ & \stackrel[0]{T}{\int}\underset{\Omega}{\int}\left(-\gamma\xi^{2}\left(q_{1}\right)_{t}-\xi^{2}\Delta q_{1}+H\left(p_{1}\right)_{t}+3\tilde{y}^{2}q_{1}-3\tilde{y}q_{1}+\frac{1}{2}q_{1}\right)\hat{v}\dd x\dd t\nonumber \\
+ & \underset{\Omega}{\int}\left(Hp_{1}-\gamma\xi^{2}q_{1}+q_{2}\right)\hat{v}\left|_{t=0}\right.\dd x-\stackrel[0]{T}{\int}\underset{\partial\Omega}{\int}\xi^{2}q_{1}\nabla\hat{v}\cdot\vec{n}\dd S\dd t\nonumber \\
+ & \stackrel[0]{T}{\int}\underset{\partial\Omega}{\int}\left(q_{3}+\xi^{2}\nabla q_{1}\cdot\vec{n}\right)\hat{v}\left|_{\partial\Omega}\right.\dd S\dd t.\label{eq:PF - final Frecher der}
\end{align}
From (\ref{eq:PF - final Frecher der}), we see that by providing
a $p_{1}$ that solves

\begin{align}
\left(p_{1}\right)_{t}+\Delta p_{1}= & \beta\xi q_{1} & \text{in }\left(0,T\right)\times\Omega,\nonumber \\
p_{1}\left|_{\partial\Omega}\right.= & 0 & \text{on }\partial\Omega\times\left[0,T\right],\nonumber \\
p_{1}\left|_{t=T}\right.= & 0 & \text{in }\Omega\label{eq:PF adjoint part I.}
\end{align}
and a $q_{1}$ that solves

\begin{align}
\gamma\xi^{2}\left(q_{1}\right)_{t}+\xi^{2}\Delta q_{1}= & H\left(p_{1}\right)_{t}+3\tilde{y}^{2}q_{1}-3\tilde{y}q_{1}+\frac{1}{2}q_{1} & \text{in }\left(0,T\right)\times\Omega,\nonumber \\
q_{1}\left|_{\partial\Omega}\right.= & 0 & \text{on }\partial\Omega\times\left[0,T\right],\nonumber \\
q_{1}\left|_{t=T}\right.= & \frac{1}{\gamma\xi^{2}}\left(\tilde{y}_{f}-\tilde{y}\left|_{t=T}\right.\right) & \text{in }\Omega\label{eq:PF adjoint part II.}
\end{align}
(in the weak sense) causes (\ref{eq:PF - final Frecher der}) to reduce
to

\begin{align}
L_{\left(y,\tilde{y}\right)}\left(y,\tilde{y},u,\lambda\right)\left[\left(v,\hat{v}\right)\right]= & \underset{\Omega}{\int}\left(p_{2}-p_{1}\right)\left|_{t=0}\right.v\left|_{t=0}\right.\dd x+\stackrel[0]{T}{\int}\underset{\partial\Omega}{\int}v\left(\nabla p_{1}\cdot\vec{n}+p_{3}\right)\dd S\dd t\\
+ & \underset{\Omega}{\int}\left(Hp_{1}-\gamma\xi^{2}q_{1}+q_{2}\right)\hat{v}\left|_{t=0}\right.\dd x+\stackrel[0]{T}{\int}\underset{\partial\Omega}{\int}\left(q_{3}+\xi^{2}\nabla q_{1}\cdot\vec{n}\right)\hat{v}\left|_{\partial\Omega}\right.\dd S\dd t.\label{eq:PF - reduced frechet derivative}
\end{align}
Lastly, we set
\begin{align}
p_{2}= & p_{1}\left|_{t=0}\right.,\nonumber \\
p_{3}= & -\nabla p_{1}\cdot\vec{n}\left|_{\partial\Omega}\right.,\nonumber \\
q_{2}= & \left(\gamma\xi^{2}q_{1}-Hp_{1}\right)\left|_{t=0}\right.,\nonumber \\
q_{3}= & -\xi^{2}\nabla q_{1}\cdot\vec{n}\left|_{\partial\Omega}\right.\label{eq: PF - setting additonal adjoint variables}
\end{align}
so that (\ref{eq:PF - reduced frechet derivative}) becomes the zero
operator. Next, we notice that the problems (\ref{eq:PF adjoint part I.})
and (\ref{eq:PF adjoint part II.}) can be transformed using

\[
t\rightarrow T-t.
\]
We use the following notation for the transformed variables

\begin{align}
p\left(t\right)= & p_{1}\left(T-t\right),\nonumber \\
q\left(t\right)= & q_{1}\left(T-t\right),\nonumber \\
\tilde{z}\left(t\right)= & \tilde{y}\left(T-t\right).\label{eq:time-transoformations}
\end{align}
This gives rise to the well posed problem

\begin{align}
p_{t}=\Delta p- & \beta\xi q & \text{in }\left(0,T\right)\times\Omega,\nonumber \\
p\left|_{\partial\Omega}\right.= & 0 & \text{on }\partial\Omega\times\left[0,T\right],\nonumber \\
p\left|_{t=0}\right.= & 0 & \text{in }\Omega,\label{eq:PF adjoint part I. - transformed}
\end{align}
\begin{align}
\gamma\xi^{2}q_{t}=\xi^{2}\Delta q & +Hp_{t}-3\tilde{z}^{2}q+3\tilde{z}q-\frac{1}{2}q & \text{in }\left(0,T\right)\times\Omega,\nonumber \\
q\left|_{\partial\Omega}\right.= & 0 & \text{on }\partial\Omega\times\left[0,T\right],\nonumber \\
q\left|_{t=0}\right.= & \frac{1}{\gamma\xi^{2}}\left(\tilde{y}_{f}-\tilde{y}\left|_{t=T}\right.\right) & \text{in }\Omega,\label{eq:PF adjoint part II. - transformed}
\end{align}
where the equations for $p,q$ resemble the heat equation (\ref{eq:PF-heat-equation})
and the Allen-Cahn equation (\ref{eq:PF-Allen-Cahn-equation}), respectively,
but with rather different reaction terms on the right hand side. After
the system (\ref{eq:PF adjoint part I. - transformed})-(\ref{eq:PF adjoint part II. - transformed})
is solved and the additional adjoint variables are set according to
(\ref{eq: PF - setting additonal adjoint variables}), the formula
(\ref{eq:derivative calculation}) can be used to compute the gradient
of $\hat{J}$. Notice, that the first term in (\ref{eq:derivative calculation})
is zero since (\ref{eq:lagrangian adjoint condition}) is satisfied.
Let $s\in U$ be a functional variation, then the formal Fréchet derivative
at point $u\in U$ in direction $s$ reads

\begin{equation}
\hat{J}^{\prime}\left(u\right)\left[s\right]=\alpha\stackrel[0]{T}{\int}\underset{\partial\Omega}{\int}us\dd S\dd t-\stackrel[0]{T}{\int}\underset{\partial\Omega}{\int}p_{3}s\dd S\dd t.\label{eq:PF - reduced derivative computation}
\end{equation}

The numerical experiments performed in Sections \ref{subsec:Controlling-the-Extent},
\ref{subsec:Moving-a-Gap} and \ref{subsec:Moving-a-Crystal-n-t-s-lin}
show that using the linear reaction term (\ref{eq:reaction-term})
produces nonphysical results due to the condition (\ref{eq:PF physicality condition})
being violated by the optimal control. To rectify this, we additionally
consider using an alternative reaction term \cite{PF-Focusing-Latent-Heat}
\begin{equation}
f_{0}\left(y,\tilde{y};\xi\right)=2\tilde{y}\left(1-\tilde{y}\right)\left(\tilde{y}-\frac{1}{2}+\xi\beta\frac{1}{2}\Sigma\left(\tilde{y};\varepsilon_{0},\varepsilon_{1}\right)\left(y_{\text{mt}}-y\right)\right),\label{eq:sigma limiter reaction}
\end{equation}
where $\Sigma\left(p;\varepsilon_{0},\varepsilon_{1}\right)$ is a
differentiable sigmoid function (a limiter) in the form
\begin{equation}
\Sigma\left(\tilde{y};\varepsilon_{0},\varepsilon_{1}\right)=\begin{cases}
0 & \tilde{y}\leq\varepsilon_{0},\\
1 & \tilde{y}\geq\varepsilon_{1},\\
\frac{3\left(\tilde{y}-\varepsilon_{0}\right)^{2}}{\left(\varepsilon_{1}-\varepsilon_{0}\right)^{2}}-\frac{2\left(\tilde{y}-\varepsilon_{0}\right)^{3}}{\left(\varepsilon_{1}-\varepsilon_{0}\right)^{3}} & \tilde{y}\in\left(\varepsilon_{0},\varepsilon_{1}\right).
\end{cases}\label{eq:Sigma-limiter}
\end{equation}
Since (\ref{eq:sigma limiter reaction}) does not lose its physical
interpretation for any value of undercooling \cite{PF-Focusing-Latent-Heat},
replacing (\ref{eq:reaction-term}) with (\ref{eq:sigma limiter reaction})
allows us to obtain meaningful controls for a broader range of experiments
(see Section \ref{subsec:PhaseField-Num-1}).

The derivation of the corresponding adjoint equations is analogous
to the steps leading from (\ref{eq: PF strong defintion of operators})
to (\ref{eq:PF adjoint part I. - transformed}) and (\ref{eq:PF adjoint part II. - transformed}).
The resulting problem reads

\begin{align}
p_{t}=\Delta p- & \beta\xi q\left(\tilde{z}\Sigma\left(\tilde{z};\varepsilon_{0},\varepsilon_{1}\right)-\tilde{z}^{2}\Sigma\left(\tilde{z};\varepsilon_{0},\varepsilon_{1}\right)\right)\text{in }\left(0,T\right)\times\Omega, & \text{in }\left(0,T\right)\times\Omega,\nonumber \\
p\left|_{\partial\Omega}\right.= & 0, & \text{on }\partial\Omega\times\left[0,T\right],\nonumber \\
p\left|_{t=0}\right.= & 0, & \text{in }\Omega,\label{eq:adjoint-better reaction-1}
\end{align}
\begin{align}
\gamma\xi^{2}q_{t}=\xi^{2}\Delta q & +Hp_{t}-3\tilde{z}^{2}q+3\tilde{z}q-\frac{1}{2}q-qh\left(\tilde{z},z\right) & \text{in }\left(0,T\right)\times\Omega,\nonumber \\
q\left|_{\partial\Omega}\right.= & 0 & \text{on }\partial\Omega\times\left[0,T\right],\nonumber \\
q\left|_{t=0}\right.= & \frac{1}{\gamma\xi^{2}}\left(\tilde{y}_{f}-\tilde{y}\left|_{t=T}\right.\right) & \text{in }\Omega,\label{eq:adjoint-better-reaction-2}
\end{align}
where $z\left(t\right)=y\left(T-t\right)$ and
\[
h\left(\tilde{z},z\right)=\beta\left(z-y_{\text{mt}}\right)\left(\Sigma\left(\tilde{z};\varepsilon_{0},\varepsilon_{1}\right)+\tilde{z}\Sigma^{\prime}\left(\tilde{z};\varepsilon_{0},\varepsilon_{1}\right)-2\tilde{z}\Sigma\left(\tilde{z};\varepsilon_{0},\varepsilon_{1}\right)-\tilde{z}^{2}\Sigma^{\prime}\left(\tilde{z};\varepsilon_{0},\varepsilon_{1}\right)\right).
\]
The gradient computation (\ref{eq:PF - reduced derivative computation})
and the relationships (\ref{eq: PF - setting additonal adjoint variables}),
(\ref{eq:time-transoformations}) still hold in this case. Consequently,
the gradient evaluation for the phase field system with (\ref{eq:sigma limiter reaction})
as the reaction term is straightforward.

\section{General Numerical Framework\label{sec:Numerical-Framework-(General)}}

This section introduces the numerical method that is used to solve
the optimization problem (\ref{eq:PF-functional})-(\ref{eq:PF-PF-IC}).
Following the First-Optimize-Then-Discretize Paradigm \cite{Opt-Book-Hind,bookOpt2012},
optimality conditions are provided along with an overview of their
numerical treatment (Section \ref{subsec:Review-of-Opt-Discr}). This
results in several sub-problems, addressed thoroughly in the section
that follows (Section \ref{sec:Numerical-Results}).

\subsection{\label{subsec:Review-of-Opt-Discr}Numerical Treatment using the
First-Optimize-Then-Discretize Paradigm}

When following the classical first-optimize-then-discretize approach
\cite{Opt-Book-Hind,bookOpt2012} formally, optimality conditions
are derived and then discretized to give a well formulated finite
dimensional minimization problem that can then be solved numerically.
Let $\left(\overline{y},\overline{u}\right)$ be a solution of the
problem (\ref{eq:FP - functional to min})-(\ref{eq:FP - state equation})
and $\lambda_{0}$ the solution to the corresponding adjoint problem
(\ref{eq:lagrangian adjoint condition}). Then the optimality conditions
read

\begin{align}
L_{\lambda}\left(\bar{y},\bar{u},\lambda_{0}\right)=e\left(\overline{y},\overline{u}\right)= & 0,\label{eq:optimality conditions - primary}\\
L_{y}\left(\bar{y},\bar{u},\lambda_{0}\right)=e_{y}\left(\overline{y},\overline{u}\right)^{*}\lambda_{0}+J_{y}\left(\overline{y},\overline{u}\right)= & 0,\label{eq:optimality conditions - adjoint}\\
L_{u}\left(\bar{y},\bar{u},\lambda_{0}\right)s=J_{u}\left(\overline{y},\overline{u}\right)s+\left\langle \lambda_{0},e_{u}\left(\overline{y},\overline{u}\right)s\right\rangle _{Z}\geq & 0\quad\forall s\in U.\label{eq:optimality conditions - integral computation}
\end{align}
The conditions (\ref{eq:optimality conditions - primary})-(\ref{eq:optimality conditions - integral computation})
state that $\left(\overline{y},\overline{u},\lambda_{0}\right)$ satisfy
the state equation (\ref{eq:state-equation}), the adjoint problem
(\ref{eq:lagrangian adjoint condition}), and the necessary condition
for a local minimum. The equations (\ref{eq:optimality conditions - primary})
and (\ref{eq:optimality conditions - adjoint}) are then discretized
using a suitable numerical method. For example, if (\ref{eq:optimality conditions - primary})
and (\ref{eq:optimality conditions - adjoint}) represent a system
of PDE's with the requisite (initial and/or boundary) conditions,
the finite difference method (FDM) or the finite element method (FEM)
can be applied. Condition (\ref{eq:optimality conditions - integral computation})
is approximated numerically by plugging in the suitably interpolated
numerical solutions of (\ref{eq:optimality conditions - primary}),
(\ref{eq:optimality conditions - adjoint}).

For the discussion of the individual steps of the procedure, note
that particular form of the primary problem (\ref{eq:optimality conditions - primary})
is given by (\ref{eq:PF-heat-equation})-(\ref{eq:PF-PF-IC}) and
the adjoint problem (\ref{eq:optimality conditions - adjoint}) takes
the form of (\ref{eq:PF adjoint part I. - transformed})-(\ref{eq:PF adjoint part II. - transformed})
or (\ref{eq:adjoint-better reaction-1})-(\ref{eq:adjoint-better-reaction-2}).
Both these problems can be solved numerically by the FDM, using a
fixed time step $\Delta t$ and a rectangular spatial mesh. Let the
numerical solution of the state equation $y_{h}$, the control $u_{h}$,
and the adjoint variable $\lambda_{h}$ all be mesh functions, i.e.
functions only defined on a discrete (sub)set of mesh points. We denote
the spaces of these functions as $Y_{h}$, $U_{h}$, and $Z_{h}$,
respectively. The subscript $h$ merely indicates that the dimension
of these spaces is finite. The actual dimensions of these spaces depend
on the number of mesh points and time levels used. The algorithm to
solve the discrete counterpart of (\ref{eq:optimality conditions - primary})-(\ref{eq:optimality conditions - integral computation})
then reads:
\begin{enumerate}
\item \label{enu:num-general-start}Start with an initial guess $u_{h}:=u_{h,0}\in U_{h}$.
\item \label{enu:num-general-solve-primary}Use the FDM numerical solver
to compute the solution of the primary problem $y_{h}\left(u_{h}\right)\in Y_{h}$.
\item \label{enu:num-general-solve-adjoint}Using the solution of the primary
problem $y_{h}\left(u_{h}\right)$, run the FDM adjoint problem solver
to compute $\lambda_{h}\left(y_{h}\left(u_{h}\right)\right)\in Z_{h}$.
\item \label{enu:num-general-compute-gradient}Compute all components of
$\nabla\hat{J}_{h}$, i.e. the partial derivatives of the discrete
analogue of $\hat{J}$ with respect to all the basis vectors of $U_{h}$,
using a discretization of the left hand side of (\ref{eq:optimality conditions - integral computation})
with $\lambda_{h}$ obtained in the previous step.
\item \label{enu:num-general-gradient-descent}Perform one step of gradient
descent by updating $u_{h}$ as $u_{h}:=u_{h}-\varepsilon\nabla\hat{J}_{h}$
where $\varepsilon>0$ is a given step size.
\item \label{enu:num-general-end}Go to step \ref{enu:num-general-solve-primary}
unless a suitable stopping criterion is satisfied.
\end{enumerate}
Some of the possible stopping criteria include exceeding a given maximum
number of iterations or $\left|\nabla\hat{J}_{h}\right|$ falling
below a predetermined threshold \cite{Opt-Book-Hind}.

\subsection{\label{sec:Numerical-Results}Details of the Numerical Method}

In this section, the detailed numerical treatment of each of the problems
outlined in Section \ref{subsec:Adjoint-Method-Derivative-PF} is
laid out.

As the numerical results presented later in Section \ref{sec:Numerical-Results-1}
are limited to one and two spatial dimensions, it is enough to follow
the notation for the two-dimensional case. The one-dimensional case
then arises as a straightforward simplification.

\subsubsection{Finite Difference Scheme}

Let $\Omega=\left(0,L_{x_{1}}\right)\times\left(0,L_{x_{2}}\right)\subset\R^{2}$,
$x=\left(x_{1},x_{2}\right)$, and $T>0$. Both the adjoint and the
primary problems are solved using the finite difference method on
a uniform mesh. The time step and the spatial mesh resolution are
\[
\Delta t=\frac{T}{N_{t}-1},\;\Delta x_{1}=\frac{L_{x_{1}}}{N_{x_{1}}-1},\;\Delta x_{2}=\frac{L_{x_{2}}}{N_{x_{2}}-1},
\]
where $N_{t}$ denotes the number of time layers and $N_{x_{1}},N_{x_{2}}$
the number of mesh points in the $x_{1}$ and $x_{2}$ directions,
respectively. The explicit Euler scheme is employed. Let $f$ be a
real-valued function defined on $\left(0,T\right)\times\Omega$ and
$f_{h}$ be its approximation by the respective mesh function. Then
we replace the time derivative of $f$ by
\begin{align}
f_{t}\left(k\Delta t,i\Delta x_{1},j\Delta x_{2}\right)\approx f_{h,t}\left(k\Delta t,i\Delta x_{1},j\Delta x_{2}\right)\equiv & \frac{f_{h}\left(\left(k+1\right)\Delta t,i\Delta x_{1},j\Delta x_{2}\right)-f_{h}\left(k\Delta t,i\Delta x_{1},j\Delta x_{2}\right)}{\Delta t}.\label{eq:FDM-time-derivative}
\end{align}
For the discretization of the Laplacian of $f$, the central difference
quotient
\begin{align}
\Delta f\left(k\Delta t,i\Delta x_{1},j\Delta x_{2}\right)\approx\Delta_{h}f_{h}\left(k\Delta t,i\Delta x_{1},j\Delta x_{2}\right)\equiv\label{eq:FDM-Laplacian}\\
\frac{f_{h}\left(k\Delta t,\left(i+1\right)\Delta x_{1},j\Delta x_{2}\right)-2f_{h}\left(k\Delta t,i\Delta x_{1},j\Delta x_{2}\right)+f_{h}\left(k\Delta t,\left(i-1\right)\Delta x_{1},j\Delta x_{2}\right)}{\left(\Delta x_{1}\right)^{2}}\\
+\frac{f_{h}\left(k\Delta t,i\Delta x_{1},\left(j+1\right)\Delta x_{2}\right)-2f_{h}\left(k\Delta t,i\Delta x_{1},j\Delta x_{2}\right)+f_{h}\left(k\Delta t,i\Delta x_{1},\left(j-1\right)\Delta x_{2}\right)}{\left(\Delta x_{2}\right)^{2}}
\end{align}
is used.

\subsubsection{Numerical Integration}

The integral in (\ref{eq:PF - reduced derivative computation}) is
evaluated using a piece-wise constant interpolation as follows. Let
$\varPi$ be a subset of the Cartesian rectangular mesh representing
the boundary ($\left[0,T\right]\times\partial\Omega$) and let $f_{h}:\varPi\rightarrow\mathbb{R}$,
$s_{h}:\varPi\rightarrow\mathbb{R}$ be mesh functions. To improve
readability, the subscript $h$ is dropped for both $f_{h}$ and $s_{h}$
in this section.

Let $S_{\varPi}$ be an interpolation operator that transforms a mesh
function on $\varPi$ to a piecewise (a.e.) constant function on $\left[0,T\right]\times\partial\Omega$.
More precisely, the definition of $S_{\varPi}$ on the bottom edge
($x_{2}=0$) of a rectangular 2D domain reads
\[
S_{\varPi}f\left(t,x_{1},0\right)=f_{k,i,0}\text{ for }\left(t,x_{1},0\right)\in M_{k,i,0},
\]
where
\begin{align*}
M_{k,i,0}= & \left(\max\left(0,\left(k-1\right)\Delta t+\frac{\Delta t}{2}\right),\min\left(T,k\Delta t+\frac{\Delta t}{2}\right)\right)\\
\times & \left(\max\left(0,\left(i-1\right)\Delta x_{1}+\frac{\Delta x}{2}\right),\min\left(L_{x_{1}},i\Delta x_{1}+\frac{\Delta x_{1}}{2}\right)\right)\times\left\{ 0\right\} ,
\end{align*}
$k$ denotes the time level and $i$ represents the mesh point position
along the $x_{1}$ axis. Intuitively, the set $M_{k,i,0}$ can be
viewed as a rectangle in $\left(0,T\right)\times\partial\Omega$ centered
in $\left(k\Delta t,i\Delta x_{1},0\right)$ except for $k=0,N_{t},i=0,N_{x}$,
where there is a ``cut-off'' at the boundary. The operator $S_{\varPi}$
is defined analogously for the other edges. Since the numerical scheme
described by (\ref{eq:FDM-time-derivative})-(\ref{eq:FDM-Laplacian})
does not require the points in the corners of the spatial domain $\Omega$,
they are left out in the approximation of (\ref{eq:PF - reduced derivative computation}),
giving rise to
\begin{align}
\stackrel[0]{T}{\int}\underset{\partial\Omega}{\int}\left(S_{\varPi}f\right)\left(S_{\varPi}s\right)\dd S\dd t\approx & \stackrel[k=0]{N_{t}-1}{\sum}\stackrel[i=1]{N_{x_{1}}-1}{\sum}\left[f_{k,i,0}s_{k,i,0}+f_{k,i,N_{x_{2}}}s_{k,i,N_{x_{2}}}\right]\Delta t\Delta x_{1}\nonumber \\
+ & \stackrel[k=0]{N_{t}-1}{\sum}\stackrel[j=1]{N_{x_{2}}-1}{\sum}\left[f_{k,0,j}s_{k,0,j}+f_{k,N_{x_{1}},j}s_{k,N_{x_{1}},j}\right]\Delta t\Delta x_{2},\label{eq:integral-approx}
\end{align}
where $j$ represents the mesh point position along the $x_{2}$ axis
and $\alpha=0$ is chosen for simplicity. Notice that omitting the
corners of the spatial domain does not affect the convergence of (\ref{eq:integral-approx})
as $N_{t},N_{x_{1}},N_{x_{2}}\to+\infty$.

In particular, (\ref{eq:integral-approx}) will be used along with
a finite difference approximation $p_{3,h}$ of $p_{3}$ in (\ref{eq: PF - setting additonal adjoint variables})
to get the computation rule for the $k,i,j$-th component of the gradient
as

\begin{equation}
\stackrel[0]{T}{\int}\underset{\partial\Omega}{\int}\left(S_{\varPi}p_{3,h}\right)\left(S_{\varPi}e_{k,i,j}\right)\dd S\dd t=\begin{cases}
\left(p_{3,h}\right)_{k,i,j}\Delta t\Delta x_{1} & \text{if }j=\left\{ 1,N_{x_{2}}-1\right\} ,\\
\left(p_{3,h}\right)_{k,i,j}\Delta t\Delta x_{2} & \text{if }i=\left\{ 1,N_{x_{1}}-1\right\} ,
\end{cases}\label{eq:grad-computation}
\end{equation}
where $e_{k,i,j}$ is the characteristic function of $M_{k,i,j}.$
Since $p$ (see (\ref{eq:PF adjoint part I. - transformed})) is subject
to the homogeneous Dirichlet boundary condition, the approximation
$p_{3,h}$ in the gradient computation (\ref{eq:grad-computation})
reduces to the interior values of $p_{1}$ (see (\ref{eq:time-transoformations}),
(\ref{eq: PF - setting additonal adjoint variables})) adjacent to
the boundary $\partial\Omega$.

\section{\label{sec:Numerical-Results-1}Numerical Results}

Utilizing the adjoint formulations derived in Section \ref{subsec:Adjoint-Method-Derivative-PF}
and the corresponding numerical treatments detailed in Section \ref{sec:Numerical-Framework-(General)},
the problem (\ref{eq:PF-functional})-(\ref{eq:PF-PF-IC}) is solved
numerically. In Section \ref{subsec:PhaseField-Num}, simulations
in one spatial dimension using the linear reaction term (\ref{eq:reaction-term})
are performed and the effects of regularization, changes in final
time and different initial guesses are discussed. Several scenarios
where the bound (\ref{eq:PF physicality condition}) is not violated
are shown. Some of the experiments, however, show that the optimal
control violates the bound (\ref{eq:PF physicality condition}) in
order to achieve the desired crystal shape (see some of the experiments
in Section \ref{subsec:PhaseField-Num} and the experiment of Section
\ref{subsec:Moving-a-Gap}). When this is the case, we call a simulation
or associated optimal control \emph{non-realistic.}

As expected, this inadequacy can be observed in experiments with two
spatial dimensions as well (see Section \ref{subsec:PhaseField-Num-1}).
In Section \ref{subsec:Moving-a-Crystal-n-t-s-lin}, an experiment
that directly compares the behaviors obtained using the linear reaction
term (\ref{eq:reaction-term}) and the alternative reaction term (\ref{eq:sigma limiter reaction})
is described. The last experiment of the section shows how the reaction
term (\ref{eq:sigma limiter reaction}) can be used to find a realistic
optimal control that separates a crystal into two.

Before discussing each of the simulations in detail, some terminology
and notes will be presented. When commenting on experiments, we recall
the definitions (\ref{eq:solid-subdomain})-(\ref{eq:phase-interface})
and the related terminology introduced in Section \ref{subsec:Adjoint-Method-Derivative-PF},
using the terms ``crystal'' and ``solid subdomain'' interchangeably
to refer to $\Omega_{\text{s}}\left(t\right)$. It is well known \cite{PF-Focusing-Latent-Heat,Benes-Math_comp_aspects_solid,CaginalpAnalysis86}
that inside of $\Omega_{\text{s}}\left(t\right)$ and $\Omega_{\text{l}}\left(t\right)$,
the value of the phase field $\tilde{y}$ is very close to $1$ and
$0$, respectively, except for a thin transition layer between the
phases near $\Gamma\left(t\right)$. For both the models considered,
this transition layer (diffuse interface) has a characteristic profile
with a thickness proportional to $\xi$. Specifically, in the one-dimensional
case, let $x_{0}\left(t\right)$ denote the position of the phase
interface (i.e. $\left\{ x_{0}\left(t\right)\right\} =\Gamma\left(t\right)$)
and $\Omega_{\text{s}}\left(t\right)=\left[0,x_{0}\left(t\right)\right)$,
$\Omega_{\text{l}}\left(t\right)=\left(x_{0}\left(t\right),1\right]$.
Then
\begin{equation}
\tilde{y}\left(t,x\right)=\frac{1}{2}\left[1-\tanh\left(\frac{x-x_{0}\left(t\right)}{2\xi}\right)\right]+o\left(\xi\right).\label{eq:PF-characteristic-profile}
\end{equation}
This explanation can be extended to multiple dimensions, then the
profile (\ref{eq:PF-characteristic-profile}) is maintained in the
direction normal to the interface \cite{PF-Focusing-Latent-Heat,Benes-Asymptotics}.

Denote the numerical approximation of $y$ as $y_{h}$, the numerical
approximation of $\tilde{y}$ as $\tilde{y}_{h}$ and let $P_{h}:Y\to Y_{h}$
be the projection operator onto the mesh.

\subsection{\label{subsec:PhaseField-Num}Dirichlet Boundary Condition Control
for the Phase Field Problem in 1D}

All the experiments in this section are performed in one spatial dimension.
Let all the physical parameters be set according to the values in
Table \ref{tab:Settings-for-the-physical-params-PF}. These settings
do not necessarily correspond to any real material, they do however
serve to illustrate the capabilities and deficiencies of the optimization
problem (\ref{eq:PF-functional})-(\ref{eq:PF-PF-IC}) with the linear
reaction term (\ref{eq:reaction-term}). Additionally, the influence
of regularization $\alpha$, changes in final time $T$ and initial
guess for the control are discussed in this section also.

\begin{table}
\caption{\label{tab:Settings-for-the-physical-params-PF}Parameter settings
for the phase field simulations in Section \ref{subsec:PhaseField-Num}.}

\centering{}%
\begin{tabular}{ccl}
\toprule
Param. & Value & Physical Meaning\tabularnewline
\midrule
$\gamma$ & 1 & coefficient of attachment kinetics\tabularnewline
$\beta$ & 2 & dimensionless representation of supercooling\tabularnewline
$\xi$ & 0.005 & interface thickness scaling\tabularnewline
$y_{\text{mt}}$ & 0.5 & melting temperature\tabularnewline
$H$ & 1 & latent heat\tabularnewline
$L_{x}$ & 1 & spatial dimension in the $x$ direction\tabularnewline
\bottomrule
\end{tabular}
\end{table}

\subsubsection{\label{subsec:Controlling-the-Extent}Controlling the Extent of Crystal
Growth}

\begin{figure}
\begin{centering}
\includegraphics[width=0.95\figwidth]{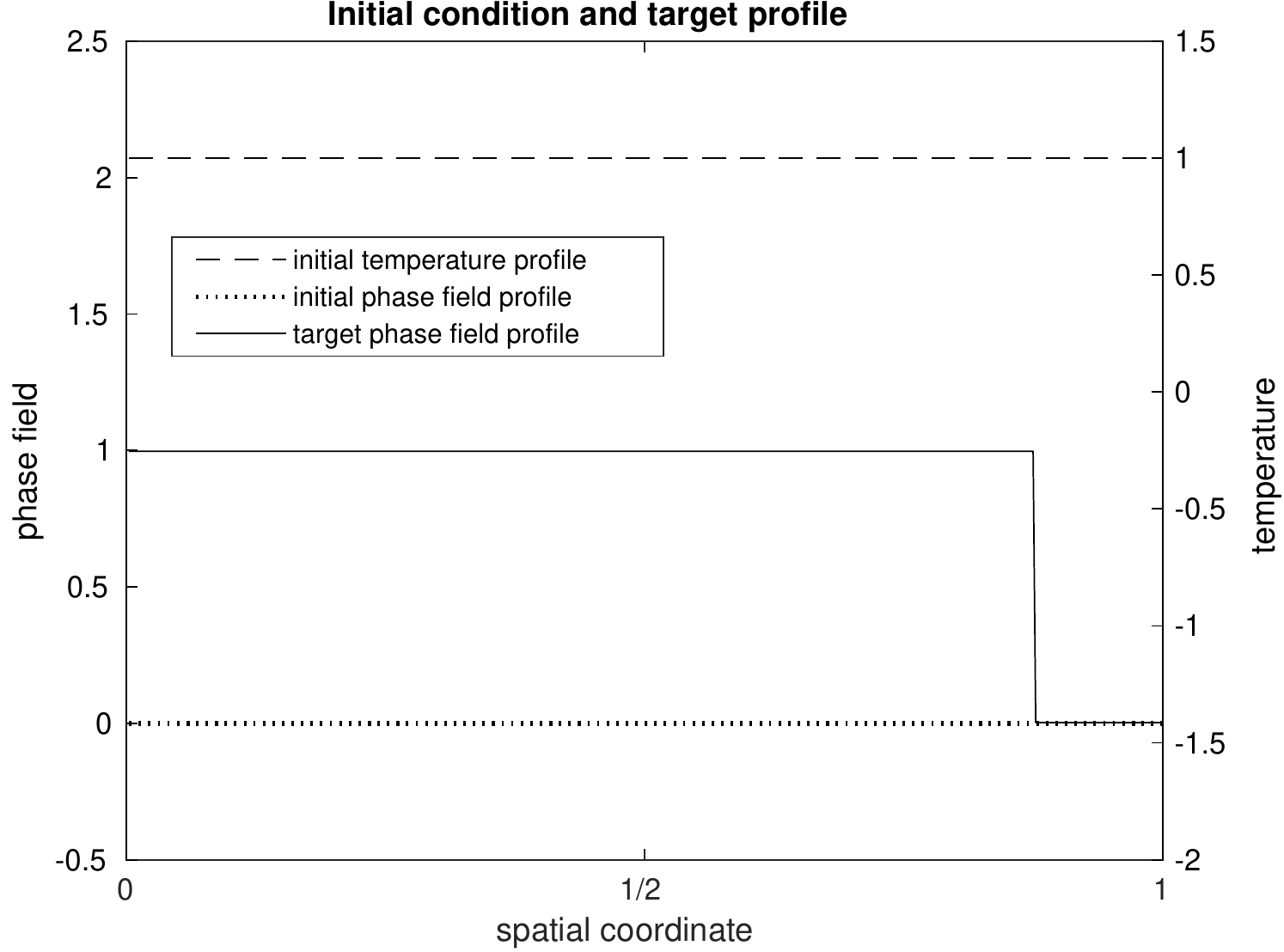}
\par\end{centering}
\caption{\label{fig:Initial conditions and target profile ex1-3}The initial
temperature and phase field spatial profiles $\tilde{y}_{\text{ini}},y_{\text{ini}}$
along with the target profile $\tilde{y}_{f}$ for experiments 1,
2, 3. The values of the boundary condition $\tilde{y}_{\text{bc}}$
are given by (\ref{eq:experiments-123-PF-BCs}).}
\end{figure}

We attempt to find a control that will produce a crystal $\Omega_{\text{s},f}$
of prescribed length inside of the spatial domain $\Omega$ at the
fixed final time $T>0$. We set the phase field boundary condition
to

\begin{align}
\tilde{y}_{\text{bc}}\left(t,0\right)= & 1 & \forall t\in\left[0,T\right),\nonumber \\
\tilde{y}_{\text{bc}}\left(t,1\right)= & 0 & \forall t\in\left[0,T\right),\label{eq:experiments-123-PF-BCs}
\end{align}
which creates a nucleation site at $x=0$. The initial conditions
$\tilde{y}_{\text{ini}},y_{\text{ini}}$ and the target profile $\tilde{y}_{f}$
are depicted in Figure \ref{fig:Initial conditions and target profile ex1-3}.\textbf{
}In this case, $\tilde{y}_{f}$ is merely the characteristic function
of the target crystal shape $\Omega_{\text{s},f}$. Alternatively,
$\tilde{y}_{f}$ could be chosen as a continuous function with the
characteristic shape across the interface (\ref{eq:PF-characteristic-profile}),
as demonstrated further in Sections \ref{subsec:Keeping-Crystal-Separation}
and \ref{subsec:Moving-a-Gap}. The initial guess for the control
is

\begin{align}
u_{0}\left(t,0\right)= & 0 & \forall t\in\left[0,T\right),\nonumber \\
u_{0}\left(t,1\right)= & 1 & \forall t\in\left[0,T\right).\label{eq:initial guess asy}
\end{align}

In this setting, three numerical experiments are presented. The parameters
including the spatial mesh resolution, the number of time steps as
well as the difference of the final solution from the prescribed phase
field profile are summarized in Table \ref{tab:Controlling-the-Extent}.
\begin{table}
\caption{\label{tab:Controlling-the-Extent}Settings for experiments 1, 2,
3 and the respective values of the difference (error) from the prescribed
profile.}

\centering{}%
\begin{tabular}{cccc}
\toprule
 & \multicolumn{3}{c}{Simulation number}\tabularnewline
Parameter & 1 & 2 & 3\tabularnewline
\midrule
number of time steps $N_{t}$ & $4\cdot10^{5}$ & $4\cdot10^{5}$ & $4\cdot10^{5}$\tabularnewline
number of grid points $N_{x}$ & $400$ & $400$ & $400$\tabularnewline
initial control given by & (\ref{eq:initial guess asy}) & (\ref{eq:initial guess asy}) & (\ref{eq:initial guess asy})\tabularnewline
final time $T$ & $0.1$ & $0.05$ & $0.05$\tabularnewline
regularization parameter $\alpha$ & $0$ & $0$ & $5\cdot10^{-11}$\tabularnewline
gradient descent step size $\varepsilon$ & $3\cdot10^{15}$ & $2\cdot10^{16}$ & $2\cdot10^{16}$\tabularnewline
number of iterations & $100$ & $100$ & $100$\tabularnewline
$\left\Vert \tilde{y}_{h}-P_{h}\tilde{y}_{f}\right\Vert _{2}$ at
$t=T$ & $1.117846$ & $1.276584$ & $13.64326$\tabularnewline
\bottomrule
\end{tabular}
\end{table}

\begin{figure}
\begin{centering}
\includegraphics[width=0.95\figwidth]{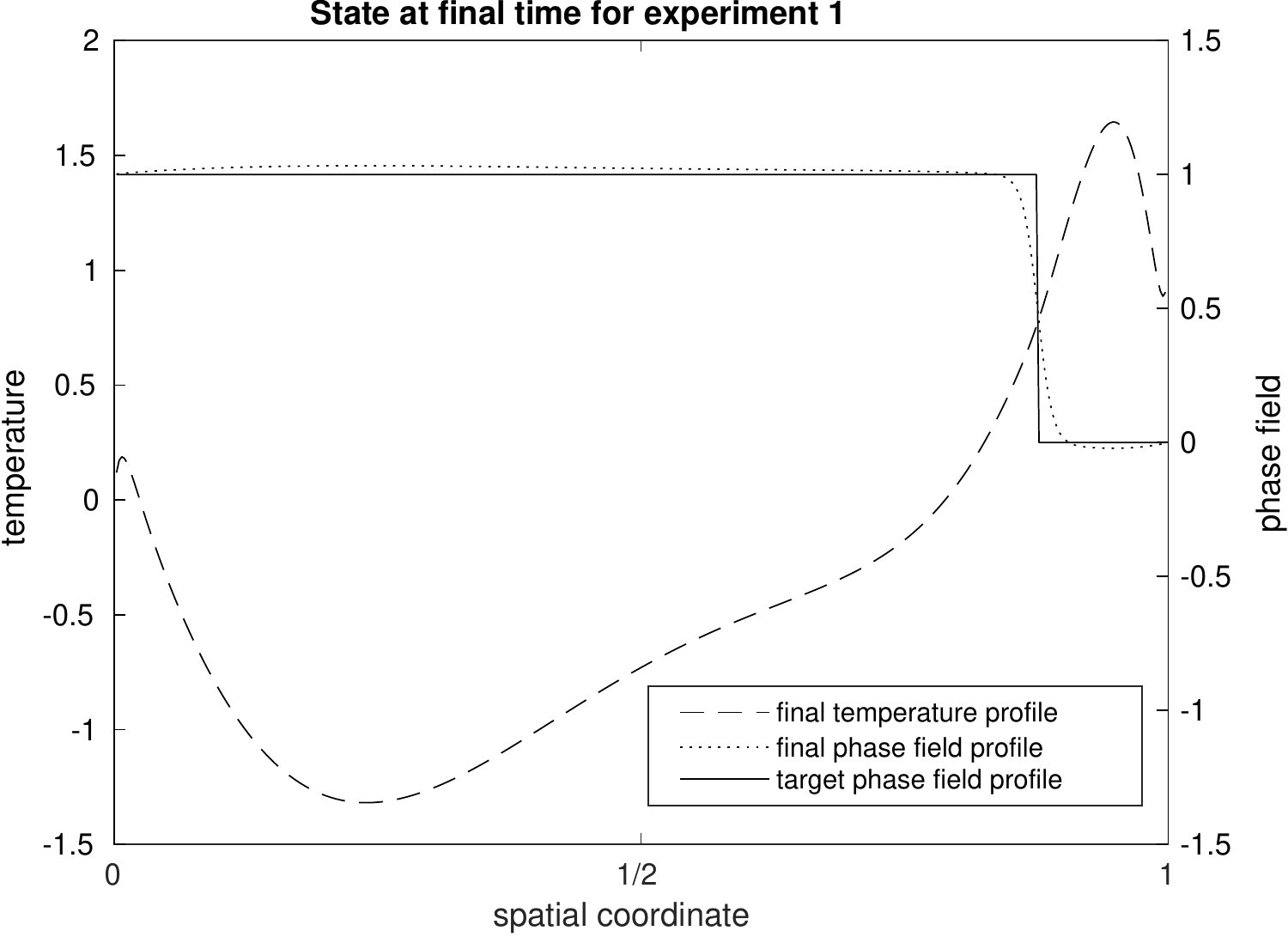}\medskip{}
\par\end{centering}
\begin{centering}
\includegraphics[width=0.95\figwidth]{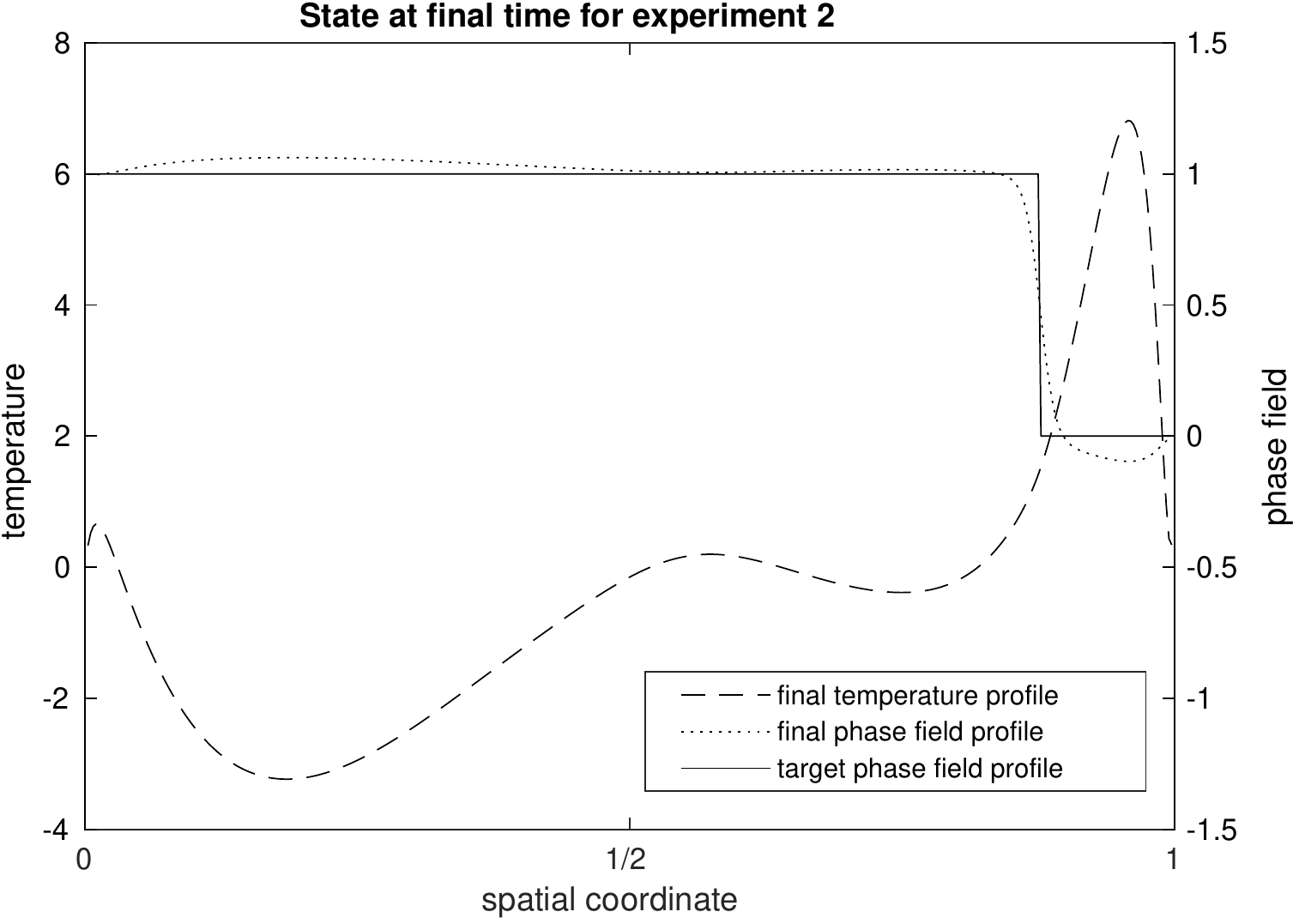}\medskip{}
\par\end{centering}
\begin{centering}
\includegraphics[width=0.95\figwidth]{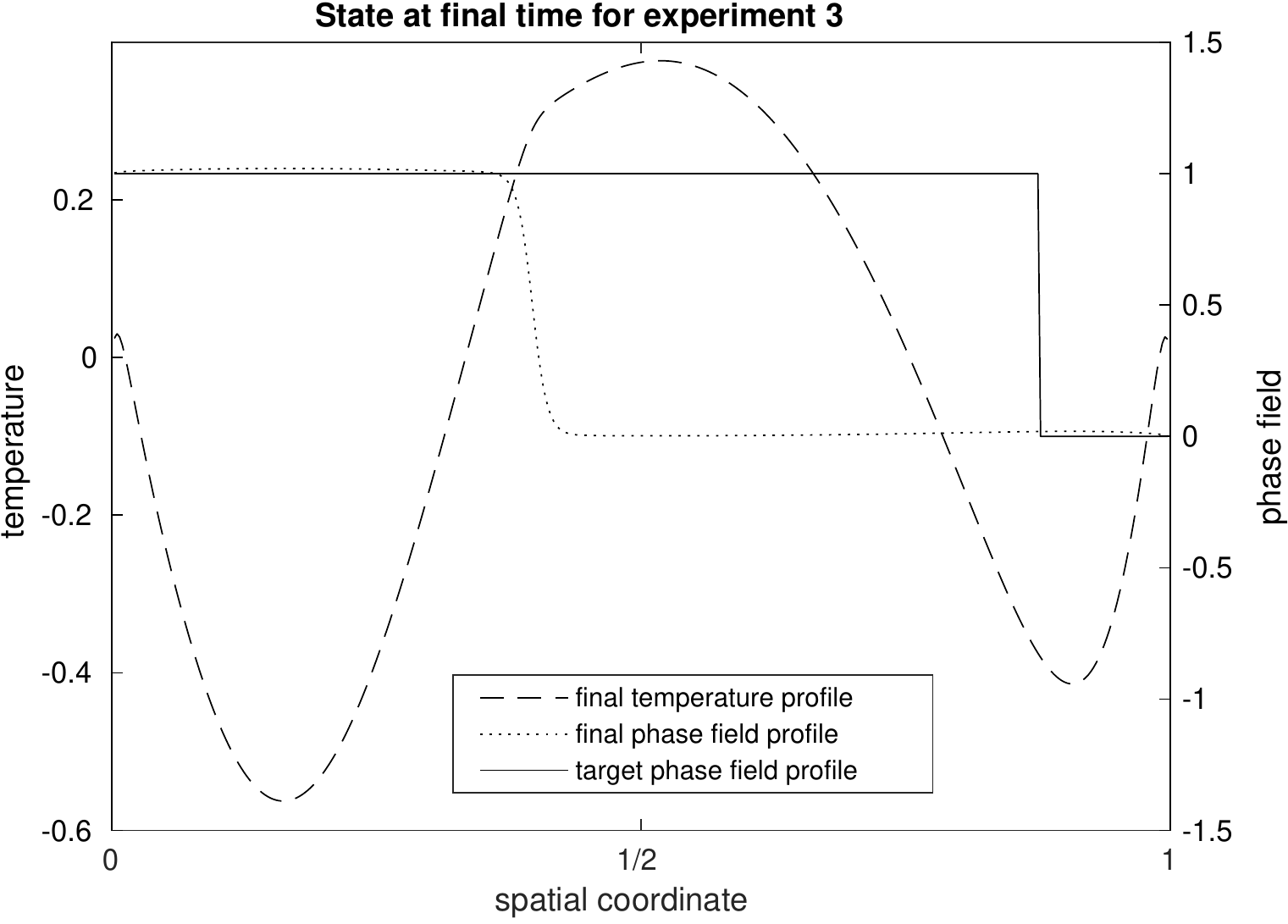}
\par\end{centering}
\caption{\label{fig:Target-Profile-of-EX1-3}Final temperature and phase field
spatial profiles of experiments 1, 2, 3. We observe that $\tilde{y}_{h}$
reaches the target $P_{h}\tilde{y}_{f}$ in experiments 1 and 2. In
experiment 3, the interface of $\tilde{y}_{h}$ does not reach its
target $P_{h}\tilde{y}_{f}$ because sufficient regularization is
added to prevent non-realistic behavior.}
\end{figure}

\begin{figure}
\begin{centering}
\includegraphics[width=0.95\figwidth]{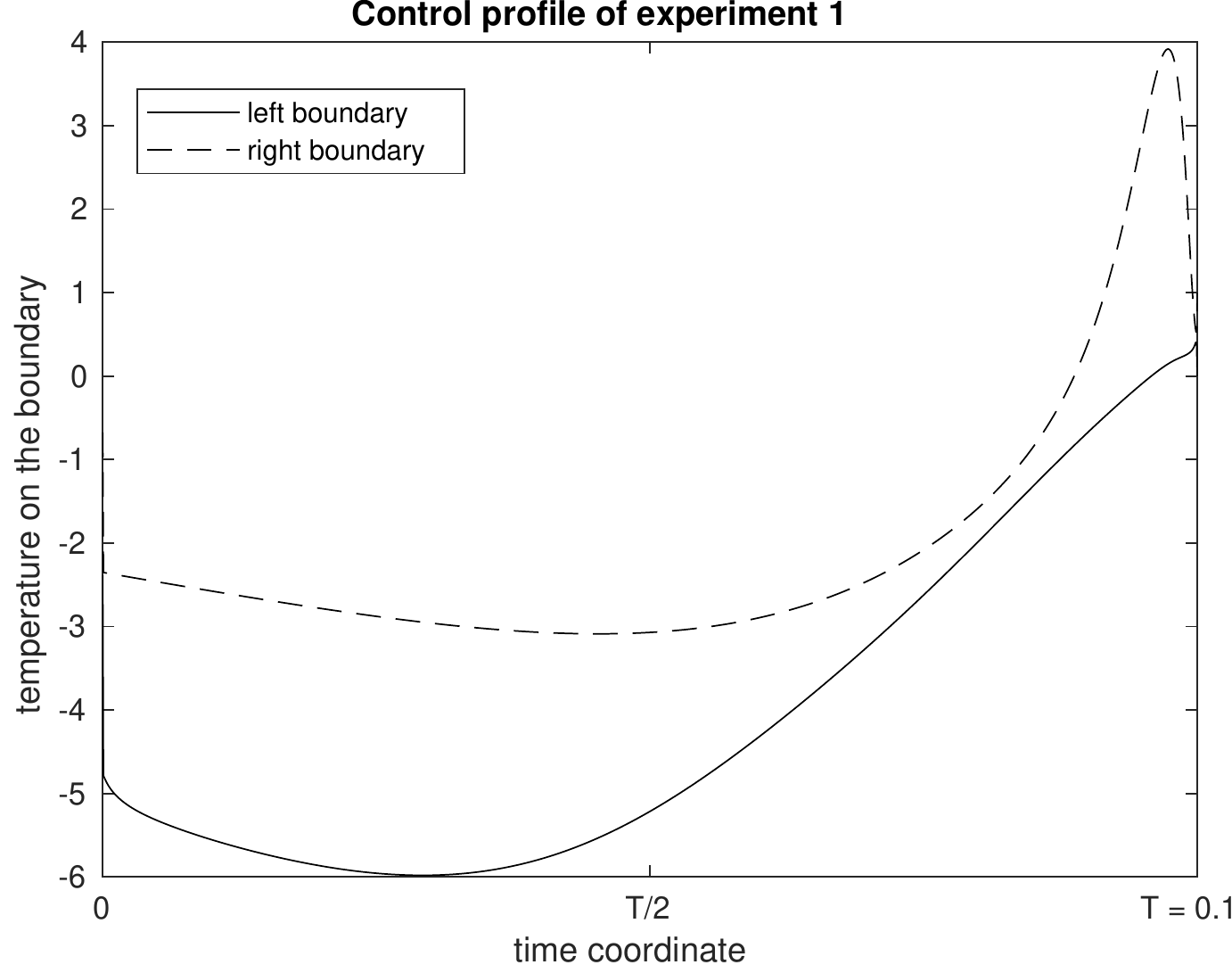}\medskip{}
\par\end{centering}
\begin{centering}
\includegraphics[width=0.95\figwidth]{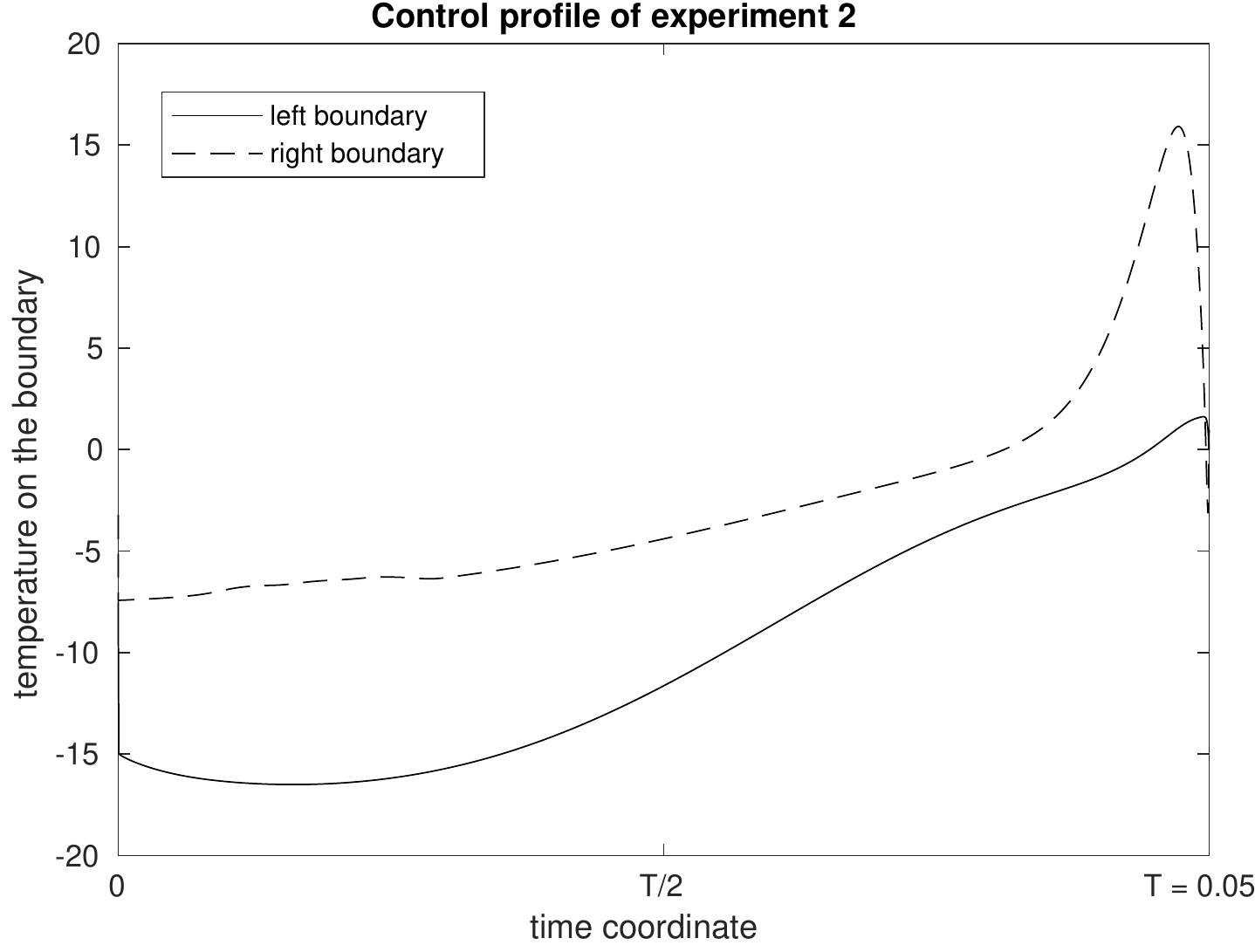}\medskip{}
\par\end{centering}
\begin{centering}
\includegraphics[width=0.95\figwidth]{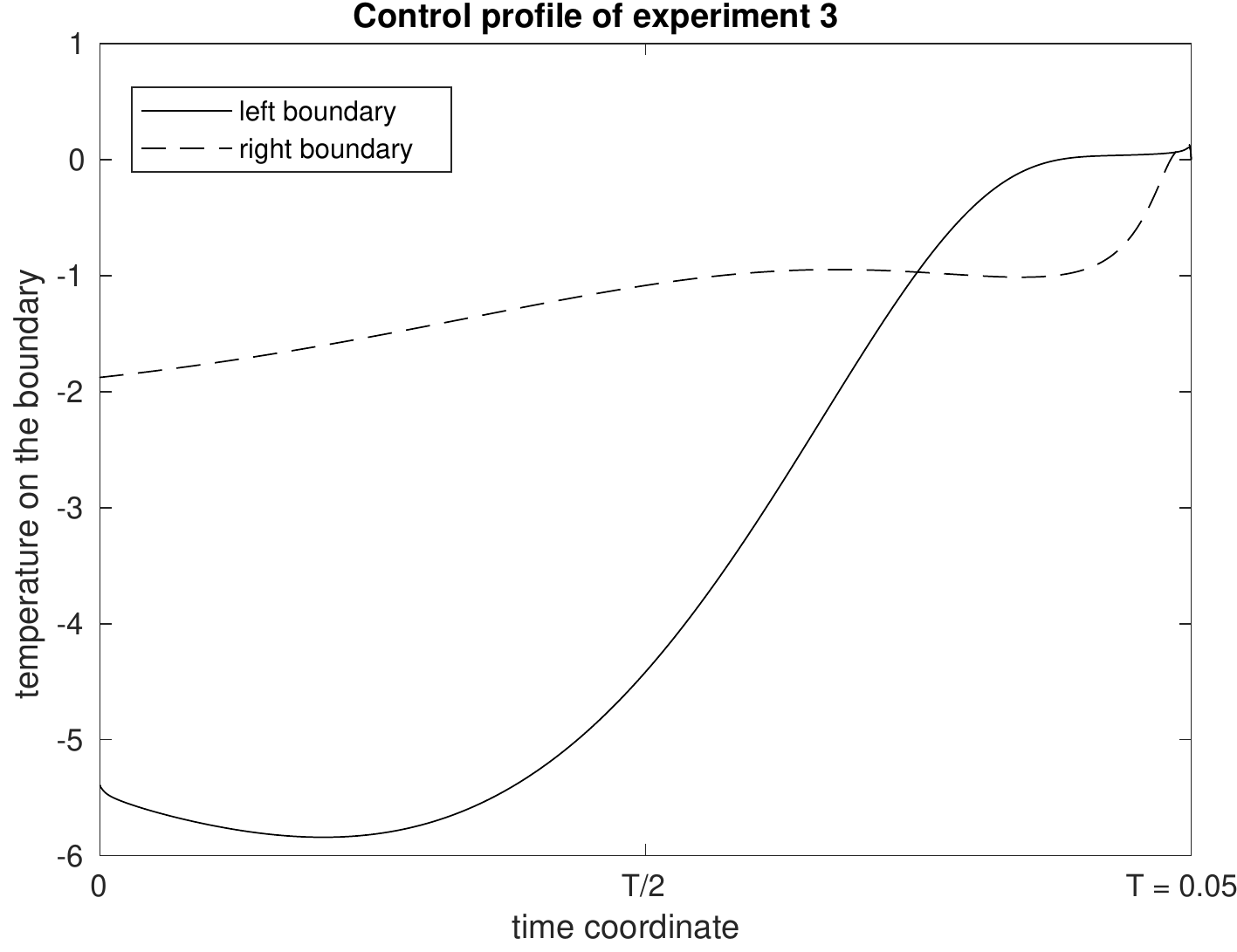}
\par\end{centering}
\caption{\label{fig:Control-profile-1-3}Optimized temporal control profiles
of experiments 1, 2, 3. In experiment 2, the values of the temperature
violate the bound (\ref{eq:PF physicality condition}) due to the
short final time $T$. The regularization applied in experiment 3
fixes this issue at the cost of not attaining the target profile (see
Figure \ref{fig:Target-Profile-of-EX1-3}).}
\end{figure}

For each of the experiments, the resulting temperature and phase field
spatial profiles at final time $T$ are depicted in Figure \ref{fig:Target-Profile-of-EX1-3}.
Figure \ref{fig:Control-profile-1-3} shows the respective temporal
control profiles of the Dirichlet boundary condition.

In experiment 1, the optimized control leads to a good agreement of
the solution with the target profile at the final time $T=0.1$. The
values of $u_{h}$ stay within the limits given by (\ref{eq:PF physicality condition})
and thus the solution is realistic. This experiment was also repeated
several time with different nonzero values of the regularization parameter
$\alpha$. Up to the values $\alpha\approx10^{-7}$, negligible impact
of the regularization was observed. In experiment 2, the final time
is halved, i.e. $T=0.05$. The target profile was still obtained at
the cost of violating the bound (\ref{eq:PF physicality condition}).
Experiment 3 differs from experiment 2 by setting the regularization
parameter to $\alpha=5\cdot10^{-11}$ (i.e. $\alpha\ll10^{-7}$, cf.
Experiment 1). This was sufficient to keep the simulation realistic
but the target profile was not obtained. Such a situation can be interpreted
as the time $T$ being too short for the target profile to be reached.
This is true only in the context of the linear model, where $f_{0}$
is defined by (\ref{eq:reaction-term}). The use of a more advanced
model, that is not limited in this way is discussed in Sections \ref{subsec:Moving-a-Crystal-n-t-s-lin}
and \ref{subsec:Separating-a-Crystal-sigma}.

\subsubsection{\label{subsec:Keeping-Crystal-Separation}Keeping Crystal Separation}

\begin{figure}
\begin{centering}
\includegraphics[width=0.95\figwidth]{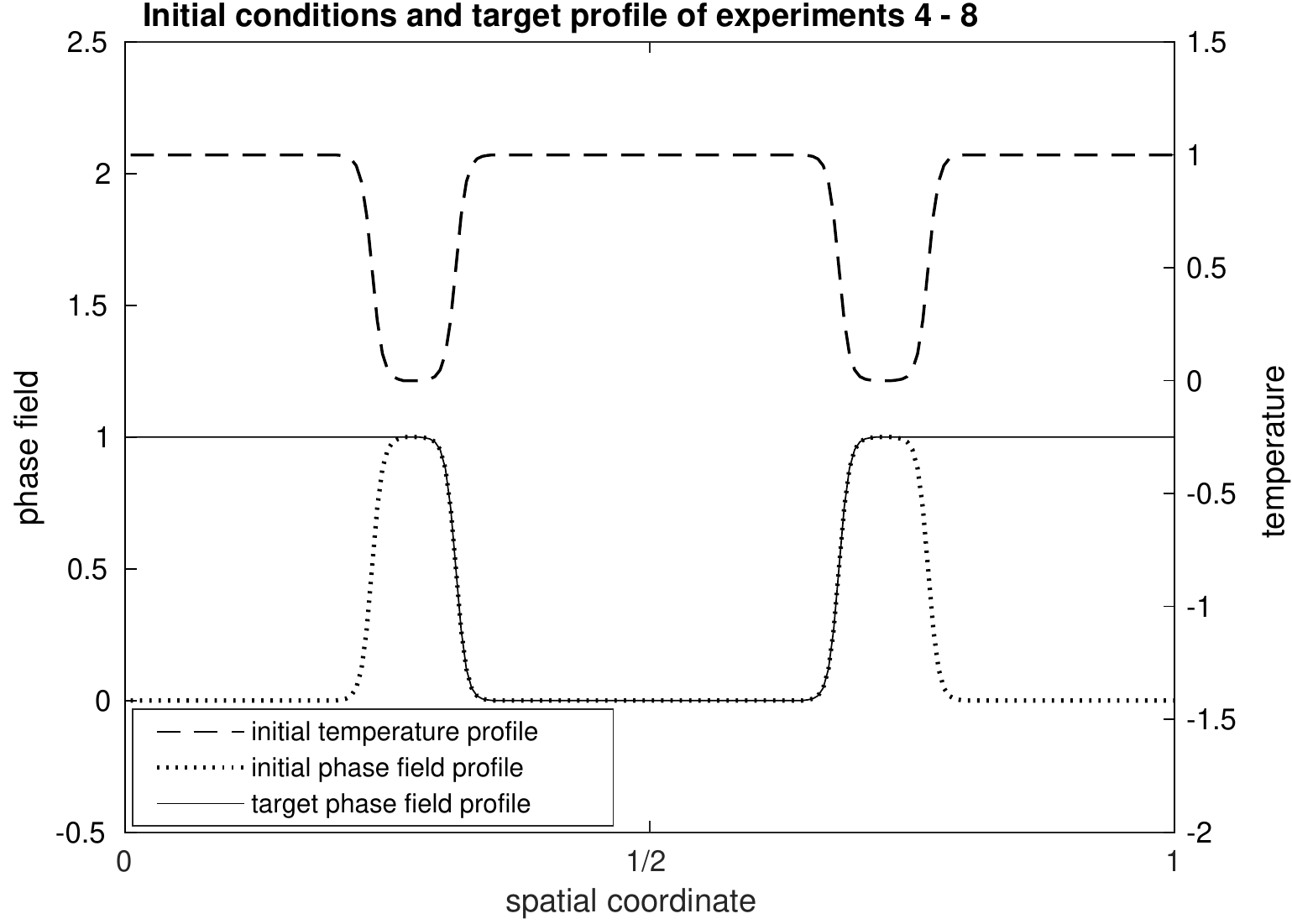}
\par\end{centering}
\caption{\label{fig:The-initial-temperature-4-8}The initial temperature and
phase field profiles $\tilde{y}_{\text{ini}},y_{\text{ini}}$ along
with the target profile $\tilde{y}_{f}$ for experiments 4 through
8. The values of the boundary condition $\tilde{y}_{\text{bc}}$ are
given by (\ref{eq:experiments-4-8-PF-BCs}).}
\end{figure}

In this series of experiments, we attempt to keep two symmetrically
placed crystals in the interior of the domain separated from each
other, while letting them grow toward the boundaries. The target profile
$\tilde{y}_{f}$ and initial conditions $y_{\text{ini}}$, $\tilde{y}_{\text{ini}}$
are depicted in Figure \ref{fig:The-initial-temperature-4-8}. In
contrast to Section \ref{subsec:Controlling-the-Extent}, $\tilde{y}_{f}$
is a continuous function with a characteristic shape across the interface
(\ref{eq:PF-characteristic-profile}). The boundary condition for
the phase field reads

\begin{align}
\tilde{y}_{\text{bc}}\left(t,x\right) & =1 & \forall x\in\left\{ 0,1\right\} ,\;\forall t\in\left[0,T\right).\label{eq:experiments-4-8-PF-BCs}
\end{align}
First, we use the non-symmetric initial control guess (\ref{eq:initial guess asy})
to observe its effects on the obtained temporal profile of $u_{h}$
and overall quality of minimization expressed by the error $\left\Vert \tilde{y}_{h}-P_{h}\tilde{y}_{f}\right\Vert _{2}$
. Then we add regularization and use the symmetric initial guess

\begin{align}
u_{0}\left(x,t\right) & =0 & \forall x\in\left\{ 0,1\right\} ,\;\forall t\in\left[0,T\right).\label{eq:initial guess sym}
\end{align}
and compare the results.

\begin{figure}
\begin{centering}
\includegraphics[width=0.98\figwidth]{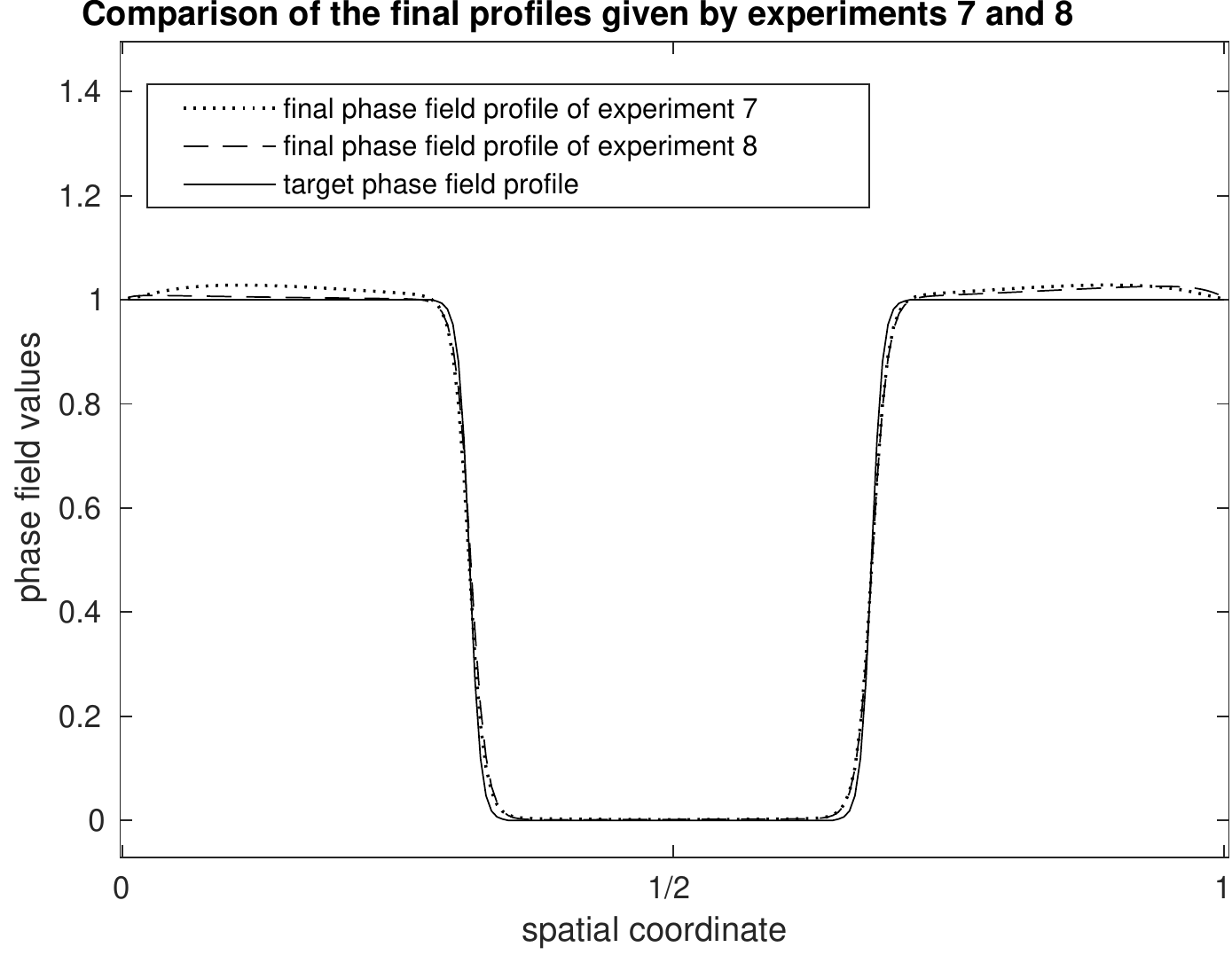}
\par\end{centering}
\caption{\label{fig:The-comparison-of-profiles-7-8}Comparison of the best
and the worst (in terms of the error $\left\Vert \tilde{y}_{h}-P_{h}\tilde{y}_{f}\right\Vert _{2}$)
phase field profile estimations for experiments 4 through 8.}
\end{figure}

The settings of these experiments can be found in Table \ref{tab:The-settings-for-4-8}.
First of all, let us emphasize that all the experiments discussed
below lead to a rather accurate reproduction of the target profile.
This is obvious from the values of the error $\left\Vert \tilde{y}_{h}-P_{h}\tilde{y}_{f}\right\Vert _{2}$
listed in Table \ref{tab:The-settings-for-4-8} as well as from Figure
\ref{fig:The-comparison-of-profiles-7-8} where the best and the worst
final profiles are shown. Also, all the experiments in this section
remain realistic since $y$ is kept within the admissible bounds (\ref{eq:PF physicality condition}).
\begin{table*}
\caption{The settings for experiments 4 through 8 and the respective values
of the difference (error) from the prescribed profile.\label{tab:The-settings-for-4-8}}

\centering{}%
\begin{tabular}{cccccc}
\toprule
 & \multicolumn{5}{c}{Simulation Number}\tabularnewline
Parameter & 4 & 5 & 6 & 7 & 8\tabularnewline
\midrule
number of time steps $N_{t}$ & $10^{5}$ & $10^{5}$ & $10^{5}$ & $10^{5}$ & $10^{5}$\tabularnewline
number of grid points $N_{x}$ & $200$ & $200$ & $200$ & $200$ & $200$\tabularnewline
initial control given by & (\ref{eq:initial guess asy}) & (\ref{eq:initial guess asy}) & (\ref{eq:initial guess asy}) & (\ref{eq:initial guess asy}) & (\ref{eq:initial guess sym})\tabularnewline
$T$ (final time) & $0.05$ & $0.4$ & $0.4$ & $0.4$ & $0.4$\tabularnewline
regularization parameter $\alpha$ & $0$ & $0$ & $5\cdot10^{-10}$ & $10^{-9}$ & $0$\tabularnewline
gradient descent step size $\varepsilon$ & $2\cdot10^{14}$ & $3\cdot10^{13}$ & $10^{13}$ & $10^{13}$ & $3\cdot10^{14}$\tabularnewline
number of iterations & $150$ & $100$ & $100$ & $125$ & $100$\tabularnewline
$\left\Vert \tilde{y}_{h}-P_{h}\tilde{y}_{f}\right\Vert _{2}$ & $0.3359370$ & $0.3392653$ & $0.4725518$ & $0.8469198$ & $0.3293709$\tabularnewline
\bottomrule
\end{tabular}
\end{table*}

\begin{figure*}
\begin{centering}
\includegraphics[width=0.49\textwidth]{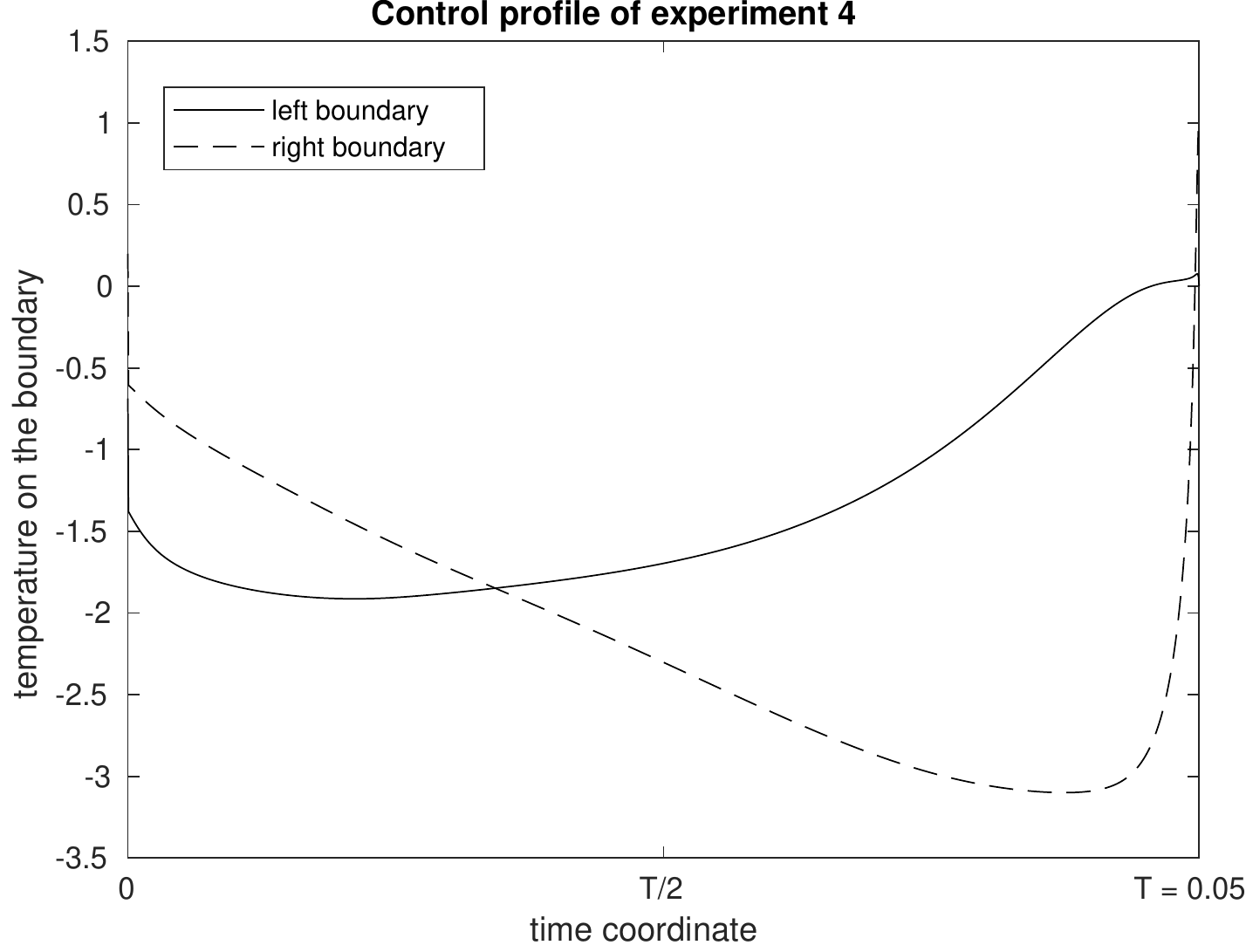}
\includegraphics[width=0.49\textwidth]{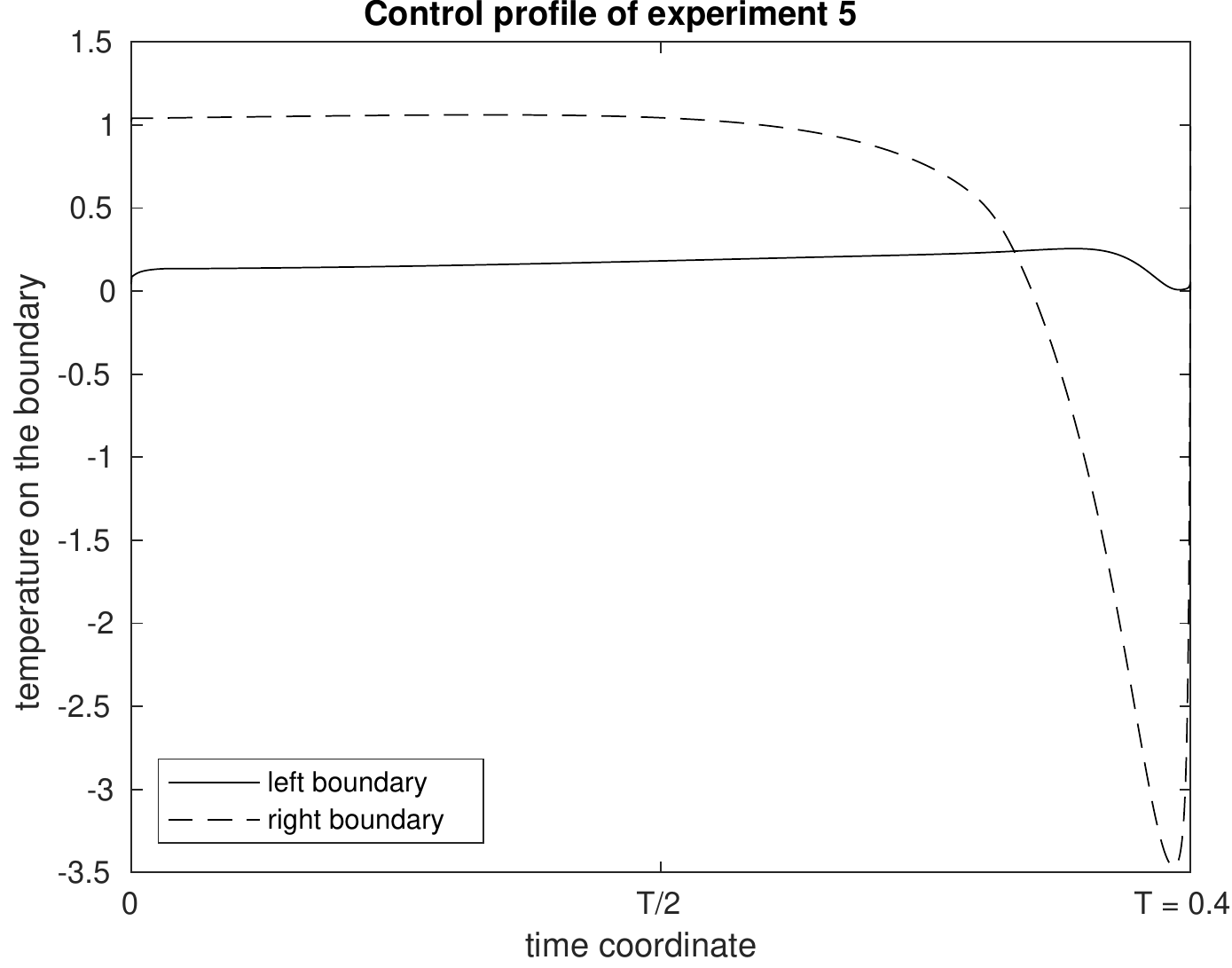}\medskip{}
\par\end{centering}
\begin{centering}
\includegraphics[width=0.49\textwidth]{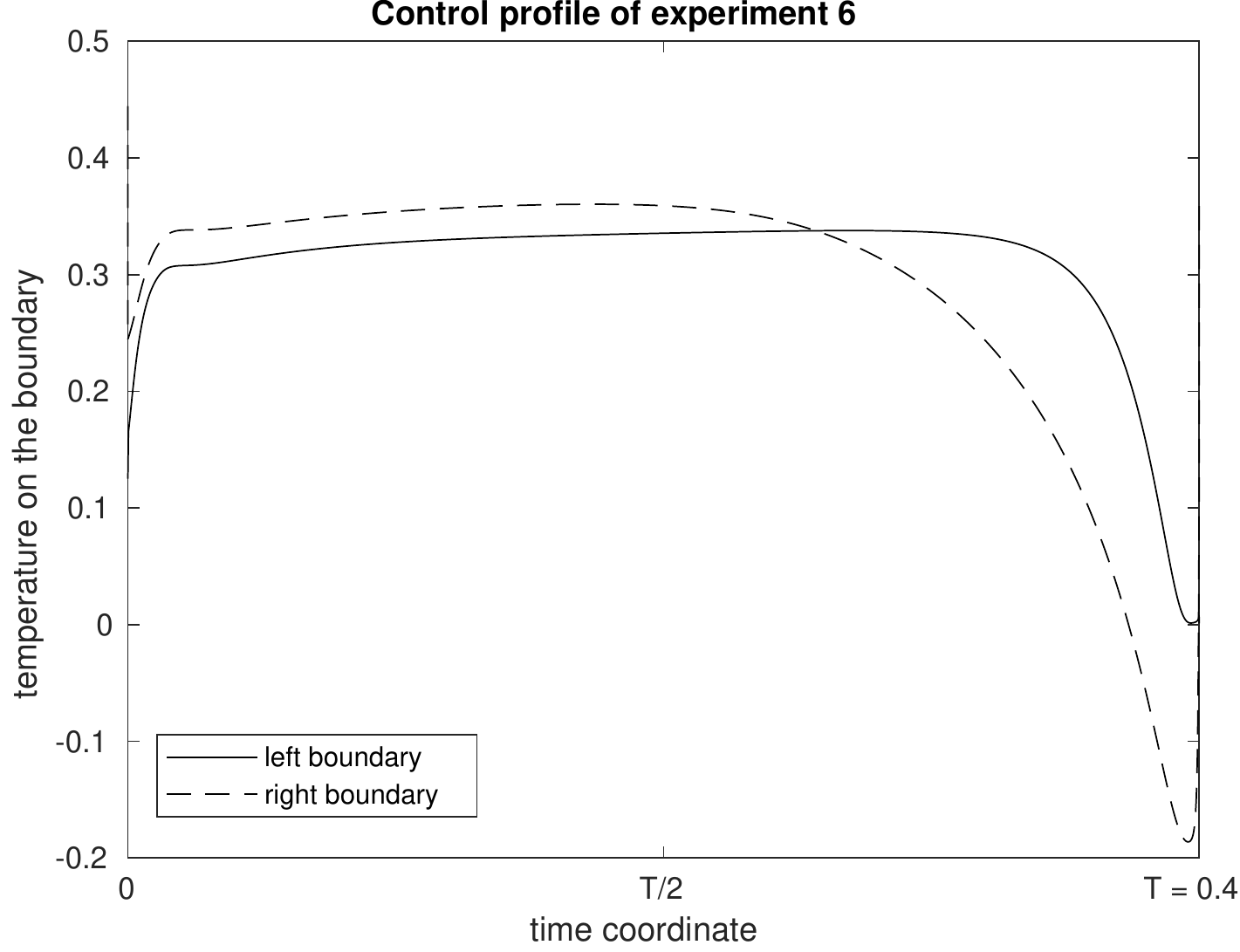}
\includegraphics[width=0.49\textwidth]{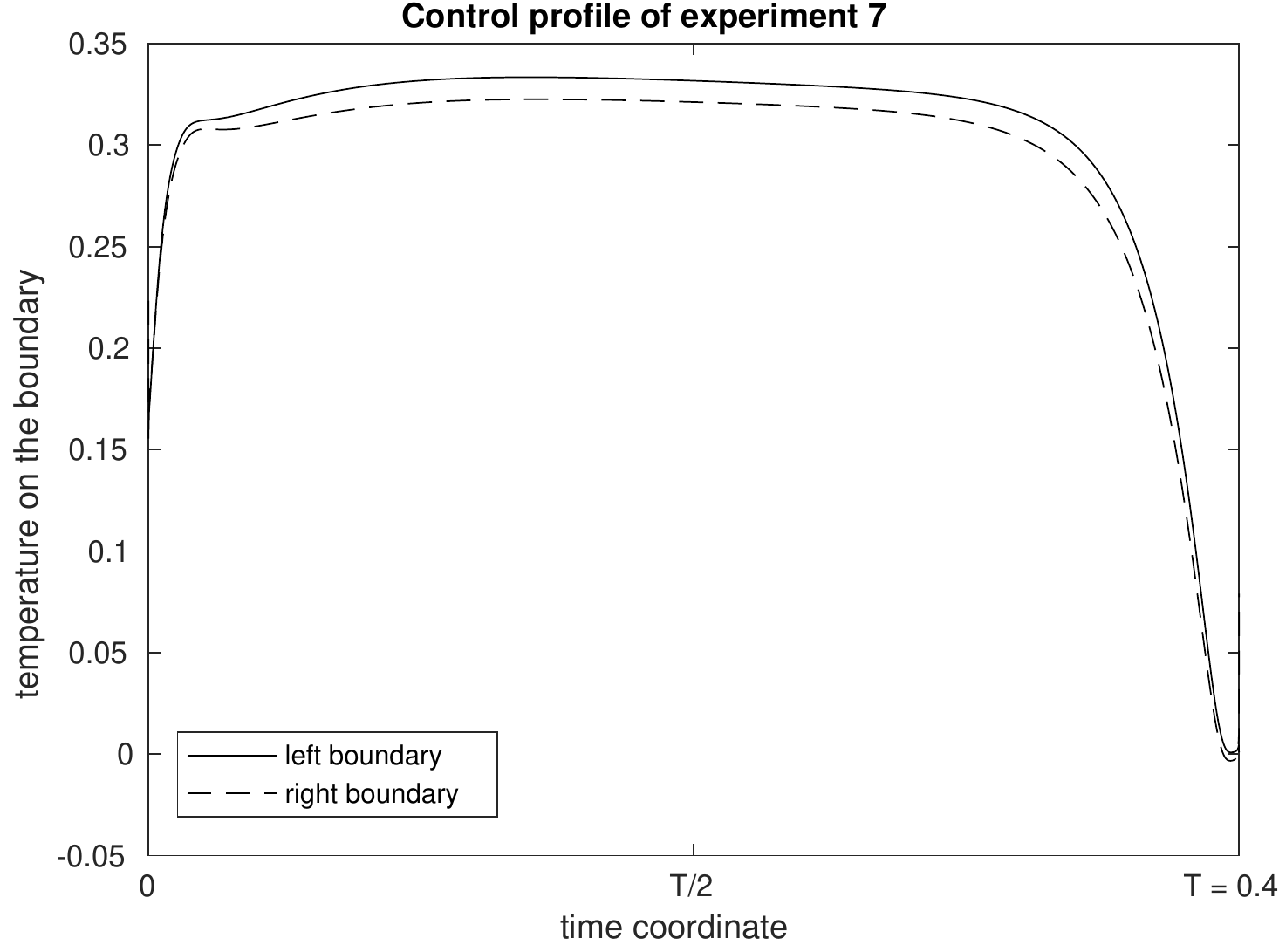}\medskip{}
\par\end{centering}
\begin{centering}
\includegraphics[width=0.49\textwidth]{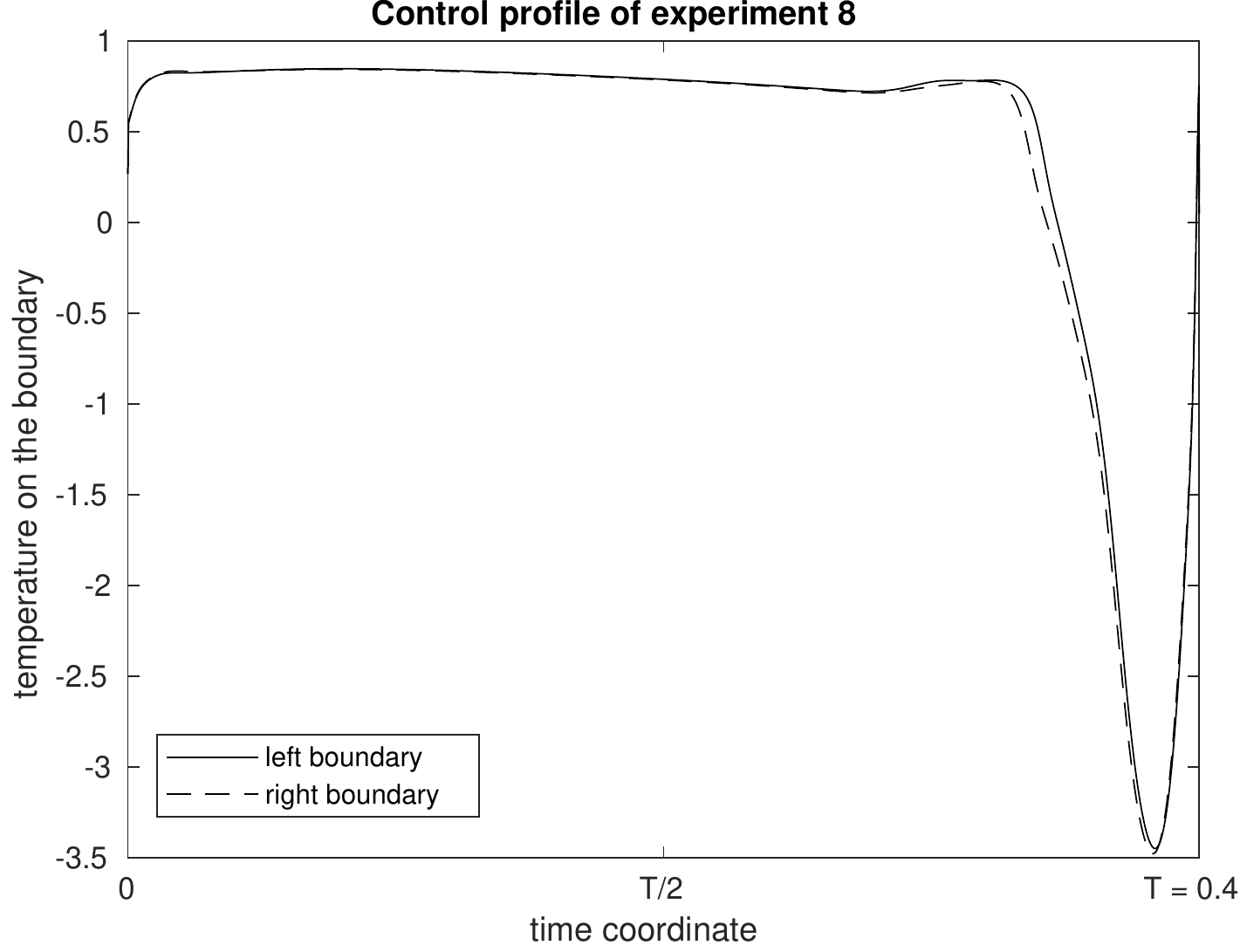}
\par\end{centering}
\caption{\label{fig:Optimized-control-profiles-4-8}Optimized temporal control
profiles of experiments 4 through 8. Experiments 4 and 5 show how
an optimization without regularization responds to an increase of
final time $T$. Regularization is then added in experiments 6 and
7. The control profiles flatten and the control becomes more evenly
distributed. Lastly, experiment 8 shows a different non-regularized
control given by the symmetric initial condition $\tilde{u}_{\text{ini}}$.}

\end{figure*}

Let us now focus on the differences in the obtained temporal control
profiles $u_{h}$ between the individual experiments (see Figure \ref{fig:Optimized-control-profiles-4-8}).
\begin{itemize}
\item Experiments 4 and 5 demonstrate how the choice of very different final
times $T$ ($T=0.05$ and $T=0.4$) affects the control $u_{h}$.
When the control is given more time (experiment 5, $T=0.4$), its
action is delayed to the final part of the interval $\left(0,T\right)$.
Nevertheless, the range of $u_{h}$ is similar in both cases.
\item Experiments 5, 6, and 7 show the effects of regularization $\alpha$
of different amplitudes with the fixed final time $T=0.4$. As expected,
increasing $\alpha$ significantly reduces the range of the control
$u_{h}$. However, this does not affect the overall minimization quality
substantially. This is because the action of the control is distributed
over a longer time period.
\item Experiment 5 and 8 show the effect of the initial guess $u_{0}$.
It is not surprising that in experiment 8, a spatially symmetric problem
with a symmetric initial guess (\ref{eq:initial guess sym}) result
in (almost) identical left and right temporal control profiles $u_{h}\left(\cdot,0\right)$
and $u_{h}\left(\cdot,1\right)$. Experiment 5 differs from experiment
8 only in the non-symmetry of the initial guess (\ref{eq:initial guess asy}).
As a consequence, the left and right temporal control profiles are
completely different.
\item Experiments 5, 6, and 7 also indicate that the regularization partially
symmetrizes the left and the right temporal control profiles.
\end{itemize}

\subsubsection{\label{subsec:Moving-a-Gap}Moving a Gap Between Crystals}

\begin{figure}
\begin{centering}
\includegraphics[width=0.98\figwidth]{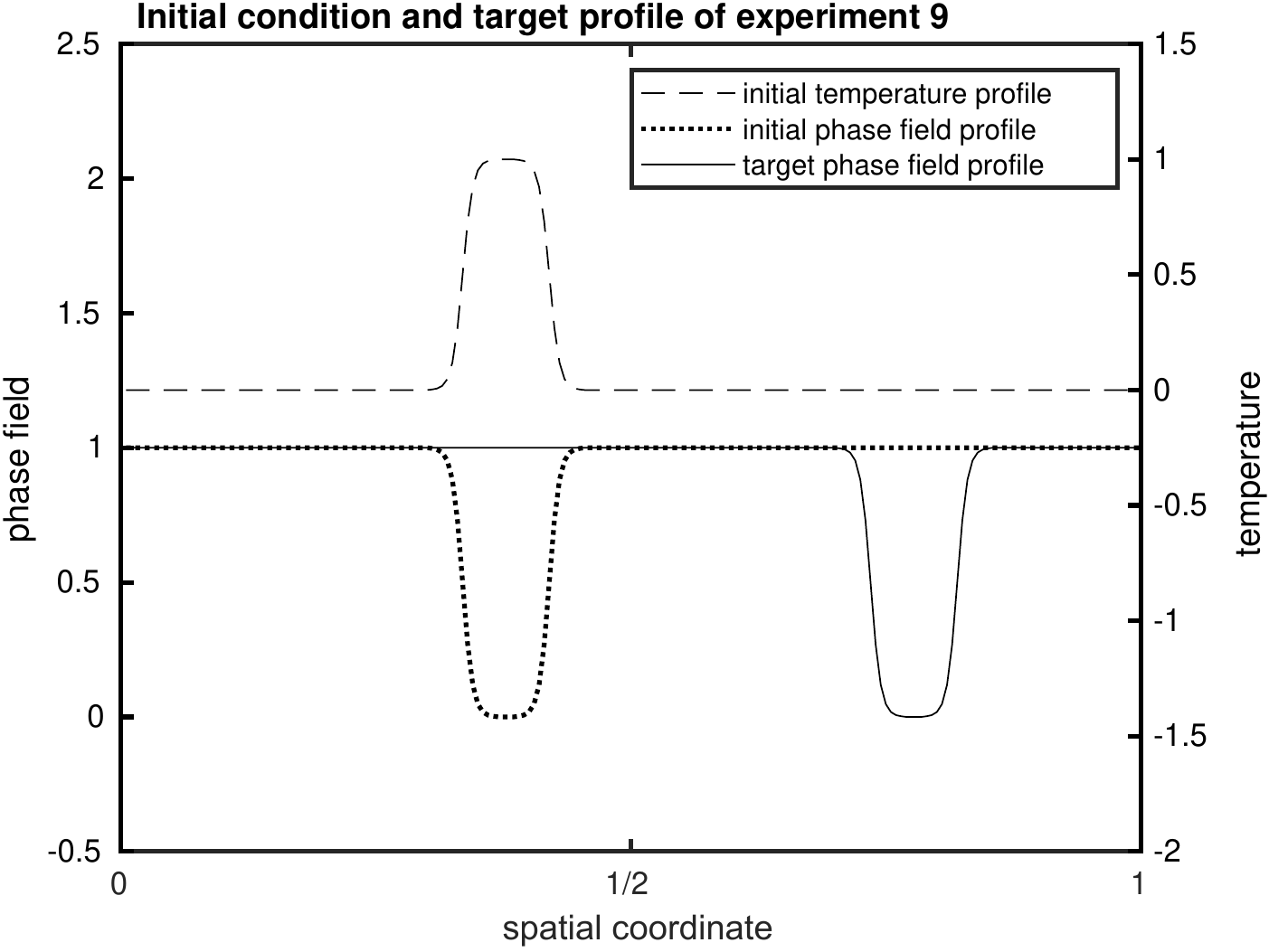}
\par\end{centering}
\caption{\label{fig:Initial-setting-for-ex9}The initial temperature and phase
field profiles $\tilde{y}_{\text{ini}},y_{\text{ini}}$ along with
the target profile $\tilde{y}_{f}$ for experiment 9. The values of
the boundary condition $\tilde{y}_{\text{bc}}$ are given by (\ref{eq:experiments-4-8-PF-BCs}).}
\end{figure}

Experiment 9 showcases how even highly non-trivial control can be
obtained. Consider a situation in which two crystals occupy the whole
domain except for a comparatively small gap between them. We aim to
move the liquid gap to a different position in the domain. This is
reflected by the settings of $y_{\text{ini}}$, $\tilde{y}_{\text{ini}}$
and the target profile $\tilde{y}_{f}$ shown in Figure \ref{fig:Initial-setting-for-ex9}
as well as the use of the boundary condition (\ref{eq:experiments-4-8-PF-BCs}).
The full setup of the experiment is summarized in Table \ref{tab:The-settings-of-ex9}.
Note that during optimization, the gradient descent step size was
(manually) adjusted from $\varepsilon_{1}$ to $\varepsilon_{2}$
to ensure convergence.

\begin{figure}
\begin{centering}
\includegraphics[width=0.98\figwidth]{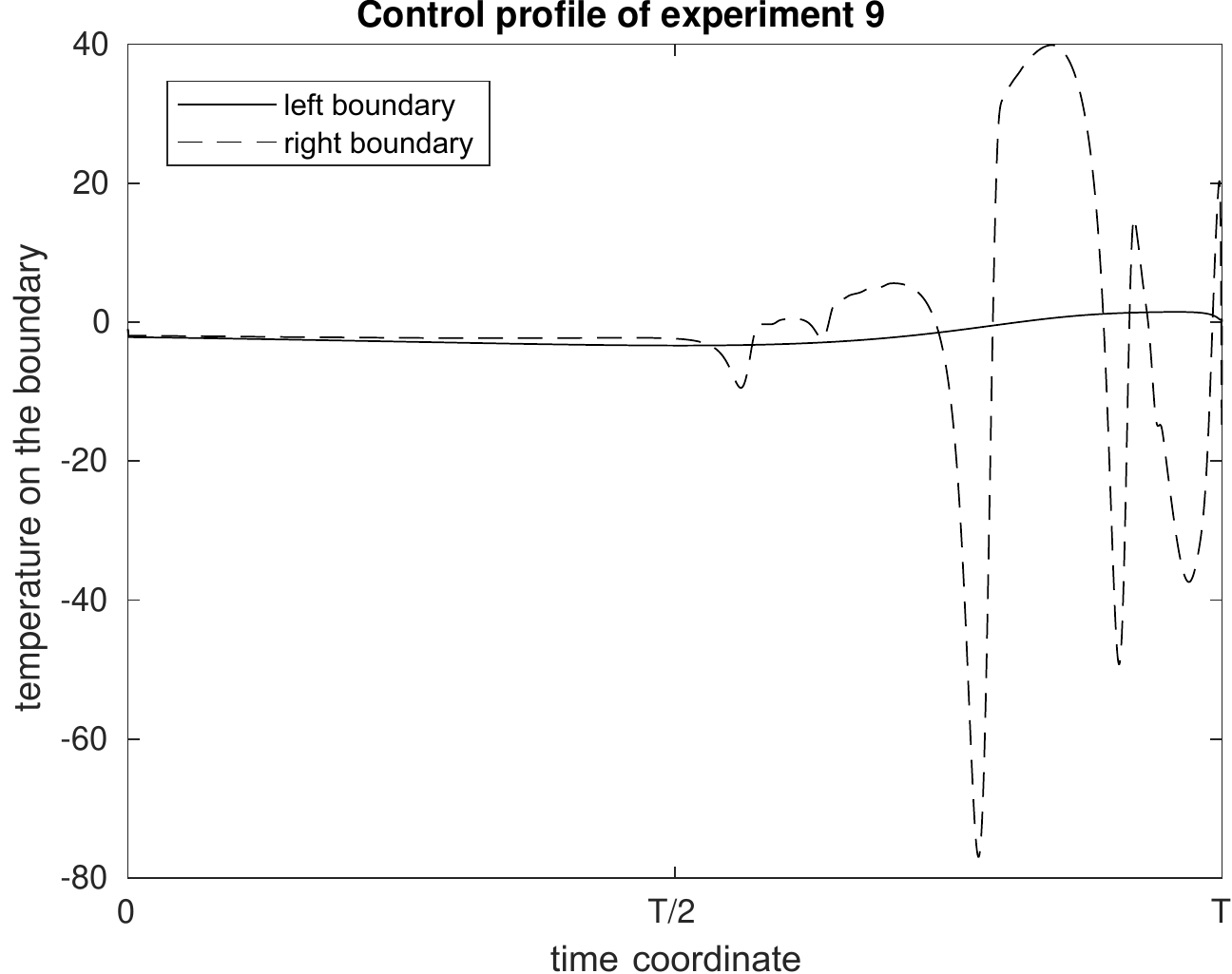}
\par\end{centering}
\caption{\label{fig:Control-profiles-for-ex9}Temporal control profiles for
experiment 9.}

\end{figure}

\begin{figure}
\begin{centering}
\includegraphics[width=0.98\figwidth]{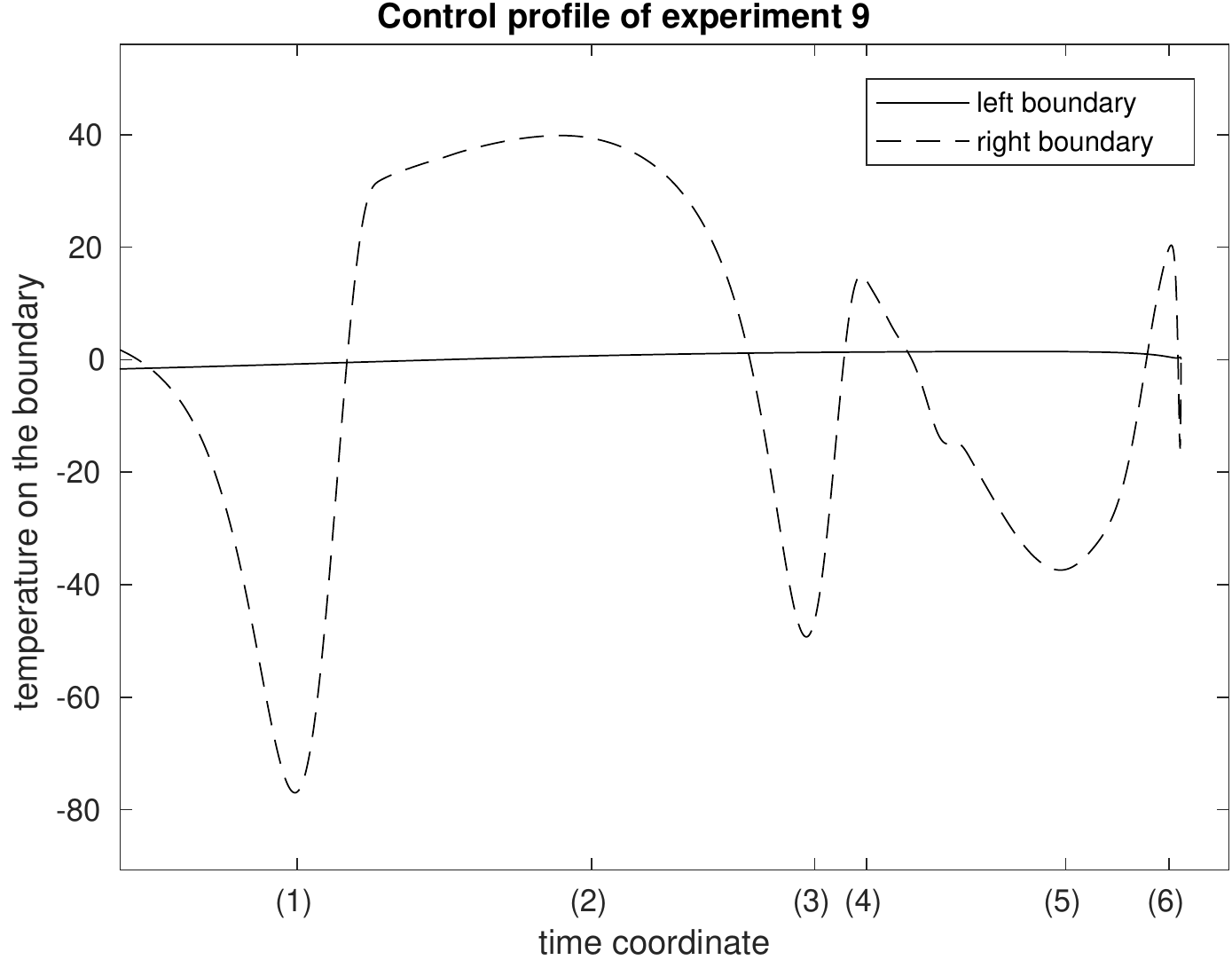}
\par\end{centering}
\caption{\label{fig:time evolution of ex 9}A detailed look at the latter part
of the temporal control profile in experiment 9.}
\end{figure}

\begin{figure*}
\begin{centering}
\includegraphics[width=0.49\textwidth]{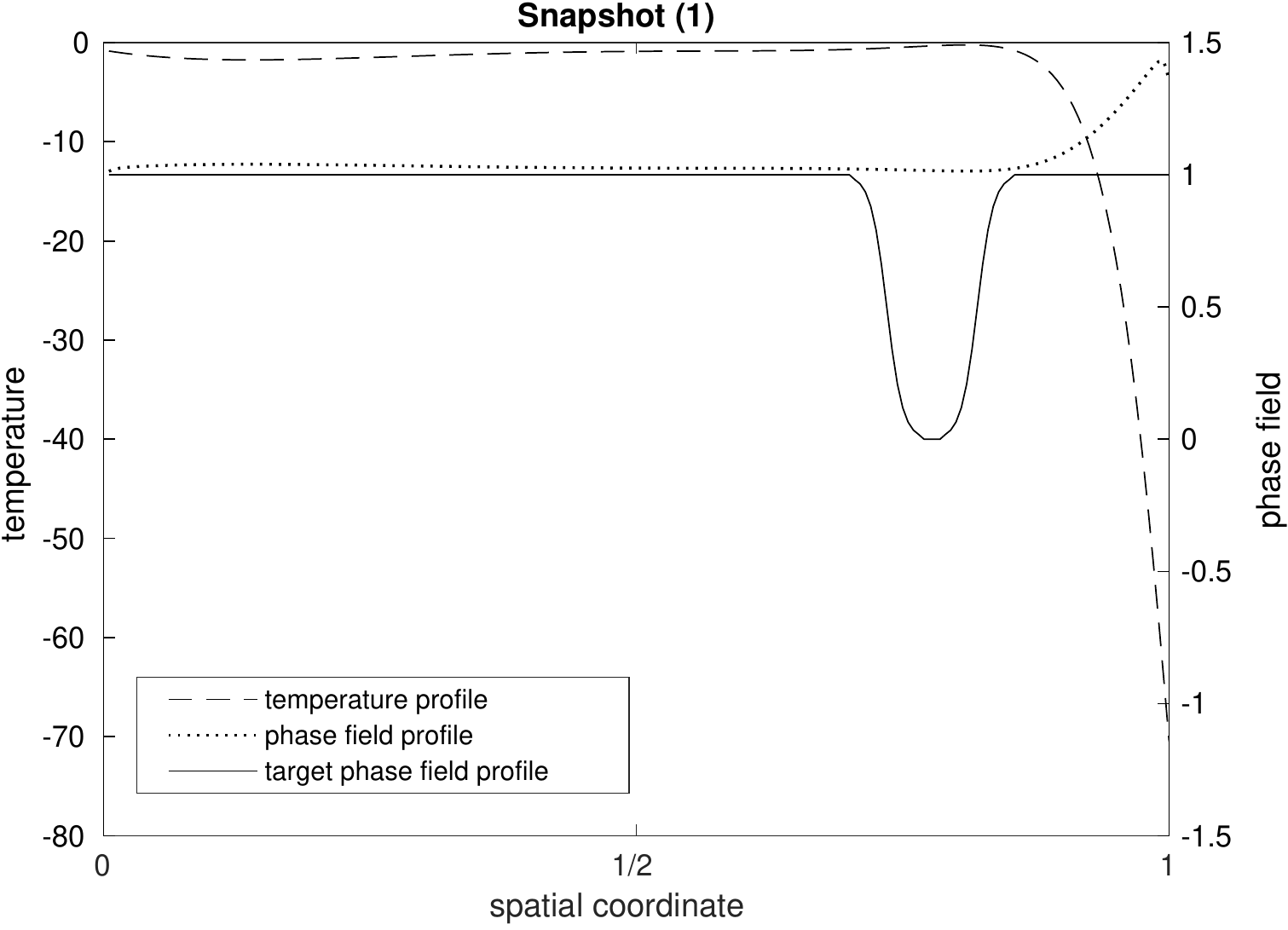}
\includegraphics[width=0.49\textwidth]{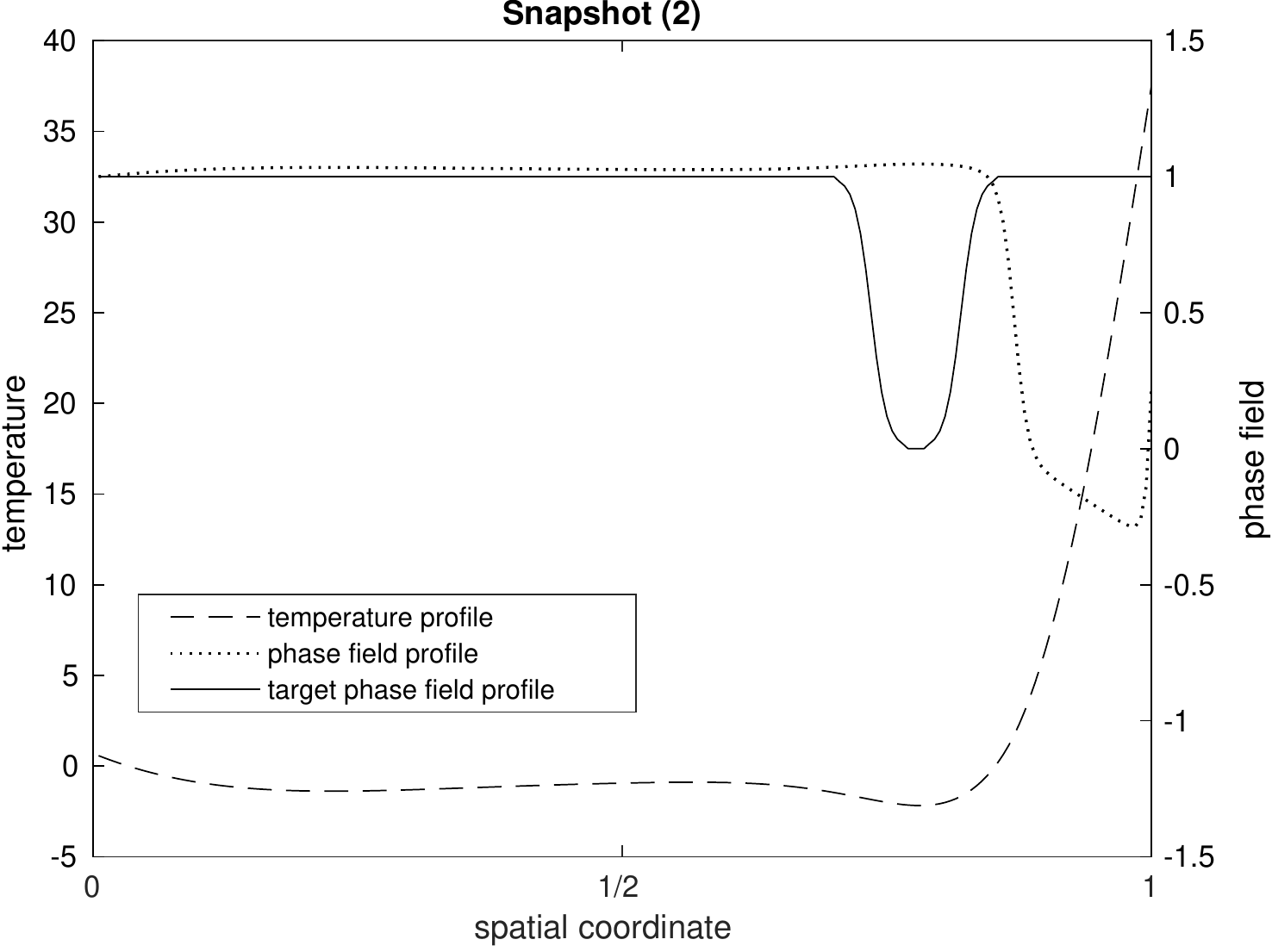}\medskip{}
\par\end{centering}
\begin{centering}
\includegraphics[width=0.49\textwidth]{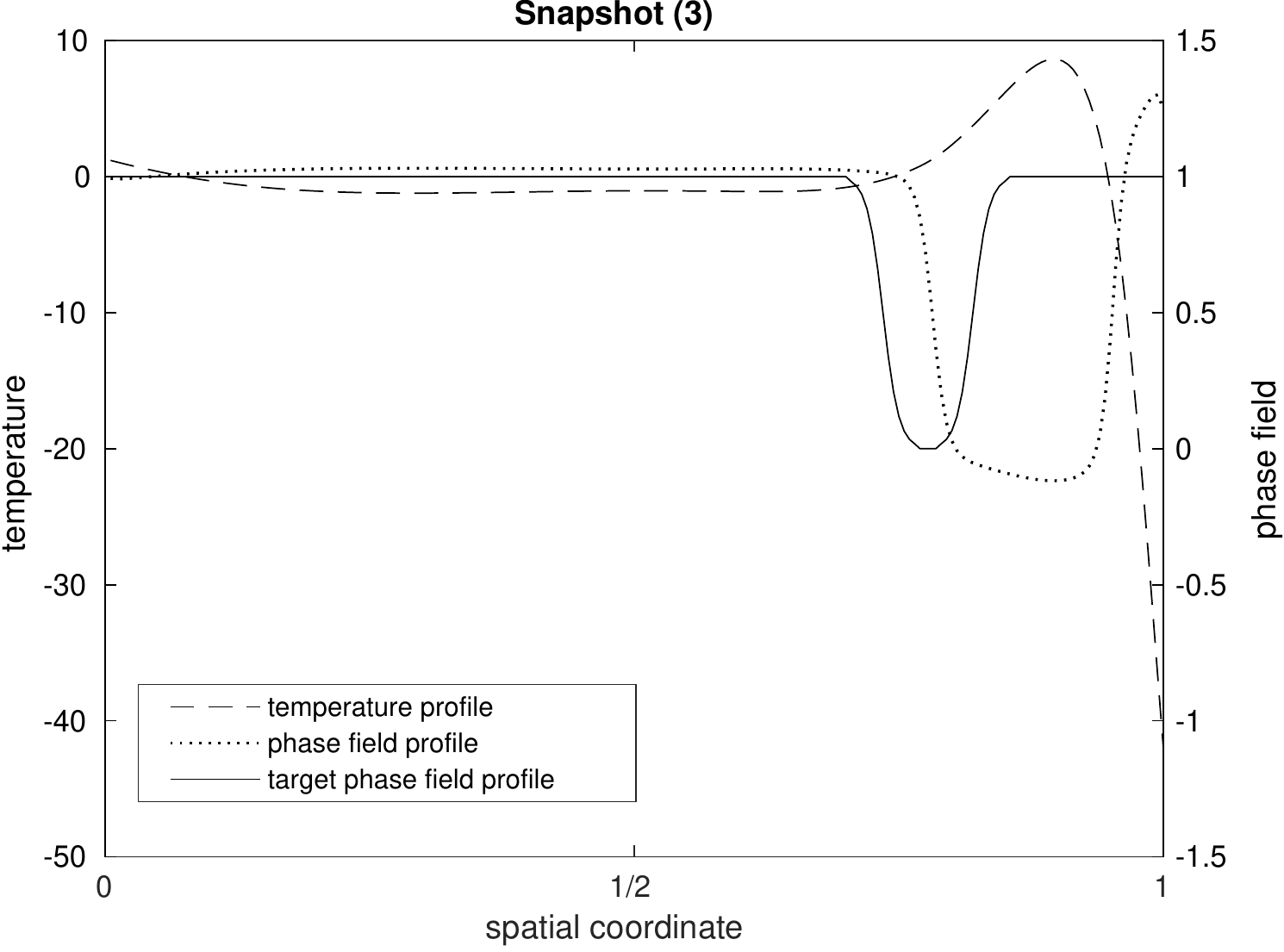}
\includegraphics[width=0.49\textwidth]{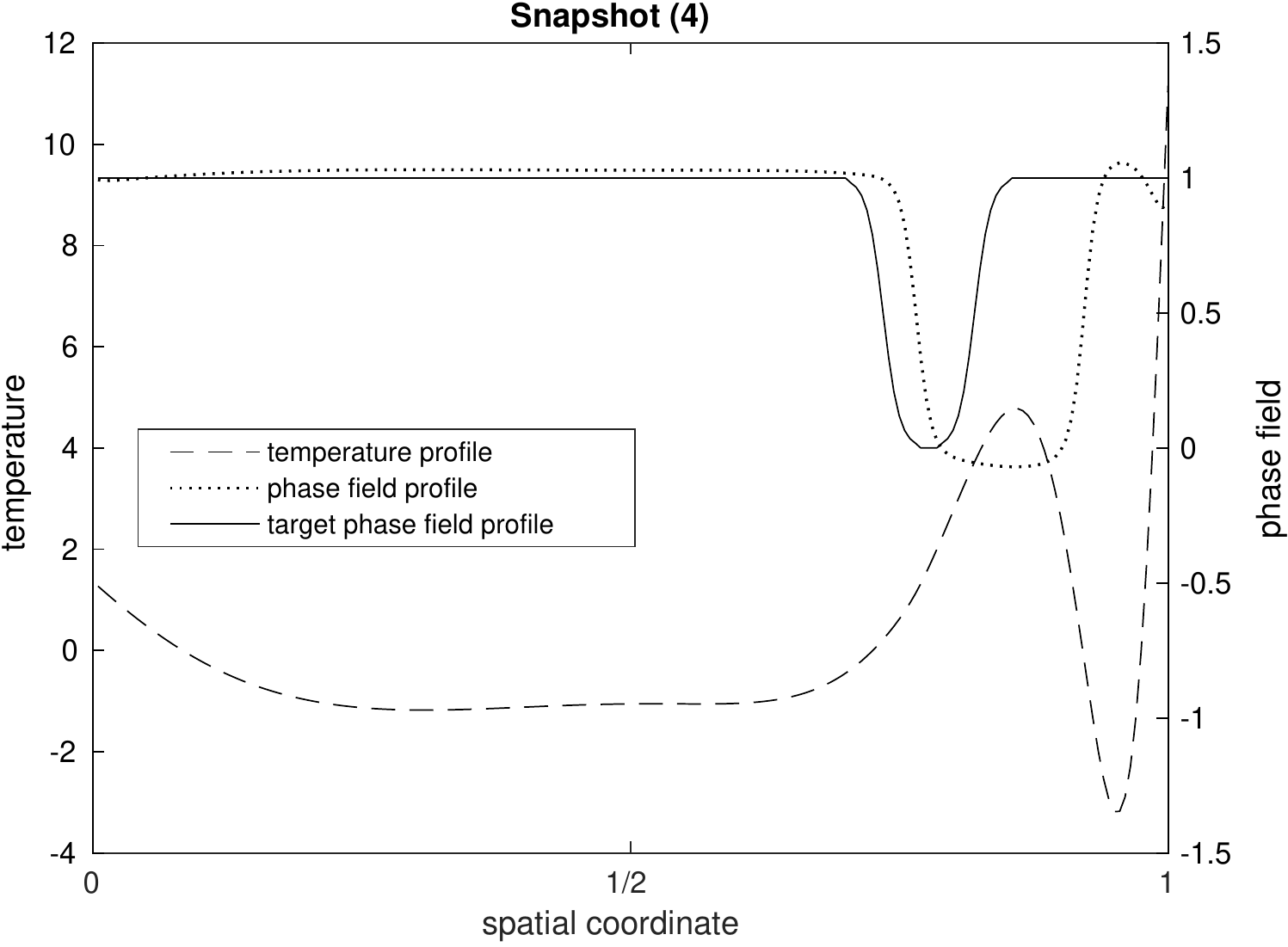}\medskip{}
\par\end{centering}
\begin{centering}
\includegraphics[width=0.49\textwidth]{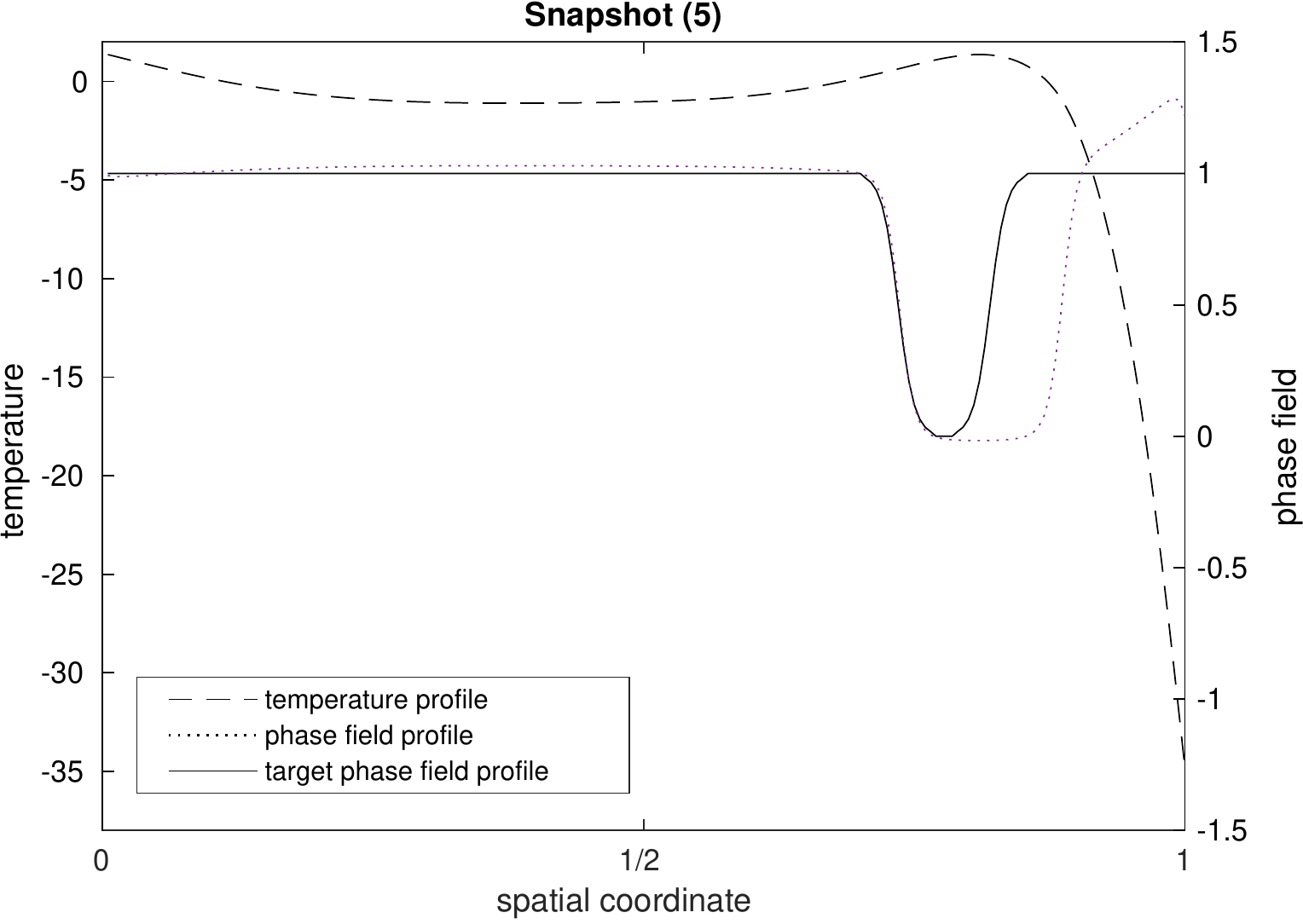}
\includegraphics[width=0.49\textwidth]{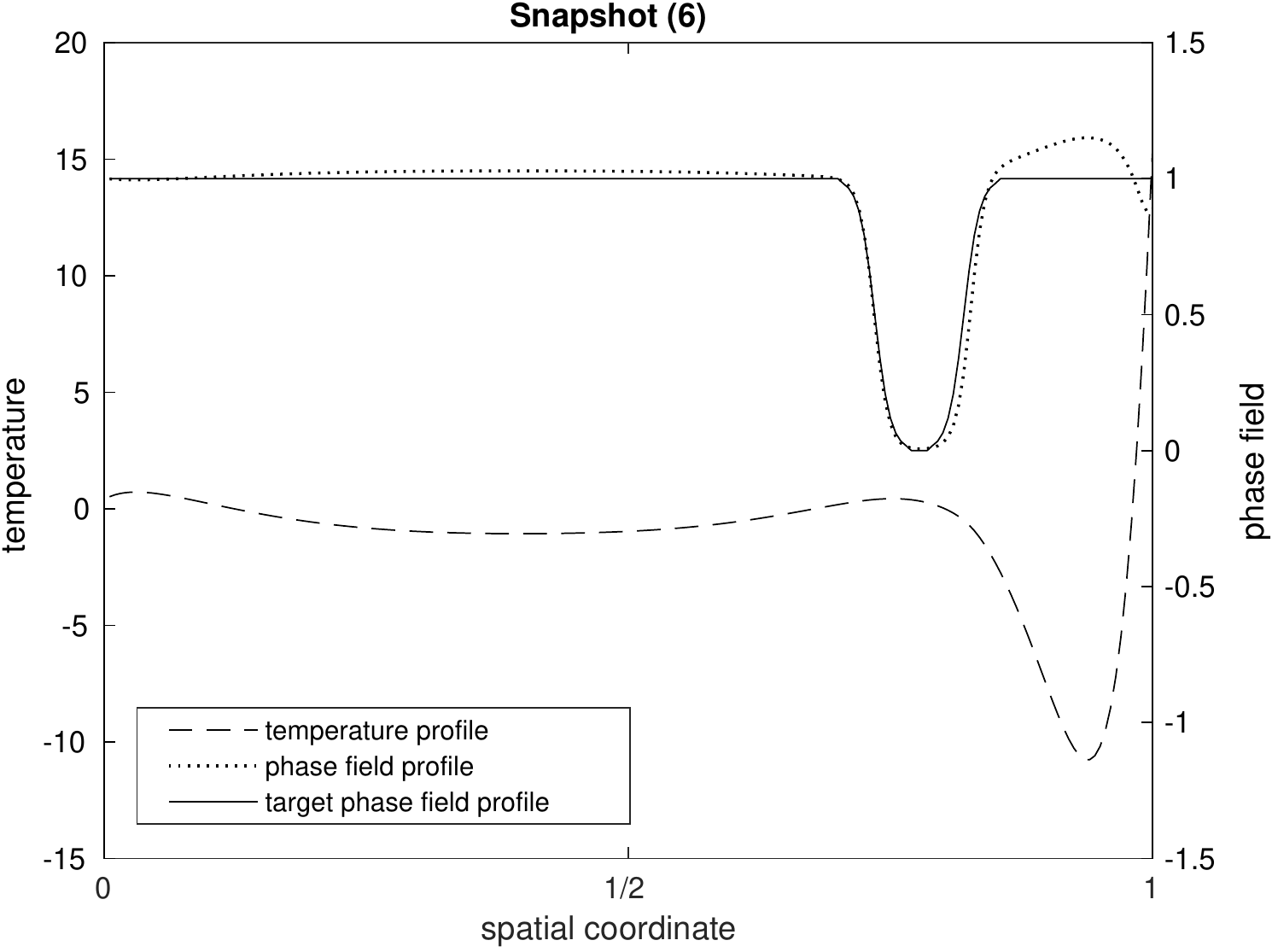}
\par\end{centering}
\caption{\label{fig:time-evolution of ex9-detail}The spatial profiles of the
phase field and temperature at the significant times (1)--(6) marked
in Figure \ref{fig:time evolution of ex 9}.}

\end{figure*}

The target profile $\tilde{y}_{f}$ is achieved by a rather surprising
process. First, the entire domain is solidified. Afterward, the molten
region is recreated next to the right boundary. As the melting proceeds
to the left, the right boundary is undercooled again to create another
solid subdomain on the right. The resulting control profiles for the
entire time interval are displayed in Figure \ref{fig:Control-profiles-for-ex9}.
In addition, six snapshots of the final phase of the solution evolution
are depicted in Figure \ref{fig:time-evolution of ex9-detail}, with
the corresponding times marked in Figure \ref{fig:time evolution of ex 9}.
In terms of the condition (\ref{eq:PF physicality condition}), the
simulation is not realistic.
\begin{table}
\caption{\label{tab:The-settings-of-ex9}Settings for experiment 9 and the
respective value of the difference (error) from the prescribed profile.}

\centering{}%
\begin{tabular}{cc}
\toprule
Parameter & Value\tabularnewline
\midrule
number of time steps $N_{t}$ & $10^{5}$\tabularnewline
number of grid points $N_{x}$ & $200$\tabularnewline
initial control given by & (\ref{eq:initial guess sym})\tabularnewline
$T$ (final time) & 0.1\tabularnewline
regularization parameter $\alpha$ & 0\tabularnewline
gradient descent step size $\varepsilon$ & $\varepsilon_{1}=10^{16},\,\varepsilon_{2}=5\cdot10^{15}$\tabularnewline
number of iterations & $225$ with step $\varepsilon_{1}$, $25$ with $\varepsilon_{2}$\tabularnewline
$\left\Vert \tilde{y}_{h}-P_{h}\tilde{y}_{f}\right\Vert _{2}$ & 0.5932381\tabularnewline
\bottomrule
\end{tabular}
\end{table}

\subsection{\label{subsec:PhaseField-Num-1}Dirichlet Boundary Condition Control
for the Phase Field Problem in 2D}

In the preceding Section, a bevy of examples detailed the utility
and the possible shortcomings that come with using a linear reaction
term (\ref{eq:reaction-term}) in (\ref{eq:PF-functional})-(\ref{eq:PF-PF-IC}).
Among other things, the effects of regularization, changes in final
time and initial guess for the control were discussed in detail. In
this section, simulations in two spatial dimensions are performed.
The focus will no longer be on tweaking the parameters of the simulations.
Instead, we focus on the effect of using the reaction term (\ref{eq:sigma limiter reaction})
in (\ref{eq:PF-functional})-(\ref{eq:PF-PF-IC}) as well as showing
non-trivial realistic controls that arise in some situations.

Parameters that are common to all the experiments detailed in this
section are listed in Table \ref{tab:Settings-for-the-physical-params-PF-2D}.
\begin{table}
\caption{\label{tab:Settings-for-the-physical-params-PF-2D}Parameter settings
for the phase field simulations in Section \ref{subsec:PhaseField-Num-1}.}

\centering{}%
\begin{tabular}{ccl}
\toprule
Param. & Value & Physical Meaning\tabularnewline
\midrule
$\gamma$ & 3.0 & coefficient of attachment kinetics\tabularnewline
$\beta$ & 300 & dimensionless representation of supercooling\tabularnewline
$\xi$ & 0.0101 & interface thickness scaling\tabularnewline
$y_{\text{mt}}$ & 1.0 & melting temperature\tabularnewline
$H$ & 2.0 & latent heat\tabularnewline
$\varepsilon_{0}$ & 0 & parameter of the sigmoid function in (\ref{eq:sigma limiter reaction})\tabularnewline
$\varepsilon_{1}$ & 0.2 & parameter of the sigmoid function in (\ref{eq:sigma limiter reaction})\tabularnewline
$L_{x_{1}}$ & 0.6 & spatial dimension in the $x_{1}$ direction\tabularnewline
$L_{x_{2}}$ & 1.0 & spatial dimension in the $x_{2}$ direction\tabularnewline
\bottomrule
\end{tabular}
\end{table}

\subsubsection{\label{subsec:Moving-a-Crystal-n-t-s-lin}Moving a Crystal from North
to South with Different Reaction Terms}

The aim of this set of numerical experiments is to move a crystal
from one position in the domain to another, while maintaining its
shape and size. Two experiments are performed. One using the linear
reaction term (\ref{eq:reaction-term}) and the other uses the alternative
more advanced reaction term (\ref{eq:sigma limiter reaction}).

\begin{figure}
\begin{minipage}[t]{0.47\columnwidth}%
\begin{center}
\includegraphics[width=1\textwidth]{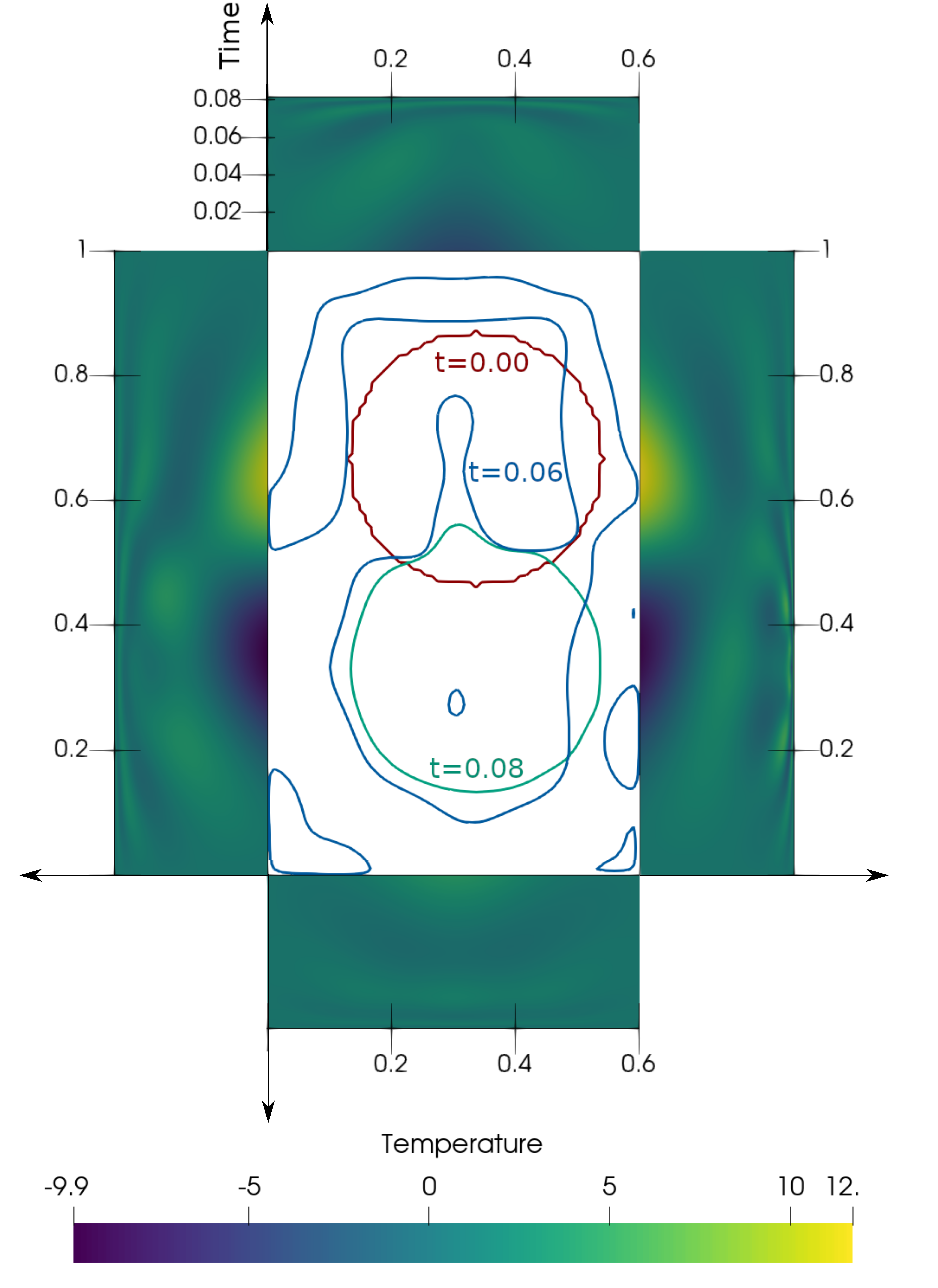}\caption{\label{fig:2d_test4_lin}Simulation with the linear reaction term
(\ref{eq:reaction-term}). For details, see Section \ref{subsec:Moving-a-Crystal-n-t-s-lin}.}
\par\end{center}%
\end{minipage}\hfill{}%
\begin{minipage}[t]{0.47\columnwidth}%
\begin{center}
\includegraphics[width=1\textwidth]{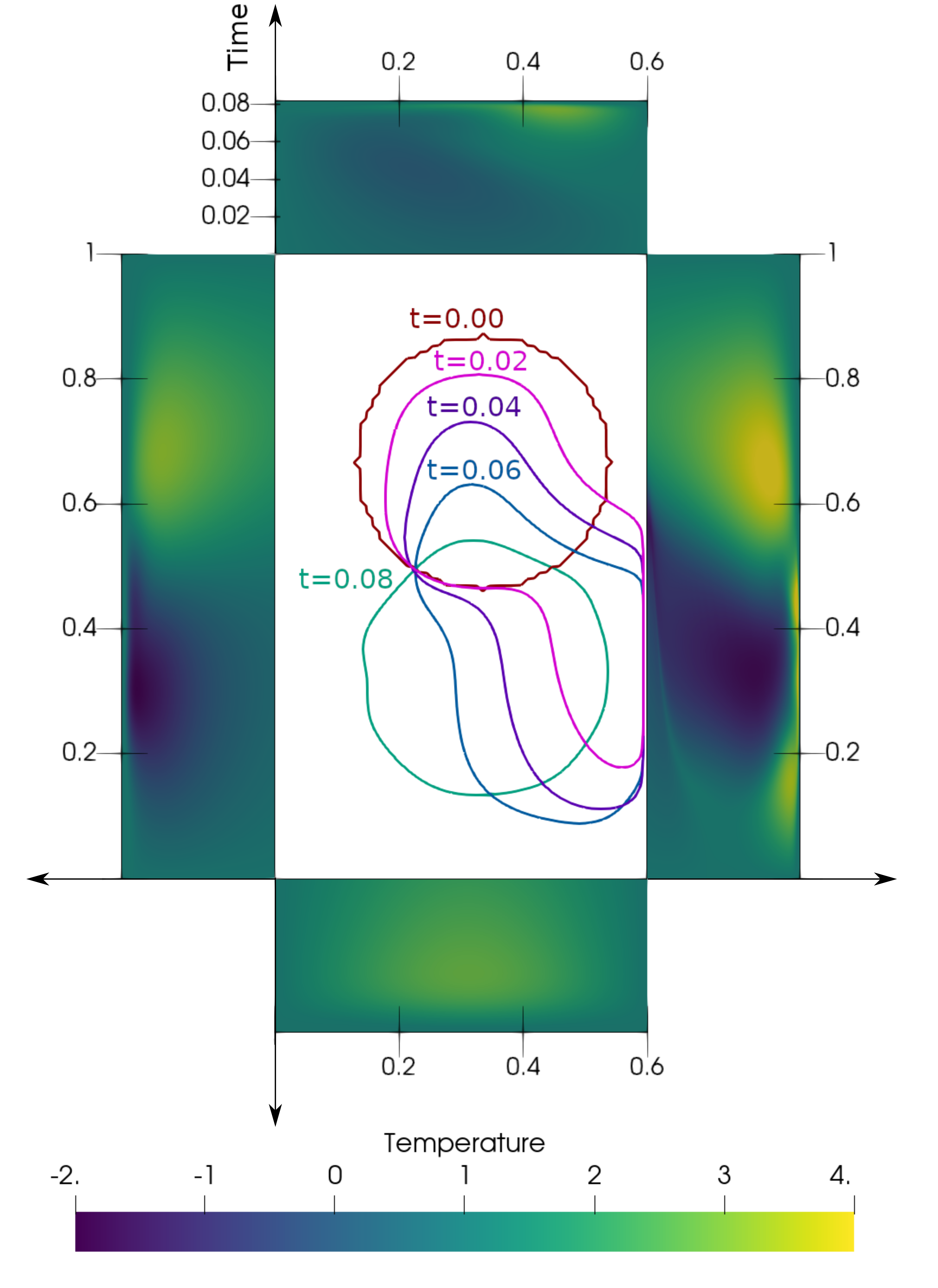}\caption{\label{fig:2d_test4_nonlin} Simulation with the alternative reaction
term (\ref{eq:sigma limiter reaction}). For details, see Section
\ref{subsec:Moving-a-Crystal-n-t-s-lin}.}
\par\end{center}%
\end{minipage}
\end{figure}

For both of the experiments, the initial condition $\tilde{y}_{\text{ini}}$
along with the target profile $\tilde{y}_{f}$ can be found in Figure
\ref{fig:2d_test4_lin}. The boundary condition for the phase field
$\tilde{y}$ is given by

\begin{equation}
\tilde{y}_{\text{bc}}\left(t,x\right)=0,\;\forall x\in\partial\Omega,\;\forall t\in\left[0,T\right).\label{eq:2d_test4_bc}
\end{equation}
The initial guess for the Dirichlet control on the boundary reads

\begin{equation}
u_{0}\left(t,x\right)=1,\;\forall x\in\partial\Omega,\;\forall t\in\left[0,T\right).\label{eq:2d_test4_init}
\end{equation}
Other data for the experiments (some of which is not common for both
experiments), like the mesh resolution or the final error, can be
found in Table \ref{tab:Settings-for-the-physical-params-PF-2D} and
Table \ref{tab:2d_num_param}.

\begin{figure}
\centering{}\includegraphics[width=0.6\figwidth]{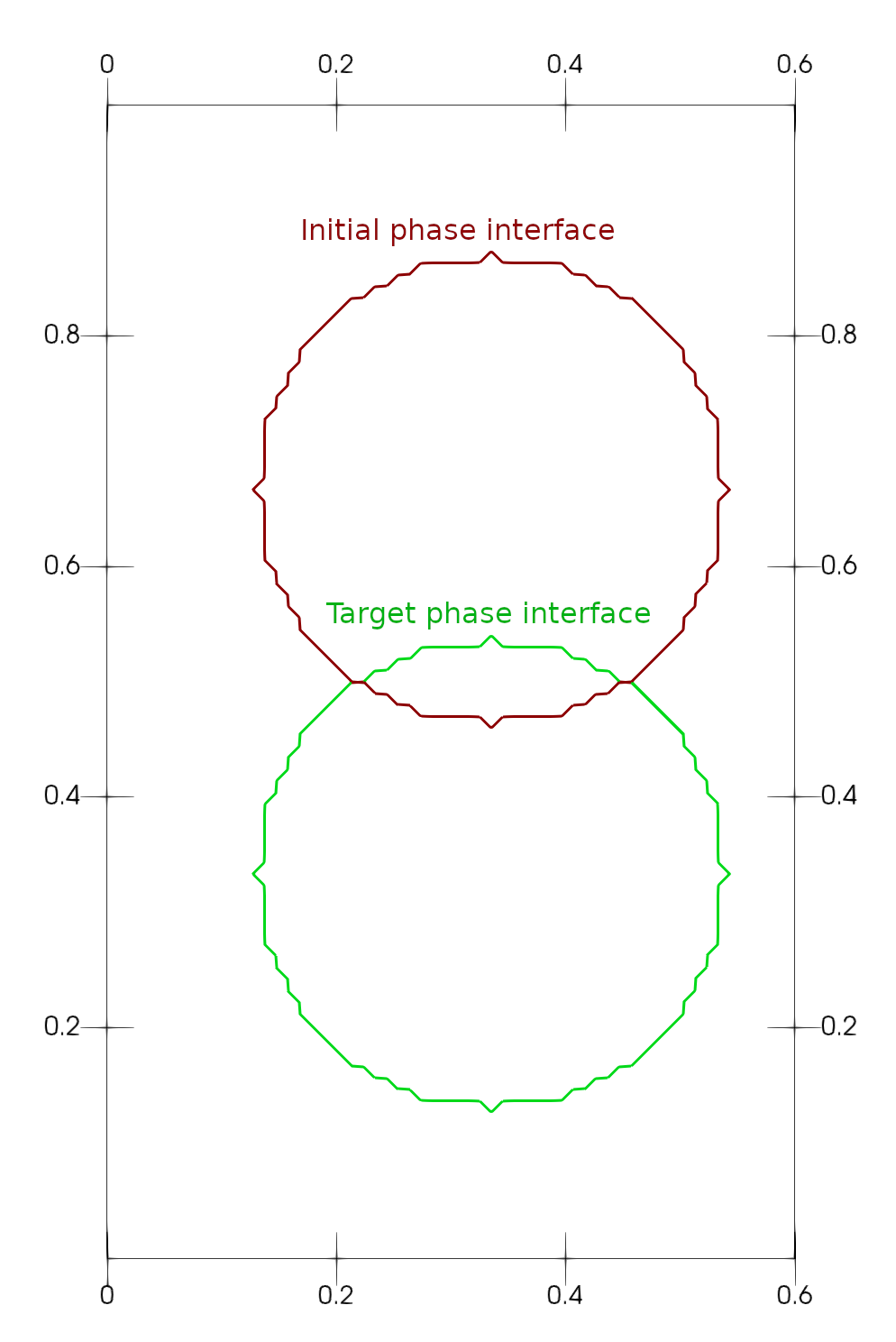}\caption{\label{fig:2d_test4_init}Initial and target phase field interface
$\Gamma$ for the simulation in Section \ref{subsec:Moving-a-Crystal-n-t-s-lin}.}
\end{figure}

The resulting temporal control profiles along with the time evolution
of the level set $\Gamma\left(t\right)$ (the shape of the crystal)
for both experiments can be reviewed in Figures \ref{fig:2d_test4_lin},
for the linear reaction term (\ref{eq:reaction-term}), and Figures
\ref{fig:2d_test4_nonlin} for the alternative term (\ref{eq:sigma limiter reaction}).
Comparing the aforementioned figures shows that the two controls obtained
are qualitatively and quantitatively different.

This is unsurprising, since the state equation (\ref{eq:PF-functional})-(\ref{eq:PF-PF-IC})
with (\ref{eq:reaction-term}) allows for spontaneous nucleation when
the bound (\ref{eq:PF physicality condition}) is violated \cite{PF-Focusing-Latent-Heat}.
It has been shown that the term (\ref{eq:sigma limiter reaction})
does not suffer from such a deficiency \cite{PF-Focusing-Latent-Heat}
and the obtained control reflects this. As a result, only the control
depicted in Figure \ref{fig:2d_test4_nonlin} can be interpreted as
solidification.

\begin{figure}
\begin{centering}
\begin{minipage}[t]{0.32\columnwidth}%
\begin{center}
\includegraphics[viewport=1000bp 0bp 2000bp 1488bp,clip,width=1\textwidth]{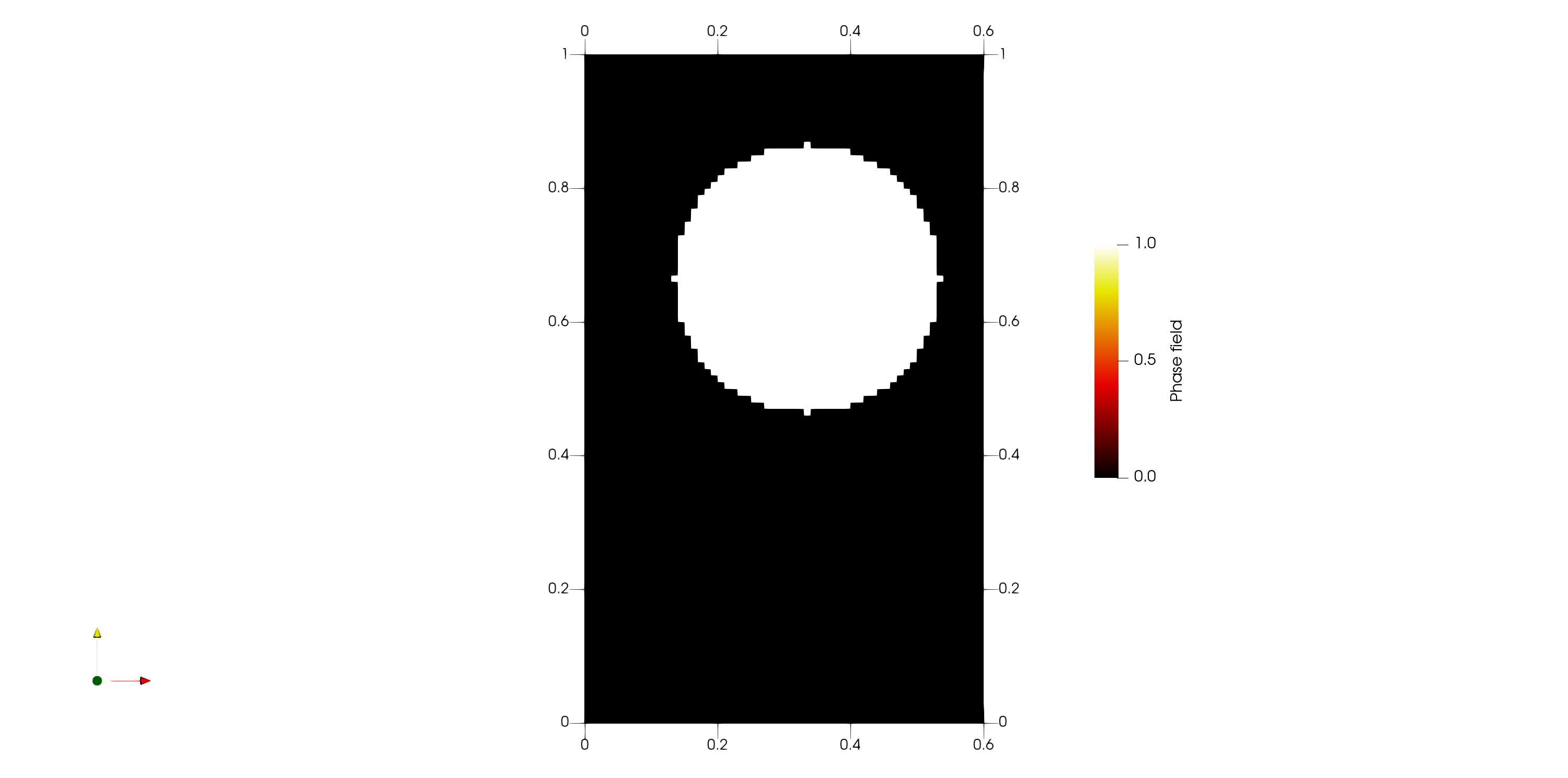}
\par\end{center}
\vspace{-0.8cm}

\begin{center}
$\tilde{y}(0)$
\par\end{center}%
\end{minipage}\hfill{}%
\begin{minipage}[t]{0.32\columnwidth}%
\begin{center}
\includegraphics[viewport=1000bp 0bp 2000bp 1488bp,clip,width=1\textwidth]{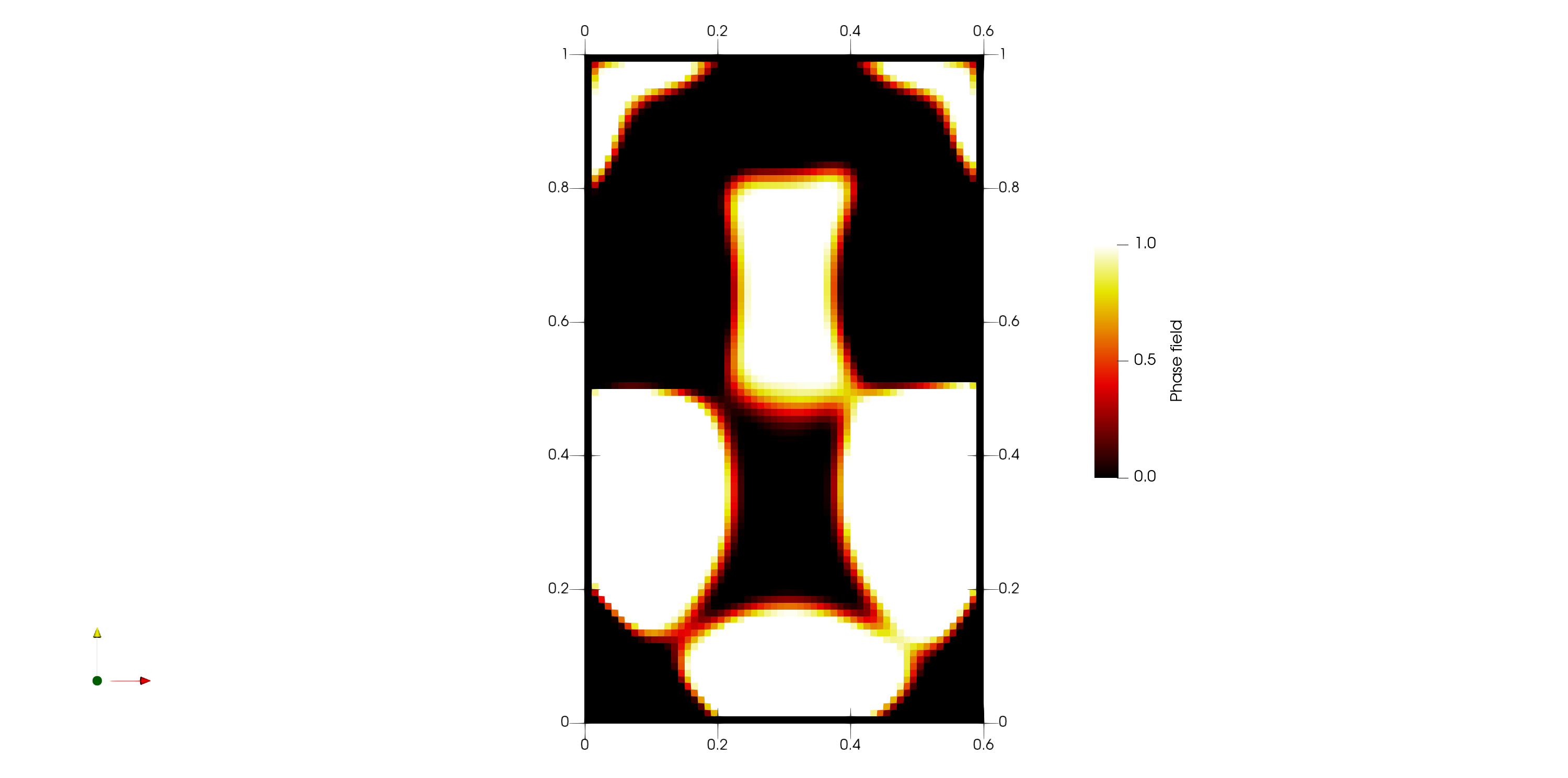}
\par\end{center}
\vspace{-0.8cm}

\begin{center}
$\tilde{y}(0.02)$
\par\end{center}%
\end{minipage}\hfill{}%
\begin{minipage}[t]{0.32\columnwidth}%
\begin{center}
\includegraphics[viewport=1000bp 0bp 2000bp 1488bp,clip,width=1\textwidth]{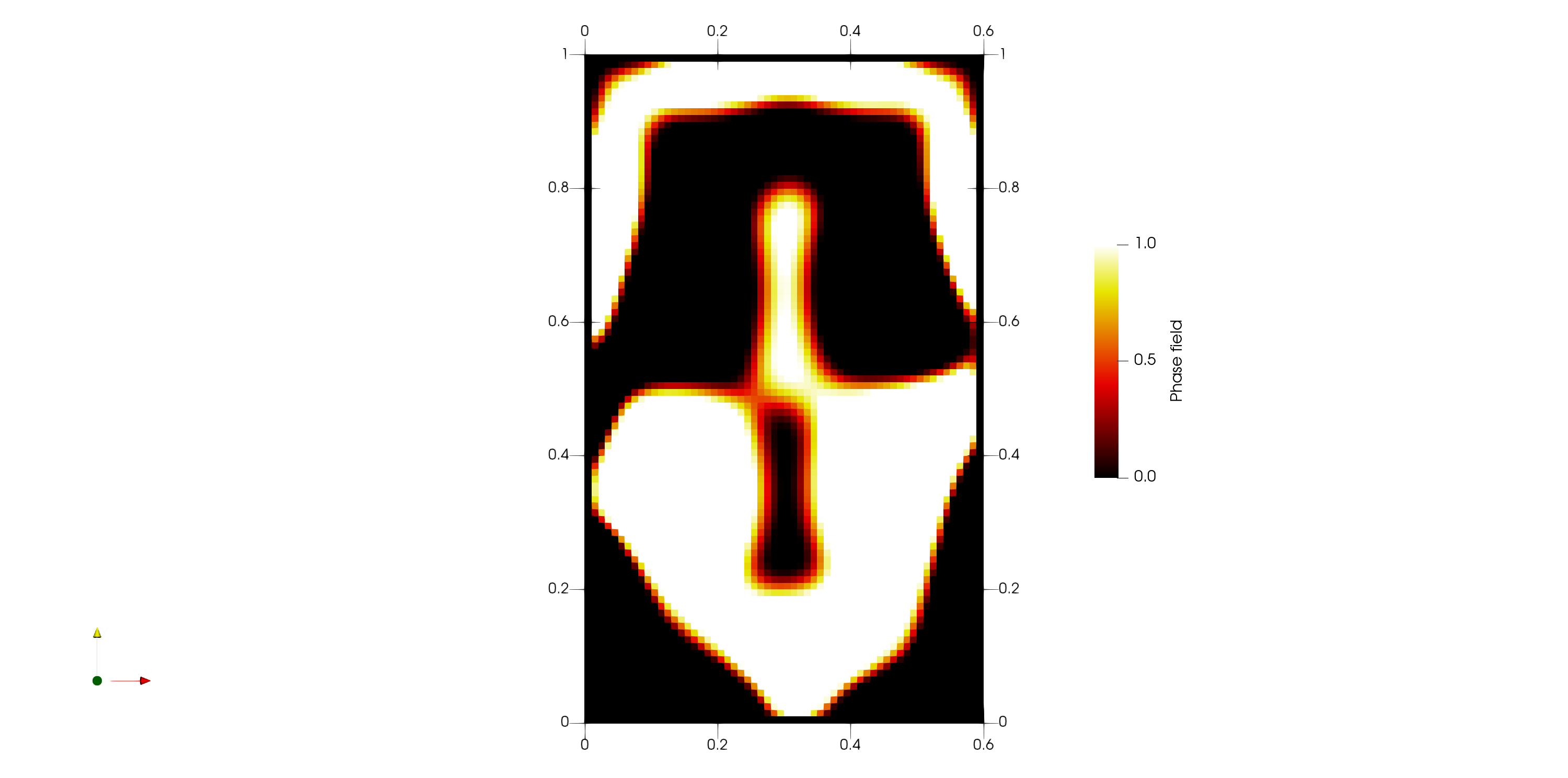}
\par\end{center}
\vspace{-0.8cm}

\begin{center}
$\tilde{y}(0.04)$
\par\end{center}%
\end{minipage}\hfill{}
\par\end{centering}
\vspace{0.3cm}

\begin{centering}
\hfill{}%
\begin{minipage}[t]{0.32\columnwidth}%
\begin{center}
\includegraphics[viewport=1000bp 0bp 2000bp 1488bp,clip,width=1\textwidth]{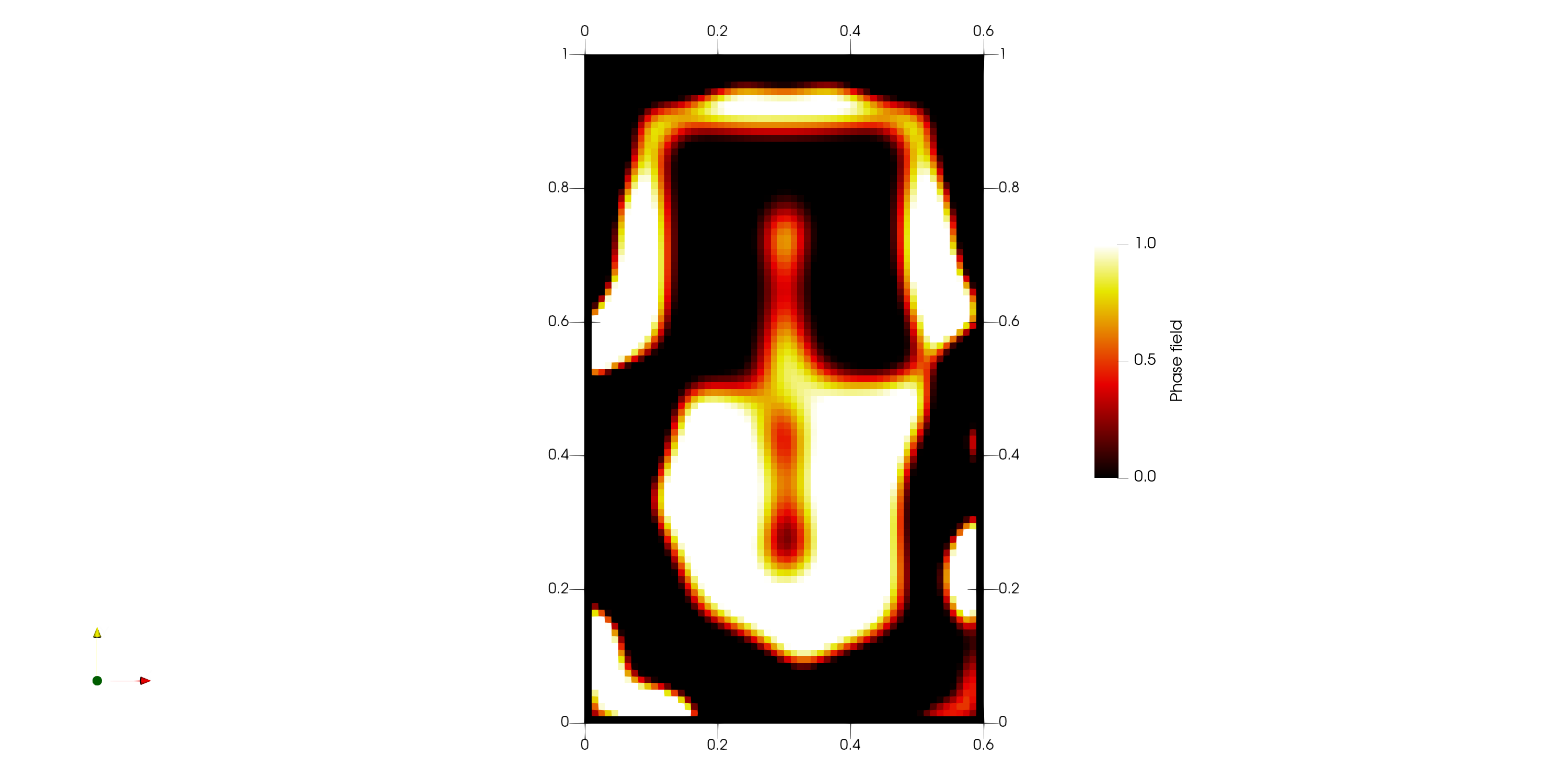}
\par\end{center}
\vspace{-0.8cm}

\begin{center}
$\tilde{y}(0.06)$
\par\end{center}%
\end{minipage}\quad{}%
\begin{minipage}[t]{0.32\columnwidth}%
\begin{center}
\includegraphics[viewport=1000bp 0bp 2000bp 1488bp,clip,width=1\textwidth]{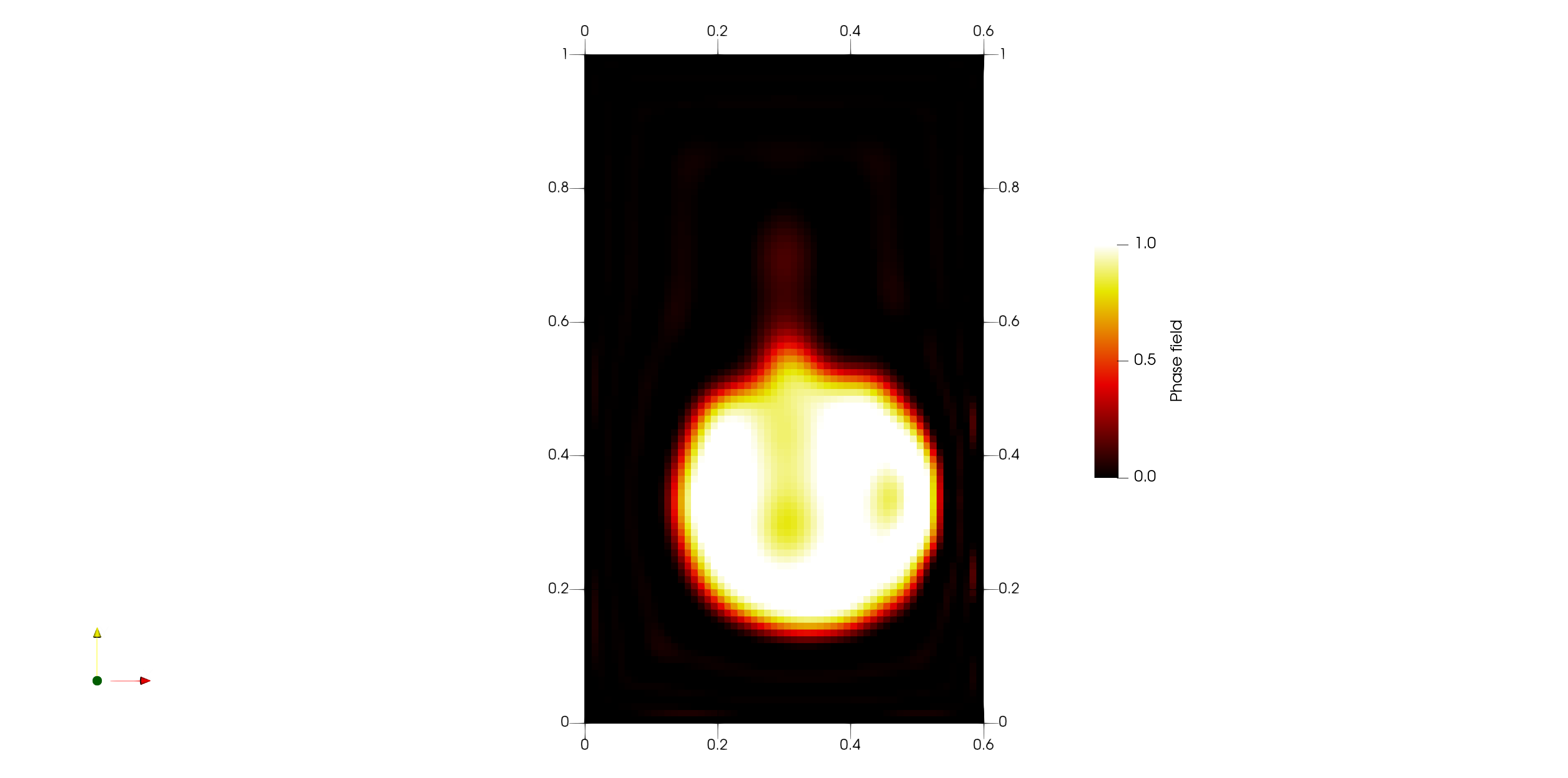}
\par\end{center}
\vspace{-0.8cm}

\begin{center}
$\tilde{y}(0.08)$
\par\end{center}%
\end{minipage}\quad{}%
\noindent\begin{minipage}[t]{0.1\columnwidth}%
\begin{center}
\includegraphics[viewport=2050bp 300bp 2280bp 1350bp,clip,width=1\textwidth]{pictures/2d_viz/pf_0}
\par\end{center}%
\end{minipage}\hfill{}
\par\end{centering}
\caption{\label{fig:test4_lin_pf}Evolution of phase field in simulation with
the linear reaction term (\ref{eq:reaction-term}). For details, see
Section \ref{subsec:Moving-a-Crystal-n-t-s-lin}.}
\end{figure}

\begin{figure}
\begin{minipage}[t]{0.47\columnwidth}%
\begin{center}
\includegraphics[viewport=800bp 0bp 1800bp 1488bp,clip,width=1\textwidth]{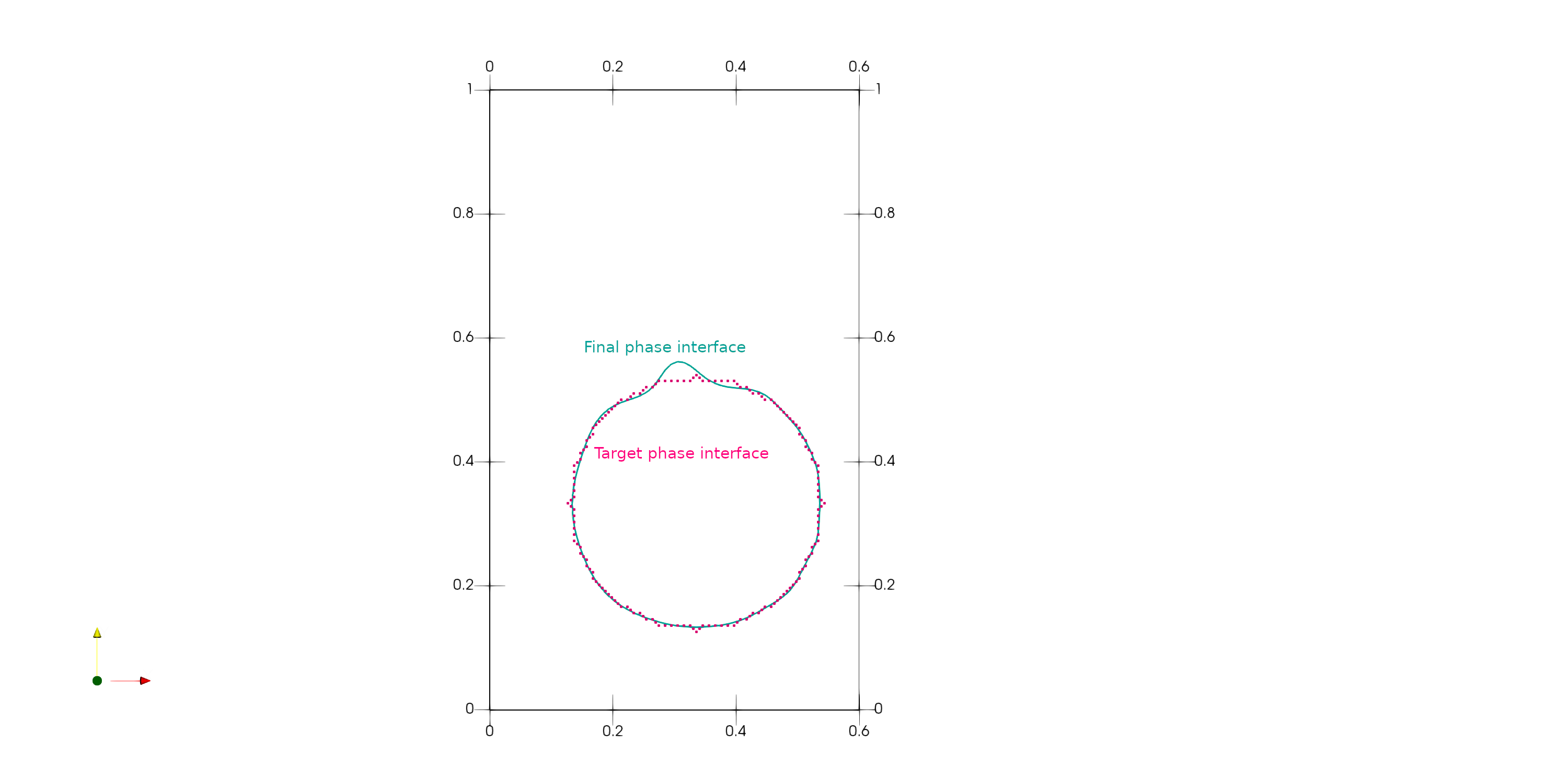}\caption{\label{fig:2d_test4_final_lin} Target and final phase field obtained
in the simulation with the linear reaction term (\ref{eq:reaction-term}).
For details, see Section \ref{subsec:Moving-a-Crystal-n-t-s-lin}.}
\par\end{center}%
\end{minipage}\hfill{}%
\begin{minipage}[t]{0.47\columnwidth}%
\begin{center}
\includegraphics[viewport=800bp 0bp 1800bp 1488bp,clip,width=1\textwidth]{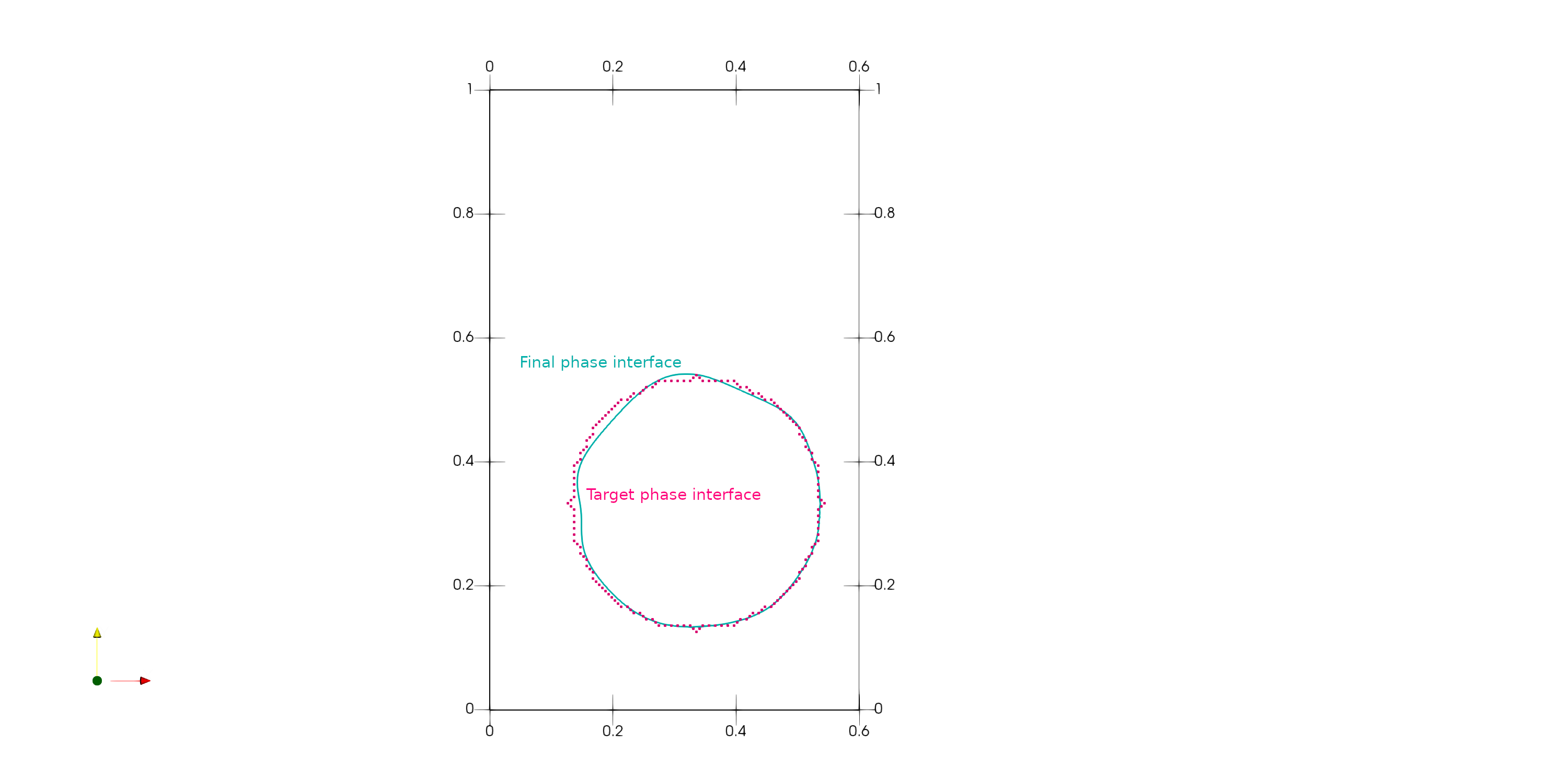}\caption{\label{fig:2d_test4_final_nonlin} Target and final phase field obtained
in the simulation with the alternative reaction term (\ref{eq:sigma limiter reaction}).
For details, see Section \ref{subsec:Moving-a-Crystal-n-t-s-lin}.}
\par\end{center}%
\end{minipage}
\end{figure}

In Figure \ref{fig:2d_test4_lin}, highly complex shape of the interface
$\Gamma$ can be observed at time $t=0.06$ and the crystal assumes
its final shape very close to the final time $t=0.08$. More snapshots
of the phase field evolution are shown in Figure \ref{fig:test4_lin_pf}.
On the other hand, in Figure \ref{fig:2d_test4_nonlin}, the crystal
keeps shape close to the original in all of the shown snapshots. In
this simulation, the evolution is affected by the non-symmetric initial
position of the crystal, as the evolution is driven mainly by the
heating and cooling of the right boundary. A preference to keep the
crystal close to the right boundary, where the Dirichlet's boundary
condition has the strongest influence, can be observed. Close to the
final time $t=0.08$, the crystal separates from the domain boundary
and fine adjustments in the control shape it to match the target profile
$\tilde{y}_{f}$. Figures \ref{fig:2d_test4_final_lin} and \ref{fig:2d_test4_final_nonlin}
show the correspondence between the target phase interface and the
final phase interface obtained as the result of the simulation with
both reaction terms and the corresponding optimal boundary control.
It may be noticed that in the case of linear reaction term, the top
portion of the final crystal does not have the optimal shape while
in the other case the crystal is pushed to the right. This difference
is caused by the different process leading to the final state.

\subsubsection{\label{subsec:Separating-a-Crystal-sigma}Separating a Crystal with
the Improved Reaction Term}

The last numerical simulation aim to address the question of crystal
separation. The more advanced reaction term (\ref{eq:sigma limiter reaction})
is used. Namely, considering a single rectangular crystal in the spatial
domain at initial time, a control that separates this crystal into
two circular ones is sought after.

\begin{figure}
\begin{minipage}[t]{0.47\columnwidth}%
\begin{center}
\includegraphics[width=0.97\textwidth]{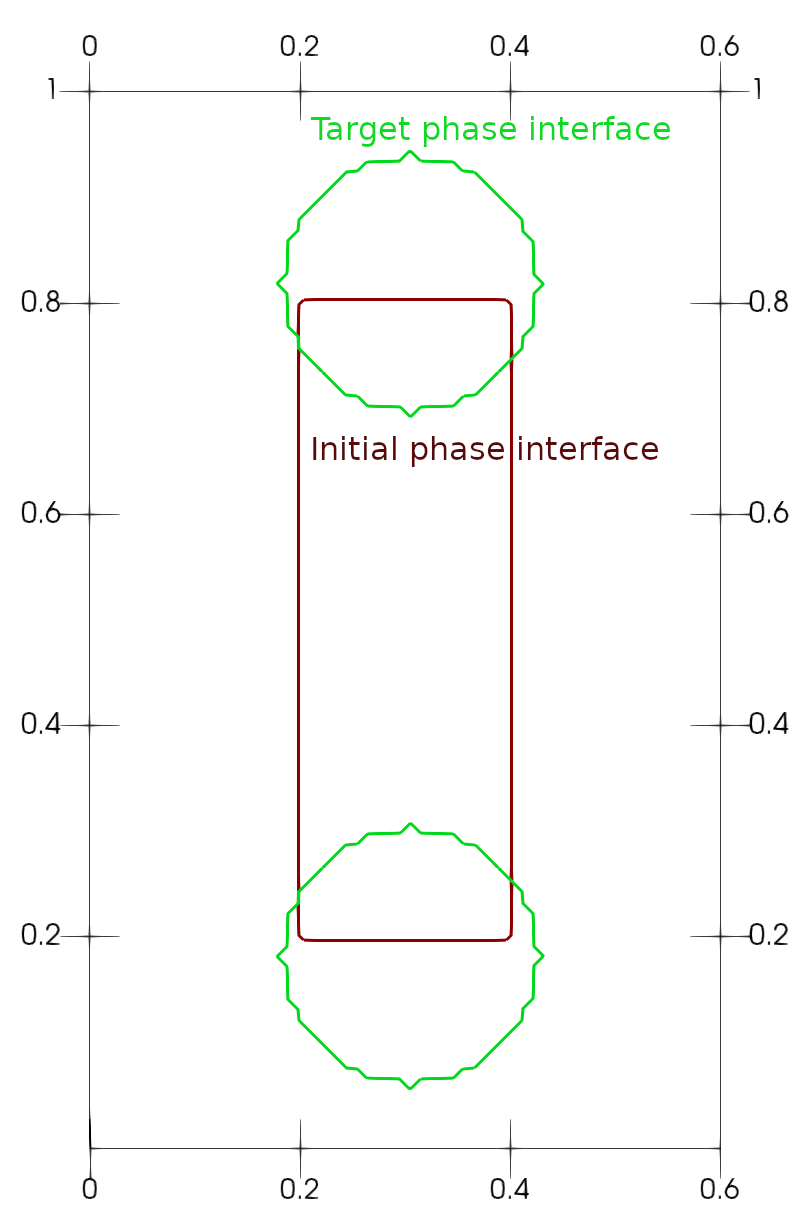}\caption{\label{fig:2d_test2_init}The initial and target position of the phase
interface $\Gamma$ for the simulation in Section \ref{subsec:Separating-a-Crystal-sigma}.}
\par\end{center}%
\end{minipage}\hfill{}%
\begin{minipage}[t]{0.47\columnwidth}%
\begin{center}
\includegraphics[width=0.97\textwidth]{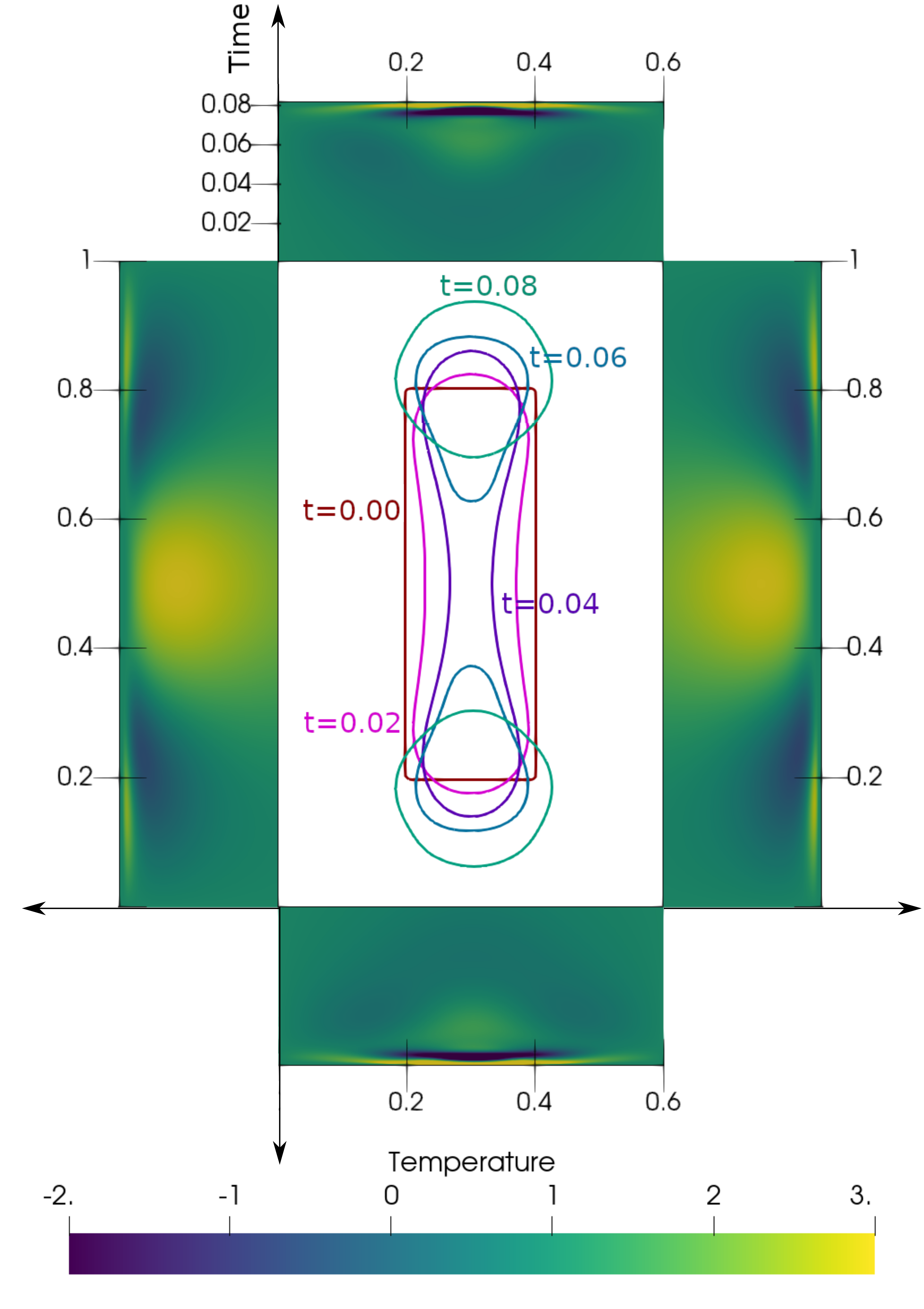}\caption{\label{fig:2d_test2_vyvoj} The heating of the domain boundary in
time and the snapshots of the phase interface $\Gamma(t)$ during
the evolution and in the final time T = 0.08. For details on this
simulation, see Section \ref{subsec:Separating-a-Crystal-sigma}.}
\par\end{center}%
\end{minipage}

\end{figure}

Figure \ref{fig:2d_test2_init} depicts the initial conditions $\tilde{y}_{\text{ini}}$
and the target profile $\tilde{y}_{f}$. The boundary condition for
the phase field $\tilde{y}$ is given by

\begin{equation}
\tilde{y}_{\text{bc}}\left(t,x\right)=0,\;\forall x\in\partial\Omega,\;\forall t\in\left[0,T\right).\label{eq:2d_test2_bc}
\end{equation}
The initial guess for the Dirichlet control on the boundary reads

\begin{equation}
u_{0}\left(t,x\right)=1,\;\forall x\in\partial\Omega,\;\forall t\in\left[0,T\right).\label{eq:2d_test2_init}
\end{equation}
Additional data for the experiment are listed in Table \ref{tab:2d_num_param}.

The time evolution of the crystal shape $\Gamma\left(t\right)$ and
the temporal control profiles are combined in Figure \ref{fig:2d_test2_init}.

\begin{figure}
\centering{}\includegraphics[viewport=800bp 0bp 1800bp 1488bp,clip,width=0.4\textwidth]{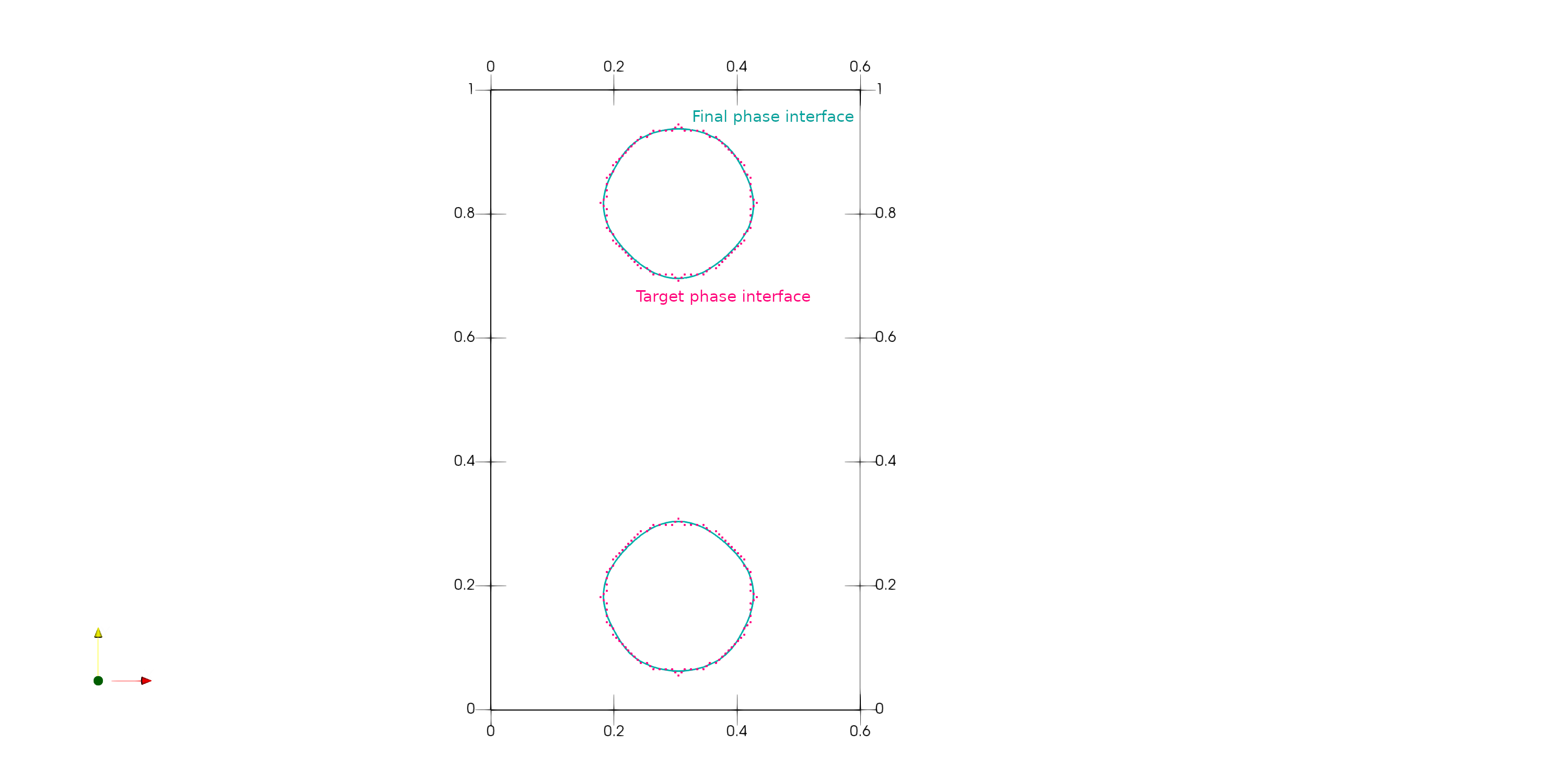}\caption{\label{fig:2d_test2_final} Target and final phase field obtained
in the simulation described in Section \ref{subsec:Separating-a-Crystal-sigma}.}
\end{figure}

On both sides of the domain, the optimal control exhibits an effort
to heat up the center of the domain in order to separate the crystal
into two. At the same time, the upper and lower parts of the domain
cool down, leading to the growth of crystals at the ends of the domain.
The two crystals end up being slightly deformed compared to the ones
prescribed by $\tilde{y}_{f}$, as can be seen in Figure \ref{fig:2d_test2_final}.

\begin{table}
\caption{\label{tab:2d_num_param}Settings for the 2D experiments and the respective
values of the difference (error) from the prescribed profile. The
first experiment was performed with two reaction terms - linear and
alternative. }

\centering{}%
\begin{tabular}{ccc}
\toprule
 & \multicolumn{2}{c}{Simulation }\tabularnewline
Parameter & Moving crystal ( Section \ref{subsec:Moving-a-Crystal-n-t-s-lin}) & Separating crystal (Section \ref{subsec:Separating-a-Crystal-sigma})\tabularnewline
\midrule
number of time steps $N_{t}$ & $8\cdot10^{3}$ & $8\cdot10^{3}$\tabularnewline
number of grid points $N_{x}$ & $60$ & $60$\tabularnewline
number of grid points $N_{y}$ & $100$ & $100$\tabularnewline
initial control given by & (\ref{eq:2d_test4_init}) & (\ref{eq:2d_test2_init})\tabularnewline
final time $T$ & $0.081$ & $0.081$\tabularnewline
regularization parameter $\alpha$ & $0$ & $0$\tabularnewline
gradient descent step size $\varepsilon$ & $5.0$ & $7.5$\tabularnewline
number of iterations & $400$ (linear), 5000 (alternative) & $1000$\tabularnewline
$\left\Vert \tilde{y}_{h}-P_{h}\tilde{y}_{f}\right\Vert _{2}$ at
$t=T$ & $7.66$ (linear), $7.39$ (alternative) & $7.21$\tabularnewline
\bottomrule
\end{tabular}
\end{table}

\subsection{Performance and Implementation Details\label{subsec:Performance-Details}}

The 1D and 2D solvers used in Section \ref{subsec:Moving-a-Crystal-n-t-s-lin}
and \ref{subsec:Separating-a-Crystal-sigma}, respectively, were implemented
in MATLAB and C++. This section intends to give the reader an idea
about the performance and not to serve as a comparison of the two
implementations. Therefore, the performance analysis is only presented
for the more computationally intensive 2D case. The C++ solver used
in Section \ref{subsec:PhaseField-Num-1} was executed on a single
CPU core. Simple benchmark can be reviewed in Table \ref{tab:Configuration-and-performance-2d}.

\begin{table}
\caption{\label{tab:Configuration-and-performance-2d}Time to compute 100 iterations
of gradient descent for different discretization resolutions. All
simulations were performed on desktop with Intel B360 AORUS MB, i7-8700
CPU and 16GB RAM (Fedora 31 Linux). The values of $N_{t}$ are given
in multiples of $10^{3}$.}

\centering{}%
\begin{tabular}{cccc}
\toprule
$N_{x}$ & $N_{y}$ & $N_{t}$ & Computational time $[\text{s}]$\tabularnewline
\midrule
51 & 51 & 1 & 37\tabularnewline
51 & 101 & 1 & 78\tabularnewline
101 & 101 & 4 & 585\tabularnewline
101 & 101 & 8 & 1188\tabularnewline
101 & 201 & 4 & 1281\tabularnewline
\bottomrule
\end{tabular}
\end{table}

\section*{Conclusion}

The formal adjoint problem for the numerical optimization of Dirichlet
boundary condition in the phase field model was derived for two different
reaction terms. The possibilities of this formulation were explored
using several simulations performed in one and two spatial dimensions
with the help of the finite difference method. Among other things,
the influence of reaction term choice, initial control estimate and
regularization are discussed. Ultimately, several experiments performed
in two spatial dimensions show how even non-trivial control that may
be interpreted as solidification can be obtained using this method.

Using the adjoint formulation makes it possible to achieve impressive
perform\-ance figures with limited hardware and a rudimentary numerical
method. This gives a positive future outlook for the application of
Dirichlet boundary control optimization for crystal morphology estimation
and solidification control in three dimensions.

\section*{Data availability}

The datasets and computer codes are available upon request from the
authors.

\section*{Declaration of Competing Interest}

The authors declare that they have no known competing financial interests
or personal relationships that could have appeared to influence the
work reported in this paper.

\section*{CRediT authorship contribution statement}

\textbf{Aleš Wodecki:} Conceptualization, Methodology, Software, Validation,
Formal analysis, Investigation, Data Curation, Writing - Original
Draft, Writing - Review \& Editing, Visualization. \textbf{Pavel Strachota:}
Methodology, Software, Investigation, Writing - Original Draft, Writing
- Review \& Editing. \textbf{Tomáš Oberhuber:} Methodology, Writing
- Review \& Editing, Supervision. \textbf{Kateřina Škardová:} Methodology,
Writing - Review \& Editing. \textbf{Monika Balázsová}: Writing -
Review \& Editing.

\begin{acknowledgements}

This work has been supported by the projects:
\begin{itemize}
\item \emph{Centre of Advanced Applied Sciences}\\
(Reg. No. CZ.02.1.01/0.0/0.0/16\_019/0000778),
\item \emph{Research Center for Informatics}\\
(Reg. No. CZ.02.1.01/0.0/0.0/16\_019/0000765),
\end{itemize}
co-financed by the European Union. Partial support of:
\begin{itemize}
\item grant No. SGS20/184/OHK4/3T/14 of the Grant Agency of the Czech Technical
University in Prague,
\item the project No. NV19-08-00071 of the Ministry of Health of the Czech
Republic.
\end{itemize}
\end{acknowledgements}




\bibliographystyle{spmpsci}      
\bibliography{literature/optimization,literature/publications,literature/references_MATH-PHYS,literature/references_MMG}

\end{document}